\documentclass[12pt,leqno]{amsart}

\usepackage{float}
\usepackage{amssymb,amsmath,amsthm,amscd}
\usepackage{graphics}
\usepackage{mathrsfs}
\usepackage{enumerate}
\usepackage[usenames,dvipsnames]{color}

\usepackage{verbatim}

\usepackage[all]{xy}

\numberwithin{equation}{section}
\newcommand{\nc}{\newcommand}

\nc{\CC}{{c}} \nc{\RC}{\mathrm{RC}}

\newcommand{\C}{{\mathbb C}}
\newcommand{\Q}{\mathbb {Q}}

\newcommand{\Z}{{\mathbb Z}}

\newcommand{\MA}{\mathcal{A}}
\newcommand{\MK}{\mathcal{K}}

\newcommand{\R}{{\rm R}}

\newcommand{\shc}{{\mathscr{C}}}

\newcommand{\seteq}{\mathbin{:=}}

\newcommand{\m}{\mathfrak{m}}

\theoremstyle{plain}
\newtheorem{lemma}{Lemma}[section]
\newtheorem{prop}[lemma]{Proposition}
\newtheorem{theorem}[lemma]{Theorem}
\newcommand{\Prop}{\begin{prop}}
\newcommand{\enprop}{\end{prop}}
\newcommand{\Lemma}{\begin{lemma}}
\newcommand{\enlemma}{\end{lemma}}
\newcommand{\Th}{\begin{theorem}}
\newcommand{\enth}{\end{theorem}}
\newtheorem{corollary}[lemma]{Corollary}
\newcommand{\Cor}{\begin{corollary}}
\newcommand{\encor}{\end{corollary}}
\newtheorem{definition}[lemma]{Definition}
\newtheorem{question}[lemma]{Question}
\newcommand{\Question}{\begin{question}}
\newcommand{\enquestion}{\end{question}}
\newcommand{\Def}{\begin{definition}}
\newcommand{\edf}{\end{definition}}

\theoremstyle{definition}
\newtheorem{remark}[lemma]{Remark}
\newtheorem{example}[lemma]{Example}

\nc{\Rem}{\begin{remark}}
\nc{\enrem}{\end{remark}}

\nc{\emprule}[1]{\rule{#1}{0pt}}
\newcommand{\g}{{\mathfrak{g}}}

\newcommand{\Hom}{\operatorname{Hom}}
\newcommand{\End}{\operatorname{End}}

\newcommand{\isoto}[1][]{\mathop{\xrightarrow[#1]%
{\rule{0pt}{.9ex}%
{\raisebox{-.6ex}[0ex][-.7ex]{$\mspace{4mu}\sim\mspace{3mu}$}}}}}
\newcommand{\isofrom}[1][]{\mathop{\xleftarrow[#1]%
{\rule{0pt}{.9ex}%
{{\raisebox{-.6ex}[0ex][-.6ex]{$\mspace{3mu}\sim\mspace{4mu}$}}}}}}
\renewcommand{\hom}{\operatorname{\it \mathscr{H}\kern-.25em om}}

\newcommand{\M}{{\mathscr M}}
\newcommand{\N}{{\mathscr N}}
\newcommand{\eq}{\begin{eqnarray}}
\newcommand{\eneq}{\end{eqnarray}}

\newcommand{\eqn}{\begin{eqnarray*}}
\newcommand{\eneqn}{\end{eqnarray*}}
\newcommand{\op}{{\operatorname{op}}}
\newenvironment{tenumerate}{
  \begin{enumerate}
  
  }{\end{enumerate}}

\nc{\bnum}{\begin{enumerate}[{\rm (i)}]}
\nc{\enum}{\ee}
\nc{\bna}{\begin{enumerate}[{\rm (a)}]}
\nc{\ena}{\ee}

\newcommand{\on}{\operatorname}
\newcommand{\Ker}{\ker}

\newcommand{\bni}{\begin{tenumerate}}
\newcommand{\eni}{\end{tenumerate}}

\newcommand{\QED}{\end{proof}}
\newcommand{\Proof}{\begin{proof}}

\newcommand{\soplus}{\mathop{\mbox{\normalsize$\bigoplus$}}\limits}

\newcommand{\cl}{\colon}
\newcommand{\To}[1][{\emprule{1.8ex}}]{\xrightarrow{\,#1\,}}
\newcommand{\from}[1][{\emprule{1.8ex}}]{\xleftarrow{\,#1\,}}
\newcommand{\id}{\on{id}}
\newcommand{\ba}{\begin{array}}
\newcommand{\ea}{\end{array}}

\newcommand{\monoto}{\rightarrowtail}
\newcommand{\monofrom}{\leftarrowtail}

\newcommand{\set}[2]{\left\{#1 \mathbin; #2 \right\}}
\newcommand{\supp}{\operatorname{supp}}
\newcommand{\Mod}{\operatorname{Mod}}
\newcommand{\hs}{\hspace*}
\newcommand{\vs}{\vspace*}

\newcommand{\eqsub}{\begin{subequations}\begin{eqnarray}}
\newcommand{\eneqsub}{\end{eqnarray}\end{subequations}}

\newcommand{\ol}{\overline}

\newcommand{\A}{\mathscr{A}}

\nc{\la}{\lambda}
\nc{\lam}{\lambda}
\nc{\U}[1][\g]{U_q(#1)}
\nc{\te}{\tilde{e}}
\nc{\tei}{\tilde{e}_i}
\nc{\tf}{\tilde{f}}
\nc{\tfi}{\tilde{f}_i}
\nc{\tU}{\widetilde U_q(\g)}
\nc{\tE}{\tilde{E}}
\nc{\tF}{\tilde{F}}

\nc{\BZ}{{\mathbb{Z}}}
\nc{\al}{\alpha}
\nc{\qs}{{q}}
\nc{\lan}{\langle}
\nc{\ran}{\rangle}
\nc{\re}{{\mathrm{re}}}
\nc{\wt}{\operatorname{wt}}
\nc{\Uf}[1][\g]{U^-_q(#1)}
\nc{\Ue}{U^+_q(\g)}
\nc{\eps}{\varepsilon}
\nc{\vphi}{\varphi}
\nc{\sphi}{\varphi^*}
\nc{\seps}{\varepsilon^*}

\nc{\nn}{\nonumber}

\nc{\vp}{\varpi}
\nc{\cls}{{\operatorname{cl}}}

\nc{\Wt}{{\operatorname{Wt}}}
\nc{\Us}{U'_q(\g)}
\nc{\La}{\Lambda}
\nc{\ro}{{\rm(}}
\nc{\rf}{{\rm)}}
\nc{\norm}{{\mathrm{norm}}}
\nc{\qbox}{\quad\mbox}
\nc{\braid}{{\mathfrak{B}}}
\nc{\Ad}{\operatorname{Ad}}
\nc{\Aut}{\operatorname{Aut}}
\nc{\dt}[1]{\tilde{\tilde #1}}
\nc{\Sn}{S^{{\mathrm{norm}}}}
\nc{\aff}{{\mathrm{aff}}}
\nc{\rk}{{\mathrm{rk}}}
\nc{\tQ}{\widetilde{Q}}
\nc{\tP}{\widetilde{P}}
\nc{\tW}{\widetilde{\mathscr{W}}}
\nc{\Dyn}{\mathrm{Dyn}}
\nc{\tD}{\widetilde{\Delta}}
\nc{\height}{{\operatorname{ht}}}
\nc{\bl}{\bigl}
\nc{\br}{\bigr}
\nc{\Hecke}{\mathrm{H}}
\nc{\HA}{\Hecke^{\mathrm{A}}}
\nc{\HB}{\Hecke^{\mathrm{B}}}
\nc{\K}{\mathrm{K}}
\newcommand{\scbul}{{\,\raise1pt\hbox{$\scriptscriptstyle\bullet$}\,}}
\nc{\vac}{{\phi}}
\nc{\be}{\begin{enumerate}}
\nc{\ee}{\end{enumerate}}
\nc{\low}{{\mathrm{low}}}
\nc{\upper}{{\mathrm{up}}}
\nc{\Zodd}{\Z_{\mathrm{odd}}}
\nc{\Ft}[1][n]{\mathbb{P}\mathrm{ol}_{#1}}
\nc{\Ftf}[1][n]{\widetilde{\mathbb{P}\mathrm{ol}}_{#1}}
\nc{\KA}{\on{K}^{\mathrm{A}}}
\nc{\KB}{\on{K}^{\mathrm{B}}}
\nc{\Res}{\on{Res}}
\nc{\Fc}[1][{n,m}]{\mathbf{F}_{#1}}
\nc{\tphi}{\widetilde{\varphi}}
\nc{\CO}{\mathscr{O}}
\nc{\CK}{\mathscr{K}}
\nc{\disc}{\mathfrak{d}}
\nc{\tr}{\on{tr}}
\nc{\Gb}{\mathfrak{b}}
\nc{\ga}{\mathfrak{a}}
\nc{\stable}{\mathrm{stable}}
\nc{\X}{\mathfrak{X}}
\nc{\Hilb}{\mathrm{Hilb}}
\nc{\W}{\ensuremath{\mathscr{W}}}
\nc{\Ws}{\ensuremath{\rm W}}
\nc{\opp}{{\on{opp}}}
\newcommand{\prolim}[1][]{\mathop{\varprojlim}\limits_{#1}}
\nc{\corps}{{\mathbf{k}}}
\nc{\cor}{{\mathbf{k}}}
\nc{\h}{\mathrm{\hslash}}
\nc{\fL}[1][{\h}]{\C(\mspace{-1mu}(#1)\mspace{-1mu})}
\nc{\ad}{\mathrm{ad}}
\newcommand{\Endm}{\operatorname{\mathscr{E}\kern-.1pc\mathit{nd}}}
\newcommand{\Endomo}{\operatorname{\mathscr{E}\kern-.1pc\mathit{nd}}}
\nc{\bc}{\bar{\corps}}
\nc{\reg}{{\mathrm{reg}}}
\nc{\ysq}{\mathbf{y}^2}
\nc{\Ch}{\on{Ch}}

\nc{\sketch}{\Proof}
\nc{\Gm}{\mathbb{G}_{\mathrm{m}}}
\nc{\hGm}{\hat{\mathbb{G}}_{\mathrm{m}}}
\nc{\ug}{\widehat{\mathrm{G}}_{\mathrm{m}}}

\nc{\tL}{\widetilde{\mathscr{L}}}
\nc{\Fr}{\mathcal{F}}
\nc{\E}{\mathcal{E}}

\nc{\ord}{\on{ord}}
\nc{\bM}{\overset{\hs{1.5ex}\rule[-.08ex]{1.8ex}{.08ex}}{\M}}
\nc{\romano}{\mathrm{o}}
\nc{\into}{\hookrightarrow}
\nc{\good}{\mathrm{good}}
\nc{\tA}{\widetilde\A}
\nc{\Vz}{{V}\kern-1.1ex\raisebox{1.5ex}[0ex][0ex]{$\cdot$}}
\nc{\bxes}[1]{\raisebox{.9ex}{$\cdot$}{\kern#1}\raisebox{0ex}{$\cdot$}
{\kern#1}\raisebox{-.9ex}{$\cdot$}}
\nc{\ssum}{\mathop{\mbox{\normalsize{${\sum}$}}}\limits}
\nc{\ct}{{\mbox{\tiny$\mathrm{CT}$}}}
\nc{\pr}{\mathrm{pr}}
\nc{\qr}{\mathrm{rp}}

\nc{\Fs}{\ensuremath{\rm F}}

\nc{\isotf}{\overset{
{\rule{0pt}{.9ex}%
{\raisebox{-.6ex}[0ex][-.7ex]{$\mspace{3mu}\sim\mspace{3mu}$}}}}
{\longleftrightarrow}}
\nc{\tN}{\tilde\N}
\nc{\tens}{\mathop\otimes\limits}
\nc{\super}{\mathrm{super}}
\nc{\Mods}{\on{Mod}_\super}
\nc{\rev}{\mathrm{rev}}
\nc{\clif}{\mathfrak{C}}
\nc{\clifm}{\clif^{-}}
\nc{\Fct}{\mathrm{Fct}}
\nc{\Fcts}{\mathrm{Fct}_\super}

\nc{\Ks}{\on{K^\super}}
\nc{\ts}{\widetilde{s}}
\nc{\KUGIRI}{\circ}
\nc{\Sym}{\mathfrak{S}}
\nc{\FF}{\Z/2\Z}

\nc{\cc}{\mathfrak{c}}
\nc{\SK}{\mathcal{KS}}
\nc{\noi}{\noindent}
\nc{\odd}{{\mathrm{odd}}}
\nc{\even}{{\mathrm{even}}}
\nc{\bs}{\ol{s}}
\nc{\Khc}[1][n]{\ol{\mathcal{KHC}}_{#1}}
\nc{\Ohc}[1][n]{\ol{\mathcal{OHC}}_{#1}}
\nc{\KHC}[1][n]{\mathcal{K}{\mathcal{HC}}_{#1}}
\nc{\OHC}[1][n]{\mathcal{O}{\mathcal{HC}}_{#1}}

\nc{\IODD}{I_{\odd}}
\nc{\IEVEN}{I_{\even}}
\DeclareMathOperator{\KKK}{\mathsf{K}_0}
\newcommand{\MOD}[1]{{#1\text{-}\mathsf{mod}}}
\newcommand{\F}{\mathbb{F}}
\nc{\IRED}{I}
\DeclareMathOperator{\CHAR}{char}
\newcommand{\QTL}{\tQ}

\newcommand{\MSP}[1]{\cor\mathfrak{S}_{#1}^{-}}

\nc{\MH}{\mathcal{H}}

\newcommand{\KLR}[1]{\mathsf{KLR}_{#1}}
\DeclareMathOperator{\Irr}{\mathsf{Irr}}
\newcommand{\SMOD}[1]{{#1\text{-}\mathsf{smod}}}
\newcommand{\AHCC}{{\mathcal{AHC}}}
\newcommand{\KLRR}[2]{\mathsf{KLR}^{#1}_{#2}}
\newcommand{\MAPSTO}{\longmapsto}

\newcommand{\pl}{\dag}

\newcommand{\RCP}{\RC^{\pl}}
\newcommand{\EP}{e^{\pl}}

\DeclareMathOperator{\DEG}{\mathsf{deg}}
\newcommand{\HAT}[1]{\widehat{#1}}
\newcommand{\AAA}{\mathbb{A}}
\DeclareMathOperator{\FRAC}{\mathsf{Frac}}
\newcommand{\CMH}[1]{{\overline{\mathcal{AHC}}}^{\Lambda_0}_{#1}}
\newcommand{\RP}[1]{\mathsf{RP}_3}
\newcommand{\dAHC}[1]{\overline{\mathcal{AHC}}_{#1}}
\newcommand{\dAHCC}{{\overline{\mathcal{AHC}}}}
\newcommand{\SV}[1]{V^{\mathsf{spin}}_{#1}}
\newcommand{\HV}[1]{V^{\mathsf{KLR}}_{#1}}
\newcommand{\AHC}[1][n]{{\mathcal{AHC}_{#1}}}
\newcommand{\SA}{{\mathsf{SA}}}
\newcommand{\NSA}{{\mathsf{NSA}}}
\nc{\bwr}{\mbox{\large$\wr$}}
\nc{\At}[1][{{i,j}}]{\mathscr{A}_{{#1}}}
\nc{\tAt}[1][{{i,j}}]{\widetilde{\mathscr{A}}_{{#1}}}
\nc{\Bt}[1][{{i,j}}]{\mathscr{B}_{{#1}}}
\nc{\tBt}[1][{{i,j}}]{\widetilde{\mathscr{B}}_{{#1}}}
\nc{\prt}[1]{\mathrm{par}(#1)}
\nc{\erm}{\mathrm{e}}
\nc{\ec}{\mathrm{e}^-}

\nc{\GCM}{{generalized Cartan matrix}}
\nc{\HCO}{\widehat{\CO}}
\nc{\tCO}{\widetilde{\CO}}
\nc{\T}{\mathbb{T}}
\nc{\HC}[1][n]{\mathcal{HC}_{{#1}}}
\nc{\HRC}[1][n]{\widehat{\mathrm{RC}}_{{#1}}}
\nc{\hc}[1][n]{\ol{\mathcal{HC}}_{{#1}}}
\nc{\bphi}{\bar{\phi}}
\nc{\hgt}{\mathrm{ht}}
\nc{\double}{\mathrm{double}}
\nc{\Ob}{\on{Ob}}

\begin{document}

\title
{Quiver Hecke superalgebras}

\author{Seok-Jin Kang}
\address{S-J.K.: Department of Mathematical Sciences and
Research Institute of Mathematics,
Seoul National University,
599 Gwanak-ro, Gwanak-gu, Seoul 151-747, Korea}
\email{sjkang@snu.ac.kr}
\thanks{The research of S-J.K.\
was supported by KRF Grant \# 2007-341-C00001 and NRF Grant \#
2010-0019516.}

\author{Masaki KASHIWARA}
\address{M.K.: Research Institute for Mathematical Sciences,
Kyoto University, Kyoto 606, Japan and Department of Mathematical
Sciences, Seoul National University, 599 Gwanak-ro, Gwanak-gu, Seoul
151-747, Korea }
\thanks{The research of M.K.\  was partially supported by
Grant-in-Aid for Scientific Research (B) 22340005
Japan Society for the Promotion of Science.}
\email{masaki@kurims.kyoto-u.ac.jp}

\author{Shunsuke Tsuchioka}
\address{S.T.: Institute for the Physics and
Mathematics of the Universe, University of Tokyo,
Kashiwano-ha 5-1-5, Kashiwa City, Chiba 277-8582, Japan}
\thanks{The research of S.T.\  was supported by
Grant-in-Aid for Research Activity Startup 22840026 and
Research Fellowships for Young Scientists 23$\cdot$8363,
Japan Society for the Promotion of Science.}
\email{tshun@kurims.kyoto-u.ac.jp}


\begin{abstract}

We introduce a new family of superalgebras which should be
considered as a super version of the Khovanov-Lauda-Rouquier
algebras.
Let $I$ be the set of vertices of a Dynkin diagram with a
decomposition $I=\IEVEN\sqcup\IODD$.
To this data,  we associate a
family of graded superalgebras $\R_n$, the {\em quiver Hecke
superalgebras}.  When $\IODD=\emptyset$, these algebras are nothing
but the usual Khovanov-Lauda-Rouquier algebras. We then define
another family of graded superalgebras $\RC_n$, the quiver
Hecke-Clifford superalgebras, and show that the superalgebras $\R_n$
and $\RC_n$ are weakly Morita superequivalent to each other.
Moreover, we
prove that the affine Hecke-Clifford superalgebras, as well as their
degenerate version, the affine Sergeev superalgebras, are isomorphic
to quiver Hecke-Clifford superalgebras $\RC_n$ after a completion.
\end{abstract}

\maketitle


\section{Introduction}\label{Intro}

In \cite{KL1, KL2, Rou1}, Khovanov-Lauda and Rouquier independently
introduced a remarkable family of graded algebras, the {\it
Khovanov-Lauda-Rouquier algebras} or the {\it quiver Hecke
algebras}, that categorifies the negative half of quantum groups
associated with symmetrizable Kac-Moody algebras. Moreover, the
cyclotomic quotients of Khovanov-Lauda-Rouquier algebras provide a
categorification of integrable highest weight modules over quantum
groups \cite{KL1, KK}. An important application of the
Khovanov-Lauda-Rouquier algebras is that one can derive a
homogeneous presentation of the symmetric group algebras, which
gives a quantization of Ariki's categorification theorem for the
basic $\widehat{\mathfrak{sl}}_p$-module (\cite{Ari, LLT, BK1, BK2,
Rou1}): $$V(\Lambda_0)\cong \bigoplus_{n\geq
0}\KKK(\MOD{\F_p\Sym_{n}})_{\C}.$$ Considering the long history of
the theory of symmetric groups, it is quite surprising that one can
define such a non-trivial grading on the symmetric group algebras,
which has been conjectured to exist for some time (\cite{Rou2, Tur}), only after the
discovery of Khovanov-Lauda-Rouquier algebras.
For the Dynkin diagrams of affine type $A$ or of type $A_{\infty}$,
the Khovanov-Lauda-Rouquier algebras (resp.\ the cyclotomic
Khovanov-Lauda-Rouquier algebras) are isomorphic to the affine Hecke
algebras of type $A$ (resp.\ the cyclotomic Hecke algebras or the
Ariki-Koike algebras) after we take specializations
and localizations. Thus the Khovanov-Lauda-Rouquier algebras can be
regarded as a generalization of the affine Hecke algebras of type
$A$. Since the
Khovanov-Lauda-Rouquier algebras are graded, they fit more naturally
into the categorification of quantum groups.

Indeed, in \cite{VV1}, Varagnolo and Vasserot showed that the
Khovanov-Lauda-Rouquier algebras associated with symmetric
generalized Cartan matrices are isomorphic to the Yoneda algebras of
certain complexes of constructible sheaves, and proved that the sets
of isomorphism classes of 
irreducible modules over the
Khovanov-Lauda-Rouquier algebras correspond to Kashiwara's upper
global bases\;(=Lusztig's dual canonical bases). In the same spirit, in
\cite{VV1}, Varagnolo and Vasserot (resp.\ in \cite{SVV}, Shan,
Varagnolo and Vasserot) introduced a version of
Khovanov-Lauda-Rouquier algebras corresponding to the affine Hecke
algebras of type $B$ and $C$ (resp.\ type $D$) and proved the
Lascoux-Leclerc-Thibon-Ariki type conjecture formulated in \cite{EK}
(resp.\ in \cite{KM}).

\medskip
The purpose of this paper is to introduce a new family of graded
superalgebras which should be considered as a super version of
the Khovanov-Lauda-Rouquier algebras.
Let $\cor$ be a commutative ring with 1, and
let $\bl(a_{i,j}\br)_{i,j\in I}$ be a
symmetrizable generalized Cartan matrix indexed by
a set $I$ with a decomposition
$I=\IEVEN\sqcup\IODD$ of $I$ into the set $\IEVEN$
of even vertices and the set $\IODD$ of odd vertices
such that $a_{i,j}\in 2\mspace{1mu}\Z$
for any $i\in \IODD$.
To such data we associate (see
\S\;\ref{specRC}) a family of (skew-)polynomials
$Q=(Q_{i,j}(u,v))_{i,j \in I}$ satisfying the conditions given in
\eqref{eq:Q}, and then we define a family of graded
$\cor$-superalgebras $\R_n$, the {\em quiver Hecke superalgebras}
(Definition~\ref{DefRKLR}). When $\IODD=\emptyset$, these algebras
are nothing but the usual Khovanov-Lauda-Rouquier algebras.

Let $J$ be an index set with an involution $c\cl J \to J$ and let
$\tQ=(\tQ_{i,j}(u,v))_{i,j\in J}$ be a family of polynomials
satisfying the conditions given in \eqref{cond:Q}. We define another
family of graded $\cor$-superalgebras $\RC_n$, {\it the quiver
Hecke-Clifford superalgebras} associated with $\tQ$. Let $\sim_c$ be
the equivalence relation defined by $j\sim_c j'\Leftrightarrow
\text{$j=j'$ or $j'=c(j)$}$, $I$ the set of equivalence classes, and
$I_\odd\subset I$ the image of $J^c\seteq\set{j\in J}{c(j)=j}$. Then
we show that the corresponding superalgebras $\R_n$ and $\RC_n$ are
{\it weakly Morita superequivalent} to each other (Theorem
\ref{MEQRR}). Moreover, we prove that, after a completion, the
affine Hecke-Clifford superalgebras are isomorphic to quiver
Hecke-Clifford superalgebras $\RC_n$ with a suitable choice of $\tQ$
(Theorem~\ref{th:isom_aDHC}). We  also show that the affine Sergeev
superalgebras, a degenerate version of the affine Hecke-Clifford
superalgebras, are isomorphic to quiver Hecke-Clifford superalgebras
(Theorem \ref{th:isom_aDhc}). Consequently, the affine
Hecke-Clifford superalgebras as well as the affine Sergeev
superalgebras are weakly Morita superequivalent to the quiver Hecke
superalgebras after a completion. 
We would like to add a remark
that, in \cite{Wang07}, Wang introduced the notion of {\em spin
affine Hecke algebras} (resp.\ {\em cyclotomic spin Hecke algebras})
and showed that they are {\em Morita superequivalent} to affine
Hecke-Clifford superalgebras (resp.\ cyclotomic Hecke-Clifford
superalgebras). We also note that the idea of using skew-polynomials
in super version of affine Hecke algebras has already appeared in
\cite{Wang09}. 

Let $q$ be the defining parameter for the affine Hecke-Clifford
superalgebras. In \cite{BK3}, Brundan-Kleshchev showed that when $q^2$ is a
primitive $(2\,l +1)$-th root of unity,
the representation theory of some blocks of the affine Hecke-Clifford superalgebras as
well as that of  their cyclotomic quotients is controlled by the
representation theory of quantum affine algebras of type
$A_{2l}^{(2)}$ at the crystal level.
Later, in \cite{Tsu}, the third author showed that
when $q^2$ is a primitive $(2l)$-th root of unity,
the representation
theory of some blocks is controlled by the
representation theory of quantum affine algebras of type
$D_{l}^{(2)}$.
In
other words, the affine Hecke-Clifford superalgebras and their
cyclotomic quotients with the above defining parameters give a
crystal version of the categorification theorem for the quantum
affine algebras of type $A_{2l}^{(2)}$, $D_{l}^{(2)}$ and their
highest weight modules.
In this paper, by taking all the finite-dimensional representations
of affine Hecke-Clifford superalgebras, we prove that the
corresponding algebra are isomorphic to the quiver Hecke-Clifford
superalgebras associated with (affine) Dynkin diagrams of type
$A_{\infty}$, $B_{\infty}$, $C_{\infty}$, $A_{l}^{(1)}$,
$A_{2l}^{(2)}$, $C_{l}^{(1)}$ and $D_{l}^{(2)}$.

In the forthcoming papers, we plan to prove the
conjecture that the quiver Hecke superalgebras and their cyclotomic
quotients provide a categorification of the quantum Kac-Moody
algebras and their highest weight
modules.

It is also worthwhile to note that any irreducible supermodule over a
quiver Hecke superalgebra seems to remain irreducible after forgetting its
super structure, contrary to
affine Hecke-Clifford superalgebras or quiver Hecke-Clifford superalgebras.
Hence
quiver Hecke superalgebras are more suitable for a categorification of
quantum groups
than quiver Hecke-Clifford superalgebras.

\medskip
Finally, we would like to emphasize that the superalgebras $\R_n$ and
$\RC_n$ {\it cannot} be reduced to the usual Khovanov-Lauda-Rouquier
algebras. In the following two cases,
both (cyclotomic) Khovanov-Lauda-Rouquier algebras
and (cyclotomic) quiver Hecke superalgebras should categorify the same
weight space of the representation of $U_v(\g)$, but
they are neither Morita equivalent nor weakly Morita superequivalent (see
\S\;\ref{weak superequivalence}).

\begin{enumerate}
\item (see Remark \ref{comp_calc1})
Consider the case in which $\g=A_2^{(2)}$ and $\CHAR \cor =3$ with
$I=\{0,1\}$, $I_\even=\{1\}$, $I_\odd=\{0\}$, where $0$ is a short
root.
Let $\R_\beta$ (see \eqref{eq:Rbeta})
be a direct summand of $\R_{11}$, which categorifies
$U^{-}_v(A_{2}^{(2)})_{-\beta}$ with $\beta = 8 \alpha_0 + 3 \alpha_1$.
Although the Khovanov-Lauda-Rouquier algebra
$\KLR{\beta}(A_{2}^{(2)})$ also categorifies
$U_{v}^{-}(A_{2}^{(2)})_{-\beta}$, $\Irr(\MOD{\KLR{\beta}(A_2^{(2)})})$
and $\Irr(\MOD{\R_\beta})$ correspond to different perfect bases
in the sense of \cite{BeKa} at the specialization $v=1$.

\item (see Remark~\ref{comp_calc2})
When $\g=A_1^{(1)}$, $\CHAR \cor =0$
with $I=I_\odd=\{0,1\}$,
there is no Morita equivalence nor weak Morita superequivalence
between the cyclotomic quiver Hecke superalgebra $\R_4^{\Lambda_0}$
whatever superalgebra structure we give to $\KLRR{\Lambda_0}{4}(A)$ and
whatever defining parameter of $\KLR{4}(A)$ we take.
\end{enumerate}

\vskip 3mm

This paper is organized as follows. In Section 2, we introduce the
theory of {\it supercategories} as a basic language, including the
notion of {\it Clifford twist} and {\em weak Morita superequivalence}. In
Section 3, we define the quiver Hecke superalgebras $\R_n$ and
 the
quiver Hecke-Clifford superalgebras $\RC_n$, and show that the
superalgebras $\R_n$ and $\RC_n$ are weakly Morita superequivalent to each
other. We also provide PBW-type bases for both superalgebras.
In Section 4, we
prove that, the affine Hecke-Clifford superalgebras are isomorphic
to quiver Hecke-Clifford superalgebras after a completion. The {\it
intertwiners} play an important role in the proof. In Section 5, by
a similar argument, we show that the affine Sergeev algebras are
isomorphic to quiver Hecke-Clifford superalgebras.

\vskip 5mm

{\bf Acknowledgments.} The authors would like to thank Lukas Maas
and J\"urgen M\"uller for telling us the correct dimension of the
irreducible spin representations of $\Sym_{11}$ in characteristic 3
parametrized by a partition $\lambda=(6,4,1)$ in GAP Forum. The
dimension, as explained in Remark \ref{comp_calc1}, convinced us
that $\R_n$ is a genuine generalization of the Khovanov-Lauda-Rouquier
algebras while
both of them should categorify the same Lie-theoretic objects.

Some part of this work was done when S-J.\,K.\ visited RIMS,
and when M.\,K.\ and S.\,T.\ were visiting the
Hausdorff Institute, Bonn, to participate in the trimester
``On the interactions of Representation theory with Geometry
and Combinatorics''. S-J.\,K.\ thank RIMS and M.\,K.\ and S.\,T.\
thank the organizers and HIM
for their hospitality.

\section{Supercategories}

\subsection{Clifford twist}
Let $\cor$ be a commutative ring. Let us recall that a {\em
$\cor$-linear} category is a category $\shc$ such that
$\Hom_\shc(X,Y)$ is endowed with a $\cor$-module structure for all
$X,Y\in\shc$ and the composition map
$\Hom_\shc(Y,Z)\times\Hom_\shc(X,Y)\to\Hom_\shc(X,Z)$ is
$\cor$-bilinear. A functor $F\cl\shc\to\shc'$ from a $\cor$-linear
category $\shc$ to a $\cor$-linear category $\shc'$ is called
$\cor$-linear if $F\cl \Hom_{\shc}(X,Y)\to\Hom_{\shc'}(F X,F Y)$ is
$\cor$-linear for any objects $X,Y\in\shc$. We say that an additive
 category is {\em idempotent-complete} if $\Ker(p)$
exists for any object $X$ of $\shc$ and $p\in\End_\shc(X)$ such that
$p^2=p$. Hence for any such $p$, we have $X\simeq \Ker(p)\oplus
\Ker(1-p)$.

\Def \hfill
\bnum
\item
A {\em supercategory} is a category $\shc$ equipped with an
endofunctor $\Pi$ of $\shc$ and an isomorphism
$\xi\cl\Pi^2\isoto\id_{\shc}$ such that
$\xi\circ\Pi=\Pi\circ\xi \in
\Hom(\Pi^3,\Pi)$.
\item
For a pair of supercategories $(\shc,\Pi,\xi)$ and
$(\shc',\Pi',\xi')$, a {\em superfunctor} from $(\shc,\Pi,\xi)$ to
$(\shc',\Pi',\xi')$ is a pair $(F,\alpha_F)$ of a functor $F\cl
\shc\to\shc'$ and an isomorphism $\alpha_F\cl F\circ \Pi\isoto
\Pi'\circ F$ such that the diagram \eq&&\xymatrix@C=7ex{
F\circ\Pi^2\ar[r]^-{\alpha_F\circ\Pi}\ar[d]^-{F\circ \xi}&\Pi'\circ
F\circ\Pi
\ar[r]^-{\Pi'\circ\alpha_F}&\Pi'{}^2\circ F\ar[d]_{\xi'\circ F}\\
F\ar[rr]^{\id_F}&&F}
\label{diag:comct}
\eneq
commutes.
If $F$ is an equivalence of categories,
we say that $(F,\alpha_F)$ is an {\em equivalence of supercategories}.
\item Let $(F,\alpha_F)$ and $(F',\alpha_{F'})$ be superfunctors
from a supercategory  $(\shc,\Pi,\xi)$ to $(\shc',\Pi',\xi')$.
A morphism from $(F,\alpha_F)$ to $(F',\alpha_{F'})$ is a
morphism of functors $\vphi\cl F\to F'$ such that
$$\xymatrix@C=8ex@R=4ex{
F\circ\Pi\ar[r]^{\vphi\circ\Pi}\ar[d]_{\alpha_F}&F'\circ\Pi\ar[d]_{\alpha_{F'}}\\
\Pi'\circ F\ar[r]^{\Pi'\circ\vphi}&\Pi'\circ F'
}$$
commutes.
\item The composition of morphisms of superfunctors
$\vphi\cl (F,\alpha_F)\to (F',\alpha_{F'})$ and $\vphi'\cl(F',\alpha_{F'})\to
(F'',\alpha_{F''})$ is the composition
$F\To[\vphi]F'\To[\vphi']F''$.
\item For a pair of superfunctors $F\cl(\shc,\Pi)\to(\shc',\Pi')$ and
$F'\cl(\shc',\Pi')\to(\shc'',\Pi'')$, the composition $F'\circ F\cl
\shc\to\shc'$ is endowed with a structure of superfunctor by
$$\xymatrix{
{F'\circ F\circ \Pi}\ar[rr]^-{\alpha_{F'\circ F}}\ar[dr]_-{F'\circ \alpha_F}
&&\Pi''\circ F'\circ F\\
&{F'\circ\Pi'\circ F}\ar[ur]_{\alpha_{F'}\circ F}
}$$

\item  The notion of {\em $\cor$-linear supercategories}
and {\em $\cor$-linear superfunctors} can be defined in a similar
way. \enum \edf

\noindent Throughout this paper, {\em by a supercategory we mean a
 $\cor$-linear additive supercategory}.

The identity functor $\id_\shc$ and the functor $\Pi\cl\shc\to\shc$
have  a structure of superfunctor
by $\alpha_{\id}=\id_\Pi\cl\id_{\shc}\circ\Pi\to\Pi\circ\id_\shc$
and $\alpha_\Pi=-\id_{\Pi\circ\Pi}$.
Note the sign.
Then $\alpha_F\cl F\circ \Pi\to \Pi'\circ F$,
as well as $\xi\cl \Pi\circ\Pi\to\id_{\shc}$,
is a morphism of superfunctors.

For $\lambda\in\cor^\times$, a supercategory  $(\shc,\Pi,\xi)$ is
equivalent to the supercategory $(\shc,\Pi,\lambda^2\xi)$. For a
supercategory $(\shc, \Pi,\xi)$, its {\em reversed supercategory}
$\shc^{\rev}$ is the supercategory $(\shc,\Pi,-\xi)$. If
$\sqrt{-1}$ exists in $\cor$, then $\shc^\rev$ is equivalent to
$\shc$ as a supercategory.

The {\em Clifford twist} of a supercategory $(\shc,\Pi,\xi)$ is the
supercategory $(\shc^{\ct},\Pi^\ct,\xi^\ct)$,  where $\shc^\ct$ is
the category whose set of objects is the set of pairs $(X,\vphi)$ of
objects $X$ of $\shc$ and isomorphisms $\vphi\cl \Pi X\isoto X$ such
that \eq&& \xymatrix@R=3ex@C=7ex{
\Pi^2 X\ar[dr]_{\Pi\vphi}\ar[rr]^{\xi_X}&&X\\
&\Pi X\ar[ur]_\vphi}\qquad \text{commutes.}
\label{eq:phi2}
\eneq
For a pair $(X,\vphi)$ and $(X'\vphi')$ of objects of $\shc^\ct$,
$\Hom_{\shc^\ct}\bigl((X,\vphi),(X'\vphi')\br)$ is
the subset of $\Hom_\shc(X,X')$ consisting of morphisms $f\cl X\to X'$
such that the following diagram commutes:
$$\xymatrix@C=8ex@R=4ex{
\Pi X\ar[r]^{\Pi f}\ar[d]_\vphi&\Pi X'\ar[d]^{\vphi'}\\
X\ar[r]^f&X'.}
$$
We define $\Pi^\ct : \shc^\ct \to \shc^\ct$ and $\xi^\ct\cl
(\Pi^\ct)^2\to\id_{\shc^\ct}$ by
$$
\begin{aligned}
& \Pi^\ct(X,\vphi)=(X,-\vphi), \\
& \xi^\ct_{(X, \vphi)}= \id_{(X, \vphi)}:
(\Pi^\ct)^2(X,\vphi)=(X,\vphi) \to (X,\vphi).
\end{aligned}
$$

\smallskip
Sometimes we simply write $\shc$ for $(\shc,\Pi,\xi)$,
$\Pi_\shc$ for $\Pi$, $\xi_\shc$ for $\xi$, 
$\shc^\rev$ for its reversed
supercategory and $\shc^\ct$ for its Clifford twist.

\smallskip

\Lemma\label{lem:cct} Let $(\shc,\Pi,\xi)$ be a supercategory.
\bnum
\item
There is a superfunctor $\Phi\cl \shc^\ct\to\shc^\rev$
defined by
$$
\begin{aligned}
 \Phi(X,\vphi) & =X, \\
\alpha_\Phi(X,\vphi) :  \Phi\circ \Pi^\ct(X,\vphi) & =X \To[{\ \vphi^{-1}}]\Pi
X=\Pi^\rev\circ\Phi(X,\vphi).
\end{aligned}
$$

\item
Assume that $\shc$ is a $\cor$-linear additive category
and $2$ is invertibel in $\cor$. 
\bna
\item
There is a superfunctor $\Psi\cl \shc^\rev\to \shc^\ct$
given by
$$
\begin{aligned}
 \Psi(X) & = (X\oplus\Pi X,\vphi_X), \\
\alpha_\Psi(X) :  \Psi\circ\Pi^\rev(X) & = \Pi X\oplus\Pi^2 X\To[\ {\psi_X}\ ] X\oplus\Pi
X=\Pi^\ct\circ\Psi X,
\end{aligned}
$$
where
$$
\begin{aligned}
&\vphi_X= \left(\begin{smallmatrix}0&\xi_X\\\id_{\Pi
X}&0\end{smallmatrix}\right) \cl \Pi X\oplus\Pi^2 X\to X\oplus\Pi X, \\
& \psi_X=\left(\begin{smallmatrix}0&-\xi_X\\\id_{\Pi
X}&0\end{smallmatrix}\right)\cl \Pi X\oplus\Pi^2 X\to X\oplus\Pi X.
\end{aligned}
$$

\item There exist canonical isomorphisms
$$\Phi\circ\Psi (X)\simeq X\oplus\Pi X \quad\text{and}\quad
\Psi\circ\Phi(Z)\simeq Z\oplus\Pi^\ct Z, $$ which are functorial in
$X\in\shc$ and $Z\in\shc^\ct$.
\ee
\item There are isomorphisms $\Pi\circ\Phi\simeq\
\Phi\circ\Pi^\ct\simeq\Phi$ and $\Pi^\ct\circ\Psi\simeq\Psi\circ\Pi\simeq
\Psi$.
\item $\Phi$ and $\Psi$ are biadjoint to each other \ro i.e.,
$\Phi$ is a left adjoint and a right adjoint of $\Psi$\rf.
\enum \enlemma

\Proof (i)\ Since the diagram
$$\xymatrix@R=3ex{
\Phi\circ(\Pi^\ct)^2(X,\vphi)\ar[r]^-\sim\ar@{=}[d]
&\Pi^\rev\circ\Phi\circ\Pi^\ct(X,\vphi)\ar[r]^-\sim\ar@{=}[d]&
(\Pi^\rev)^2\circ\Phi(X,\vphi)\ar@{=}[d]\\
\Phi(X,\vphi)\ar@{=}[d]\ar[r]^{-\vphi^{-1}}&\Pi^\rev\circ\Phi(X,-\vphi)\ar@{=}[d]
\ar[r]^{\Pi\circ\vphi^{-1}}
&(\Pi^\rev)^2\circ\Phi(X,\vphi)\ar@{=}[d]\\
X\ar[r]^{-\vphi^{-1}}&\Pi X\ar[r]^{\Pi\circ\vphi^{-1}}&\Pi^2X\ar@/^2pc/[ll]^{-\xi_X}
}
$$
commutes,
$\Phi$ is a superfunctor from $\shc^\ct$ to $\shc^\rev$.

\smallskip
\noindent (iia)
Set $Z=X\oplus \Pi X$.
Then $\vphi_X\cl \Pi Z\to Z$.
Since $\vphi_X\circ(\Pi\vphi_X)=
\left(\begin{smallmatrix}0&\xi_X\\\id_{\Pi X}&0
\end{smallmatrix}\right)
\left(\begin{smallmatrix}0&\Pi\xi_{ X}\\\id_{\Pi^2 X}&0\end{smallmatrix}\right)
=\left(\begin{smallmatrix}\xi_{X}&0\\0&\xi_{\Pi  X}\end{smallmatrix}\right)
=\xi_Z$,
we see that  $(Z,\vphi_X)=(X\oplus\Pi X,\vphi_X)$ is an object
of $\shc^\ct$.
Since the diagram
$$\xymatrix@C=23ex{
{\Pi^2X\oplus \Pi^3 X}
\ar[r]^{
{\Pi\psi_X=
\left(\begin{smallmatrix}0&-\Pi \xi_{X}\\\id_{\Pi^2 X}&0\end{smallmatrix}
\right)}}
\ar[d]_{\vphi_{\Pi X}=\left(\begin{smallmatrix}0&\xi_{\Pi X}\\\id_{\Pi^2 X}&0\end{smallmatrix}\right)}
&{\Pi X\oplus\Pi^2X}
\ar[d]^{-\vphi_X=\left(\begin{smallmatrix}0&-\xi_{X}\\-\id_{\Pi X}&0\end{smallmatrix}\right)}
\\
\Pi X\oplus\Pi^2X\ar[r]^{\psi_X=
\left(\begin{smallmatrix}0&-\xi_{X}\\\id_{\Pi X}&0
\end{smallmatrix}\right)}
&X\oplus\Pi X
}$$
commutes,
$\psi_X$ is a morphism in $\shc^\ct$.

Since the diagram
$$\xymatrix@C=10ex@R=3ex{
\Psi\circ\Pi^2 X\ar[r]^-\sim\ar@{=}[d]
&\Pi^\ct\circ\Psi\circ\Pi X\ar[r]^-\sim\ar@{=}[d]&
(\Pi^\ct)^2\circ\Psi X\ar@{=}[d]\\
\Pi^2 X\oplus\Pi^3 X\ar[r]^{
\left(\begin{smallmatrix}0&-\xi_{\Pi X}\\\id_{\Pi^2 X}&0
\end{smallmatrix}\right)}
\ar@/_2pc/[rr]_{\left(\begin{smallmatrix}-\xi_{X}&0\\0&-\xi_{\Pi X}
\end{smallmatrix}\right)}&\Pi X\oplus\Pi^2 X
\ar[r]^{\left(\begin{smallmatrix}0&-\xi_{X}\\\id_{\Pi X}&0
\end{smallmatrix}\right)}
&X\oplus\Pi X
}
$$
commutes,
$\Psi$ is a superfunctor from $\shc^\rev$ to $\shc^\ct$.

\medskip
\noi (iib)\quad
Setting $Z=(X,\vphi)\in\shc^\ct$,
$\Psi\circ \Phi(Z)=(X\oplus\Pi X,\vphi_X)
\isoto Z\oplus \Pi^\ct Z=\bl(X\oplus X,\vphi\oplus(-\vphi)\br)$
is given by
$\left(\begin{smallmatrix}\id_X&\vphi\\\id_X&-\vphi\end{smallmatrix}\right)$.

\smallskip
\noi (iii) $\left(\begin{smallmatrix}\id_X&0\\0&-\id_{\Pi X}
\end{smallmatrix}\right)\in \End(X\oplus\Pi X)$ gives an isomorphism
$\Pi^\ct\circ\Psi (X)\isoto\Psi (X)$. Similarly, we have the other
isomorphisms.

\smallskip
\noi (iv) is straightforward.
\QED

\medskip
In the rest of this paper, we assume that
\eq&&\hs{-10ex}
\text{\em{$2$ is invertible in $\cor$.}} \eneq Hence for a
$\cor$-linear idempotent-complete category $\shc$ and
$p\in\End_\shc(X)$ ($X\in\shc$) such that $p^2=\id_X$, we have the
decomposition $X\simeq \ker(\id_X-p)\oplus\ker(\id_X+p)$.

\Lemma\label{lem:dct} Let $\shc$ be a $\cor$-linear supercategory.
If $\shc$ is idempotent-complete, then $(\shc^\ct)^\ct$ is
equivalent to $\shc$ as a supercategory. \enlemma

\Proof The underlying category of $(\shc^\ct)^\ct$ is the category
$\shc'$ whose object is a triple $(X,\vphi,\psi)$, where
$\vphi\cl\Pi X\isoto X$ satisfies \eqref{eq:phi2},
$\psi\in\End_{\shc}(X)$ satisfies $\psi^2=\id_X$ and the diagram
$$\xymatrix{
\Pi X\ar[d]^{\Pi\psi}\ar[r]^{-\vphi}&X\ar[d]^\psi\\
\Pi X\ar[r]^{\vphi}&X}
$$ is commutative.
The involution $\Pi'$ of $(\shc^\ct)^\ct$ is given by
$\Pi'(X,\vphi,\psi)=(X,\vphi,-\psi)$ and the isomorphism
$\xi'\cl\Pi'^2\isoto \id_{\shc'}$ is given by
$\Pi'^2(X,\vphi,\psi)=(X,\vphi,\psi)\To[{\;\id\;}](X,\vphi,\psi)$.

Since $\psi^2=\id_X$, we have a direct sum decomposition
$X\simeq X_+\oplus X_-$ with $X_\pm=\ker(\id_X\mp \psi)$.
Moreover, the commutativity of the above diagram implies that
$\vphi$ induces an isomorphism $\vphi_{\pm}\cl\Pi X_\pm\isoto
X_\mp$. We define the functor $\Phi\cl\shc'\to \shc$  by
$$\Phi(X,\vphi,\psi)=X_+$$
and the morphism
$\alpha_\Phi\cl\Phi\circ\Pi'\to \Pi\circ\Phi$ by
$$\Phi\circ\Pi'\;(X,\vphi,\psi)\simeq X_-\isoto[{\;(\vphi_+)^{-1}\;}]
\Pi X_+\simeq\Pi\circ\Phi\;(X,\vphi,\psi).$$ Since the diagram
$$\xymatrix@R=3ex{
\Phi\circ{\Pi'}^2(X,\vphi,\psi)\ar[r]^-\sim\ar@{=}[d]
&\Pi\circ\Phi\circ\Pi'(X,\vphi,\psi)\ar[r]^-\sim\ar@{=}[d]&
\Pi^2\circ\Phi(X,\vphi,\psi)\ar@{=}[d]\\
\Phi(X,\vphi,\psi)\ar[r]^-\sim\ar@{=}[d]
&\Pi\circ\Phi(X,\vphi,-\psi)\ar[r]^-\sim\ar@{=}[d]&
\Pi^2\circ\Phi(X,\vphi,\psi)\ar@{=}[d]\\
X_+\ar[r]^{\vphi_-^{-1}}\ar[drr]_{\id_{X_+}}
&\Pi X_-\ar[r]^{\Pi\vphi_+^{-1}}&\Pi^2X_+\ar[d]^{\xi_{X_+}}\\
&&X_+
}
$$
commutes,
$\Phi$ is a superfunctor from $(\shc^\ct)^\ct$ to $\shc$.

Let $\Psi\cl\shc\to\shc'$ be the functor given by
$\Psi(X)=(X\oplus\Pi X,\vphi_X,\psi_X)$,
where $\vphi_X=\left(\begin{smallmatrix}0&\xi_X\\\id_{\Pi X}&0\end{smallmatrix}\right)$
and $\psi_X=\left(\begin{smallmatrix}\id_X&0\\0&-\id_{\Pi X}\end{smallmatrix}\right)$.
Then we can easily check that
$\Psi$ is a quasi-inverse of $\Phi$.
\QED

Let $(\shc,\Pi,\xi)$ and $(\shc',\Pi',\xi')$ be a pair of
supercategories. Let $\Fct(\shc,\shc')$ be the category of functors
from the category $\shc$ to $\shc'$. Then $\Fct(\shc,\shc')$ is
endowed with the structure of supercategory by
$\Pi_{\Fct}(F)\seteq\Pi'\circ F\circ \Pi$, and $\xi_{\Fct}\seteq
\xi'\circ
F\circ\xi\cl \Pi_\Fct^2(F)=\Pi'^2\circ F\circ\Pi^2\isoto F$.

Let us denote by
$\Fcts(\shc,\shc')$ the category of superfunctors from
$\shc$ to $\shc'$.
This category is endowed with a structure of supercategory
given by $\Pi(F,\alpha_F)=(\Pi'\circ F,-\Pi'\circ\alpha_F)$ and
$\xi(F,\alpha_F)=\xi'\circ F \cl\Pi{}^2(F,\alpha_F)=(\Pi'{}^2\circ
F,\Pi'{}^2\circ \alpha_F)\isoto (F,\alpha_F)$.
Hence $\Pi(F,\alpha_F)=(\Pi',\alpha_{\Pi'})\circ (F,\alpha_F)$.
Note the sign in the definition of
$\Pi(F,\alpha_F)$.

The category
$\Fct(\shc,\shc')^\ct$ is equivalent to the category
$\Fcts(\shc,\shc')$. However the supercategory strucure
on $\Fcts(\shc,\shc')$ induced by this equivalence
is given by $\Pi(F,\alpha_F)=(F,-\alpha_F)$.

Note that the category $\Fcts(\shc^\ct,{\shc'}^\ct)$
is equivalent to the category $\Fcts(\shc,\shc')$.

\Def Let $(\shc,\Pi,\xi)$ be a supercategory such that $\shc$ is an
exact category and $\Pi$ is an exact functor.
The {\em Grothendieck group} $\Ks(\shc)$ of $\shc$ is the abelian
group generated by $[X]$ \ro$X$ is an object of $\shc$\rf\ with the
defining relations: \bna
\item
if $0\to X'\to X\to X''\to0$ is an exact sequence, then $[X]=[X']+[X'']$,
\item $[\Pi X]=[X]$.
\ena
\edf
We have $\Ks(\shc^\rev)\simeq\Ks(\shc)$ and
$\Z[1/2]\otimes\Ks(\shc^\ct)\simeq\Z[1/2]\otimes\Ks(\shc)$
by Lemma~\ref{lem:cct}.

\subsection{Superalgebras and supermodules}
 Let $A$ be a {\em $\cor$-superalgebra}; 
i.e., a  $\cor$-algebra with a decomposition $A=A_{0}\oplus A_{1}$
such that $A_iA_j\subset A_{i+j}$ ($i,j\in\FF$).

Let $\phi_A$ be the involution of $A$ given by $\phi_A(a)=(-1)^\eps
a$ for $a\in A_\eps$ with $\eps=0,1$. Then the category of left
$A$-modules $\Mod(A)$ is naturally endowed with a structure of
supercategory. The functor $\Pi$ is 
induced by $\phi_A$. Namely, for $M\in\Mod(A)$,
$\Pi M\seteq\set{\pi(x)}{x\in M}$ with
$\pi(x)+\pi(x')=\pi(x+x')$ and
$ a\pi(x)=\pi(\phi_A(a)x)$ for any $a\in A$, $x,x'\in M$. The morphism
$\xi_M\cl \Pi^2M\to M$ is given by $\pi\bl(\pi(x)\br)\mapsto x$
($x\in M$).

An {\em $A$-supermodule} is an $A$-module $M$ with a decomposition
$M=M_{0}\oplus M_{1}$ such that $A_iM_j\subset M_{i+j}$ ($i,j\in\FF$).
An $A$-linear homomorphism between $A$-supermodules is called {\em even}
if it respects the $\FF$-grading.

Let $\Mods(A)$ be the category of $A$-supermodules with even
$A$-linear homomorphisms as morphisms. Then $\Mods(A)$ is also
endowed with a structure of supercategory. The functor $\Pi$ is
given by the parity shift: namely, $(\Pi
M)_\eps\seteq\set{\pi(x)}{x\in M_{1-\eps}}$  ($\eps=0,1$) and
$a\pi(x)=\pi(\phi_A(a)x)$ for $a\in A$ and $x\in M$. The
isomorphism $\xi_M\cl\Pi^2M\to M$ is given by $\pi\pi(x)\mapsto x$
($x\in M$).
Then there is a canonical superfunctor
$\Mods(A)\to\Mod(A)$. For an $A$-supermodule $M$, let $\phi_M\cl
M\to M$ be the map $\phi_M\vert_{M_\eps}=(-1)^\eps\id_{M_\eps}$.
Then we have $\phi_M(ax)=\phi_A(a)\phi_M(x)$ for any $a\in A$ and
$x\in M$.

\Lemma\label{lem:superct} Let $A$ be a superalgebra in which $2$ is invertible.
Then we have
\eqn
\Mod(A)^\ct&\simeq& \Mods(A)^\rev\quad\text{and}\\
(\Mods(A)^\rev)^\ct&\simeq& \Mod(A)\,.
\eneqn
\enlemma
\Proof
Let us define the superfunctor $(\Phi,\alpha)\cl \Mods(A)^\rev\to \Mod(A)^\ct$.
For $M\in \Mods(A)$,
set $\Phi(M)=(M,\psi_M)$, where $\psi_M\cl \Pi M\to M$ is given by
$\psi_M\pi(x)=\phi_M(x)$ ($x\in M$). We define $\al_M=\psi_M\cl\Pi M\to M$.
Then the diagram
$$\xymatrix@C=9ex{
\Pi^2 M\ar[d]_-{\psi_{\Pi M}}\ar[r]^-{\Pi\al_{ M}}&\Pi M\ar[d]^-{-\psi_M}\\\
\Pi M\ar[r]^{\al_M}& M }$$ is commutative, since for any $x \in M$
we have
$$
\begin{aligned}
 \al_M\psi_{\Pi M}(\pi\pi(x))&=\al_M(\phi_{\Pi M}\pi(x))=
\al_M(-\pi(\phi_M(x)))=-\phi_M\phi_M(x)=-x, \\
\psi_M\circ \Pi\alpha_{M}(\pi\pi(x))&=\psi_M(\pi(\alpha_M(\pi(x))))
=\phi_M\phi_M(x)=x.
\end{aligned}
$$
Hence $\al_M$ is a morphism from $\Phi\circ \Pi M$ to $\Pi^\ct\circ \Phi
M$ in $\Mod(A)^\ct$. Since the diagram
$$\xymatrix@R=3ex{
\Phi\circ{\Pi}^2 M\ar[r]^-\sim\ar@{=}[d]
&\Pi^\ct\circ\Phi\circ\Pi M\ar[r]^-\sim\ar@{=}[d]&
(\Pi^\ct)^2\circ\Phi M\ar@{=}[d]\\
(\Pi^2M,\psi_{\Pi^2M})\ar[r]^-{\alpha_{\Pi M}}\ar[d]^-{-\xi_M}
&(\Pi M,-\psi_{\Pi M})\ar[r]^-{\alpha_M}&
(M,\psi_M)\ar@{=}[d]\\
(M,\psi_M)\ar[rr]^{\id}&& (M,\psi_M)
}
$$
is commutative, we have defined a superfunctor $(\Phi,\alpha)\cl
\Mods(A)^\rev\to\Mod(A)^\ct$. It is straightforward to
prove that $\Phi$ is an equivalence.

The second equivalence follows from the first one and
Lemma~\ref{lem:dct}. \QED

\subsection{Tensor products}\label{supertensor}
Let $A$ and $B$ be superalgebras. We define the multiplication on
the tensor product $A \otimes B$  by
\begin{align*}
(a_1\otimes b_1)(a_2\otimes b_2)=(-1)^{\eps'_1\eps_2}(a_1a_2)\otimes (b_1b_2)
\end{align*}
for $a_i\in A_{\eps_i}$, $b_i\in B_{\eps'_i}$. Then $A \otimes B$ is
again a superalgebra.  Note that we have $A\otimes
B\cong B\otimes A$ as a superalgebra  by the supertwist map
\begin{align*}
A\otimes B\isoto B\otimes A,\quad a\otimes b\MAPSTO (-1)^{\eps_1\eps_2}b\otimes a
\quad \text{($a\in A_{\eps_1}$, $b\in B_{\eps_2}$).}
\end{align*}

\begin{example}
\label{ex:cliff}  \hfill \bnum
\item
For $n\in\Z_{\ge0}$, let $\clif_n$ be the Clifford superalgebra
generated by the odd generators $C_1,\ldots, C_n$ with the defining
relations:
$$\text{$C_i^2=1$ ($i=1,\ldots n$) and $C_iC_j=-C_jC_i$ for $i\not=j$.}$$
Hence $\clif_n\simeq(\clif_1)^{\otimes n}$.
\item
Let us denote by $\clifm_n$ the superalgebra with the same
definition but the relations $C_i^2=1$ are replaced by $C_i^2=-1$.

\enum
\end{example}
\Lemma\label{lem:MQ}
We have
\eq
&&\Mod(A\otimes_\cor\clif_1)\simeq\bl(\Mod(A)^\ct\br){}^\rev,\\[1ex]
&&\Mod(A\otimes_\cor\clif_1)\simeq\Mod_\super(A),\\[1ex]
&&\Mod(A\otimes_\cor\clifm_1)\simeq\bl(\Mod(A)^\rev{}\br)^\ct.
\label{ex:clifii} \eneq \enlemma \Proof Let us construct an
equivalence $\Phi\cl\Mod(A\otimes_\cor\clif_1)^\rev\to \Mod(A)^\ct$.
For $M\in\Mod(A\otimes_\cor\clif_1)$, we set $\Phi(M)=(\Pi
M\To[{\vphi_M}] M)$, where ${\vphi_M}\cl\Pi M\to M$ is given by
$\vphi_M(\pi(x))=C_1x$. Since $C_1^2=1$, we have
$$\vphi_M(\Pi\vphi_M)(\pi^2(x))=\vphi_M\bl(\pi(\vphi_M(\pi(x))\br)=\vphi_M(\pi(C_1x))=C_1^2x=x.$$
Hence $\Pi M\To[{\vphi_M}] M$ is an object of $\Mod(A)^\ct$.
Then $\vphi_{\Pi M}\cl \Pi^2 M\To\Pi M$ is given  by
$\vphi_{\Pi M}(\pi^2(x))=C_1\pi(x)=-\pi(C_1x)$.
In the diagram
\eq
&&\ba{c}\xymatrix@R=5ex{
\Phi\circ \Pi (M)\ar[d]^{\alpha_M}&{:}&\Pi^2 M\ar[d]_{\Pi\alpha_M}\ar[r]^-{\vphi_{\Pi M}}&\Pi M\ar[d]^{\alpha_M}\\
\Pi^\ct\circ\Phi(M)&{:}&\Pi M\ar[r]^-{-\vphi_M}&M
}\ea\label{eq:phict} \eneq we define $\alpha_M$ by
$\alpha_M(\pi(x))=C_1x$ so that we have
$(\Pi\alpha_M)(\pi^2(x))=\pi(C_1x)$. Hence the right
square in \eqref{eq:phict} commutes. Thus this defines the morphism
$\alpha\cl\Phi\circ \Pi \to\Pi^\ct\circ\Phi$. Then the diagram
$$\xymatrix@C=7ex{
\Pi^2M\ar[r]^-{\alpha_{\Pi M}}\ar[d]^-{\xi_M}&\Pi M
\ar[r]^-{\alpha_M}&M\ar[d]_{-\id_M}\\
M\ar[rr]^{\id_M}&&M}
$$
commutes, since $$\alpha_M\alpha_{\Pi M}(\pi^2x)=\alpha_M(C_1\pi(x))
=-\alpha_M(\pi(C_1x))=-C_1^2x=-x.$$ Hence $(\Phi,\alpha)$ gives a
superfunctor $\Mod(A\otimes_\cor\clif_1)\to
\bl(\Mod(A)^\ct\br)^\rev$.

The second equivalence follows from the first one and
Lemma~\ref{lem:superct}. The third equivalence follows from
$\Mod(A\otimes \clif_1\otimes\clifm_1)\simeq\Mod(A)$. \QED

\Rem\label{rem:center}
\bnum
\item
Let $\shc=(\shc,\Pi,\xi)$ be a supercategory.
Let $Z(\shc)$ be the commutative ring
$\End(\id_\shc)$.
It does not depend on $\Pi$.
Let $Z_\super(\shc)$ be its subring
$$\set{\psi\in\End(\id_\shc)}%
{\text{$\psi_{\Pi X}=\Pi(\psi_X)\in \End_\shc(\Pi X)$ for any $X\in\shc$}}.$$

Note that $Z(\shc)$ has a structure of a superalgebra,
where the involution $\phi_{Z(\shc)}$
($\phi_{Z(\shc)}\vert{Z(\shc)_\eps}=(-1)^\eps\id_{Z(\shc)_\eps})$, $\eps=0,1$) is  given by
\eq
\xymatrix@C=10ex{
\Pi^2 X\ar[r]^{\Pi f_{\Pi X}}\ar[d]^{\bwr}_{\xi_ X}&
\Pi^2 X\ar[d]_{\xi_X}^{\bwr}\\
X\ar[r]^{\phi_{Z(\shc)}(f)_X}&X
}
\eneq
for $f\in Z(\shc)$ and $X\in\shc$.

Then $Z_\super(\shc)$ is the even part of $Z(\shc)$.

\item
If $\shc=\Mod(A)$ for a superalgebra $A$, then we have
\eqn
Z(\shc)^{\op}&\simeq& Z(A)\seteq\set{a\in A}{\text{$ab=ba$ for any $b\in A$}},\\
Z_\super(\shc)^{\op}&\simeq& Z_0(A)\seteq A_0\cap Z(A),
\eneqn
and
\eqn
Z(\shc^\ct)^{\op}&\simeq& Z(A\otimes\clif_1)\\&\simeq&
\set{a\otimes 1+b\otimes C_1}{\text{$a$, $b\in A_{0}$ and $ax=xa$, $bx=\phi_A(x)b$
for any $x\in A$}},\\
Z_\super(\shc^\ct)^{\op}&\simeq&Z_\super(\shc)\simeq Z_0(A)
\eneqn
\ee
\enrem
\subsection{Weak superequivalence}\label{weak superequivalence}
Let $\shc$ and $\shc'$ be supercategories.
We say that $\shc$ and $\shc'$ are {\em weakly
superequivalent} if there exist a pair of supercategories $\shc_1$
and $\shc_2$ and superequivalences $\shc\simeq \shc_1\oplus\shc_2$
and $\shc'\simeq\shc_1\oplus \shc_2^\ct$.

If $\shc$ and $\shc'$ are weakly superequivalent and $\shc'$ and
$\shc''$ are weakly superequivalent, then $\shc$ and $\shc''$ are
weakly superequivalent. Indeed, if $\shc_1\oplus \shc_2\simeq \shc_3\oplus\shc_4$,
then there exist supercategories $\shc_{1,3}$, $\shc_{1,4}$, $\shc_{2,3}$, $\shc_{2,4}$
and superequivalences $\shc_{i}\simeq \shc_{i,3}\oplus \shc_{i,4}$ ($i=1,2$)
and $\shc_j\simeq \shc_{1,j}\oplus\shc_{2,j}$ ($j=3,4$)
such that the composition
$\shc_1\oplus\shc_2\simeq\shc_{1,3}\oplus\shc_{1,4}\oplus\shc_{2,3}\oplus\shc_{2,4}
\simeq\shc_3\oplus\shc_4$ is isomorphic to the given superequivalence.

\Def Let $A$ and $B$ be
superalgebras.
\bnum
\item
We say that $A$ and $B$ are {\em Morita superequivalent}
if $\Mod(A)$ and $\Mod(B)$ are superequivalent.
\item
We say that $A$ and $B$ are {\em weakly Morita superequivalent} 
if $\Mod(A)$ and $\Mod(B)$ are weakly superequivalent. 
\enum \edf
Note that superalgebras $A$ and $B$ are weakly Morita superequivalent if and only if
there are decompositions
$A=A_1\oplus A_2$ and $B=B_1\oplus B_2$ as superalgebras and  superequivalences
$\Mod(A_1)\simeq\Mod(B_1)$  and $\Mod(A_2)\simeq\Mod(B_2)^\ct$.
Note that $A$ and $B$ are Morita superequivalent if
and only if there exists an $(A,B)$-supermodule $P$ satisfying one
of the following equivalent conditions:
\bnum
\item
$N\mapsto P\otimes_BN$ is an equivalence of categories
from $\Mod(B)$ to $\Mod(A)$,
\item
$N\mapsto P\otimes_BN$ is an equivalence of categories
from $\Mods(B)$ to $\Mods(A)$,
\item
$M\mapsto \Hom_A(P,M)$ is an equivalence of categories
from $\Mod(A)$ to $\Mod(B)$,
\item
$M\mapsto \Hom_A(P,M)$ is an equivalence of categories
from $\Mods(A)$ to $\Mods(B)$,
\item
$M\mapsto M\otimes_AP$ is an equivalence of supercategories
from $\Mod(A^\op)$ to $\Mod(B^\op)$,
\item
$N\mapsto \Hom_{B^\op}(P,N)$ is an equivalence of supercategories
from $\Mod(B^\op)$ to $\Mod(A^\op)$,
\item
$P$ is a faithfully flat $A$-module of finite presentation and
$B\isoto\End_A(P)^\opp$,
\item
$P$ is a faithfully flat $B$-module of finite presentation and
$A\isoto\End_B(P)$. \enum

\Rem\label{rem:Morita}
\bnum
\item $\clif_n\otimes\clifm_n$ is Morita superequivalent to $\cor$.

\item
Let $A$ be a superalgebra and $e\in A$ a {\em full even idempotent};
i.e., $e\in A_{{0}}, e^2=e$ and
$A=AeA\seteq\set{\sum_{i=1}^{n}a_ieb_i}{a_i,b_i\in A, n> 0}$. Then,
$A$ and $eAe$ are Morita superequivalent.
\label{fulleven}

\item
Assume that $\sqrt{-1}$ exists in $\cor$. Then $\cor$-superalgebras
$A$ and $A\otimes\clif_n$ are weakly Morita superequivalent.
Two $\cor$-superalgebras $A$ and $B$ are weakly Morita superequivalent
if and only if there exist direct sum decompositions $A=A_1\oplus A_2$, $B=B_1\oplus B_2$ and
$A_1$ is Morita superequivalent to $B_1$ and $A_2$ is Morita superequivalent to $B_2\otimes\clif_1$
\ee
\enrem

\subsection{Self-associate simple objects}\label{sec:selfa}
Let us assume that $\cor$ is a field of characteristic $\not=2$,
and let $\shc$ be a $\cor$-linear abelian supercategory.

We say that a simple object $S$ of $\shc$ is {\em self-associate}
 if $\Pi S$ is isomorphic to $S$.
Let $\Irr(\shc)$ be the set of isomorphism classes of simple objects.
Then it is divided into the set $\Irr^{\SA}(\shc)$
of self-associate simple objects
and  the set $\Irr^{\NSA}(\shc)$ of simple objects which are not self-associate.
Note that $\Pi$ gives  an involution of $\Irr(\shc)$, and $\Irr^\SA(\shc)$
is the fixed point set.

\Lemma\label{lem:IMQ} Assume that
$\cor$ is algebraically closed and
$\shc$ satisfies the following conditions.
\bna
\item $\shc$ is an abelian category,
\item any object of $\shc$ has a finite length,
\item $\cor\isoto\End_{\shc}(S)$ for any simple object $S$ of $\shc$.
\ee
Then we have
\bnum
\item $\shc^\ct$ also satisfies the above conditions,
\item the superfunctors $\Phi\cl\shc^\ct\to\shc^\rev$ and
$\Psi\cl\shc^\rev\to\shc^\ct$ in {\rm Lemma~\ref{lem:cct}} induce
two-to-one maps $\Irr^{\NSA}(\shc^\ct)\to \Irr^{\SA}(\shc)$ and
$\Irr^{\NSA}(\shc)\to\Irr^{\SA}(\shc^\ct)$, respectively. \ee
\enlemma
\Proof
(i) is obvious. Let us show (ii).
Note that $\Phi$ and $\Psi$ are exact functors, and $\Psi$ is a left adjoint of $\Phi$.
Let $X\in \Irr^{\NSA}(\shc^\ct)$. Take a monomorphism  $S\monoto M\seteq\Phi(X)$
with a simple object $S$ of $\shc$.
Since $\Psi$ is a left adjoint of $\Psi$, we have
$g\cl\Psi(S)\to X$.
Since it is a non-zero morphism, $g$ is an epimorphism.
Since $\Pi^\ct\circ\Psi\simeq\Psi$ and since $X$ is not isomorphic to
$\Pi^\ct X$, $g$ is not an isomorphism.
Hence $S\oplus \Pi S \simeq\Phi\circ\Psi(S)\to M$ is
an epimorphism but not an isomorphism.
Hence $S\to M$ is an isomorphism.
Therefore $S\in\Irr^{\SA}(\shc)$.

Conversely, let $S$ be a self-associate simple object of $\shc$.
Let $\vphi\cl \Pi S\isoto S$ be an isomorphism.
Then $\vphi\circ\Pi(\vphi)=a\xi_S$ for some $a\in\cor^\times$.
Take $c\in\cor^\times$ such that $c^2=a$ and set
$\vphi'=c^{-1}\vphi$. Then we have $\vphi'\circ\Pi \vphi'=\xi_S$.
Then $X\seteq(S,\vphi')$ is an object of $\shc^\ct$.
It is obviously a simple object of $\shc^\ct$.
If $(S,-\vphi')=\Pi^\ct X\to X$ is an isomorphism,
then it is given by $f\cl S\isoto S$ such that
$f\circ (-\vphi')=\vphi'\circ \Pi f$.
Since $f=a\id_S$ for some $a\in \cor^\times$ it is a contradiction.
Hence $X$ is not self-associate.

Now it is obvious that if $X$ and $X'$ are simple objects of $\shc^\ct$
that are not self-associate and $\Phi(X)\simeq\Phi(X')$, then
$X'\simeq X$ or $X'\simeq \Pi X$.
\QED
\section{Quiver Hecke superalgebras}

\subsection {The quiver Hecke superalgebra $\R_n$}

Recall that {\em $2$ is assumed to be  invertible in the base ring $\cor$.} Let
$\IRED$ be a finite set with a decomposition
$\IRED=\IODD\sqcup\IEVEN$. We say that $i\in I_\even$ is an {\em
even vertex} and $i\in I_\odd$ is an {\em odd vertex}. For $i\in
\IRED$, we denote the parity of $i$ by $\prt{i}$; i.e.,
$$
\prt{i}=\begin{cases} 0 \ \ & \text{if} \ i \in \IEVEN, \\
1 \ \ & \text{if} \ i \in \IODD.
\end{cases}
$$

For each $i, j \in I$, consider the $\cor$-algebra $\At=\cor\langle
w,z \rangle/\langle zw-(-1)^{\prt{i}\,\prt{j}}wz\rangle$ generated
by the indeterminates $w, z$ with the defining relation
$$zw-(-1)^{\prt{i}\,\prt{j}}wz=0.$$
Let  $Q=(Q_{i,j}(w,z))_{i,j\in \IRED}$ be a family of
(skew-)polynomials satisfying the conditions
\eq
&&\hs{-10ex}\left\{\parbox{65ex}{ \bnum
\item $Q_{i,j}(w,z)\in \At$ for all $i, j \in \IRED$,

\vs{.5ex}
\item $Q_{i,j}(w,z)=0$ if $i=j$,
\vs{.5ex}
\item $Q_{i,j}(w,z)=Q_{j,i}(z,w)$ for all $i,j\in \IRED$,
\vs{.5ex}
\item $Q_{i,j}(w,z)=Q_{i,j}(-w,z)$ for all $i\in\IODD,j\in \IRED$.
\ee} \right.\label{eq:Q} \eneq
Hence $Q_{i,j}(w,z)$ belongs to the subalgebra
$\tAt$ of $\At$ generated by $w^{1+\prt{i}}$ and $z^{1+\prt{j}}$.
Note that $\tAt$ is commutative and isomorphic to the polynomial algebra.

\Def\label{DefRKLR} The {\em quiver Hecke superalgebra} $\R_n$ is
the $\cor$-superalgebra generated by the elements $\{x_p\}_{1\leq
p\leq n}$, $\{\tau_a\}_{1\leq a<n}$, $\{e(\nu)\}_{\nu\in\IRED^n}$
with parity
$$\prt{e(\nu)}=0, \ \  \prt{x_pe(\nu)}=\prt{\nu_p}, \ \
\prt{\tau_ae(\nu)}=\prt{\nu_a}\,\prt{\nu_{a+1}}$$ and the following
defining relations.
 \bnum
\item $e(\mu)e(\nu)=\delta_{\mu\nu}e(\mu)$ for all $\mu,\nu\in I^n$,
and $1=\sum_{\nu\in I^n}e(\nu)$,
\item \begin{align*}
x_px_qe(\nu) &=
\begin{cases}
-x_qx_pe(\nu) & \textrm{if $p\ne q$ and $\nu_p,\nu_q\in \IODD$,} \\
x_qx_pe(\nu) & \textrm{otherwise},
\end{cases}
\end{align*}
\item $\tau_ax_pe(\nu) =(-1)^{\prt{\nu_p}\;\prt{\nu_a}\;\prt{\nu_{a+1}}}
x_p\tau_ae(\nu) \ \text{if} \ p \ne a, a+1,$

\item
\begin{align*}&(\tau_ax_{a+1}-(-1)^{\prt{\nu_a}\;\prt{\nu_{a+1}}}x_a\tau_a)e(\nu)
=(x_{a+1}\tau_a-(-1)^{\prt{\nu_a}\;\prt{\nu_{a+1}}}\tau_ax_a)e(\nu)\\
&\hspace*{5ex}=\begin{cases}
e(\nu)&\text{if $\nu_a=\nu_{a+1}$,}\\
0&\text{otherwise,}
\end{cases}
\end{align*}
\item $\tau_a^2e(\nu) = Q_{\nu_a,\nu_{a+1}}(x_a,x_{a+1})e(\nu)$,
\item $\tau_a\tau_{b}e(\nu)=(-1)^{\prt{\nu_a}\prt{\nu_{a+1}}\prt{\nu_b}\prt{\nu_{b+1}}}
\tau_{b}\tau_{a}e(\nu)$ if $|a-b|>1$,
\item \eqn
&&(\tau_{a+1}\tau_a\tau_{a+1}-\tau_a\tau_{a+1}\tau_a)e(\nu)\\
&&\hs{5ex}=\left\{
\ba{l}
\dfrac{Q_{\nu_a,\nu_{a+1}}(x_{a+2},x_{a+1})-Q_{\nu_a,\nu_{a+1}}(x_{a},x_{a+1})}{x_{a+2}-x_{a}}e(\nu)
\hs{10ex}\textrm{if $\nu_a=\nu_{a+2}\in\IEVEN$,} \\[2ex]
(-1)^{\prt{\nu_{a+1}}}(x_{a+2}-x_{a})\dfrac{Q_{\nu_a,\nu_{a+1}}(x_{a+2},x_{a+1})
-Q_{\nu_a,\nu_{a+1}}(x_{a},x_{a+1})}{x^2_{a+2}-x^2_{a}}e(\nu)
\\[3ex]
\hs{53ex}\textrm{if $\nu_a=\nu_{a+2}\in \IODD$,} \\[2ex]
0 \hs{52ex} \textrm{otherwise.}
\ea\right.
\eneqn\label{eq:ttt}
\end{enumerate}
\edf

If $\IODD = \emptyset$, the quiver Hecke superalgebras $\R_n$ are
nothing but the usual Khovanov-Lauda-Rouquier algebras. Note that,
when $\nu_a$ is odd, $Q_{\nu_a,\nu_{a+1}}$ belongs to the
commutative ring $\cor[x_{a}^2,x_{a+1}]$, and hence we can define
$$\frac{Q_{\nu_a,\nu_{a+1}}(x_{a+2},x_{a+1})-Q_{\nu_a,\nu_{a+1}}(x_{a},x_{a+1})}{x^2_{a+2}-x^2_{a}}$$
as an element of the commutative ring $\cor[x_{a}^2,x_{a+1},x_{a+2}^2]$. For the
connection with a generalized Cartan matrix, see \S\,\ref{specRC}.

\Rem\label{rem:scale}
\bnum
\item 
By (v), $Q_{\nu_a,\nu_{a+1}}(x_a,x_{a+1})e(\nu)$ should be an even element.
Moreover (iv) show that $\tau_a^2e(\nu)$ commutes with
all $x_a$, $x_{a+1}$ when $\nu_a\not=\nu_{a+1}$.
Hence $Q_{i,j}(w,z)$ is in the center of
$\At$. Therefore $Q_{i,j}(w,z)$ should belong to $\tAt$.
\item
 For $\bl(\gamma_{i,j}\br)_{i,j\in I}$ with
$\gamma_{i,j}\in\cor^\times$, define $Q'=\bl(Q'_{i,j}\br)_{i,j\in I}$
by
$$Q'_{i,j}(w,z)=\gamma_{i,j}\gamma_{j,i}Q_{i,j}(\gamma_{i,i}w,\gamma_{j,j}z).$$
Then, the assignment $x_pe(\nu)\mapsto
\gamma_{\nu_p,\nu_p}x_pe(\nu)$, $\tau_ae(\nu)\mapsto
\gamma_{\nu_a,\nu_{a+1}}^{-1}\tau_ae(\nu)$ induces an isomorphism
$\R_n(Q)\isoto \R_n(Q')$.
\ee\enrem

The purpose of this paper is to investigate the basic properties of
the quiver Hecke superalgebras $\R_n$ and show that the affine
Hecke-Clifford superalgebras and affine Sergeev superalgebras are weakly
Morita superequivalent to some quiver Hecke algebras after a
completion.

\subsection{Twisted tensor product with symmetric groups}
Let $A$ be a ring with an action of the symmetric group $\Sym_n$
and let $Z(A)$ be the center of $A$. Let $a_i$ ($i=1,\ldots,n-1$) be
elements of $Z(A)$ satisfying the following conditions:
 \eq &&\text{$wa_i=a_j$ for any $w\in \Sym_n$, $1\le
i,j<n$ such that $w\{i,i+1\}=\{j,j+1\}$.} \label{eq:syma} \eneq It
is easy to see that this condition is equivalent to the conditions
\eq &&\left\{ \ba{l}
\text{$s_ja_i=a_i$ if $j\not=i-1,i+1$,}\\[1ex]
\text{$s_{i+1}a_i=s_ia_{i+1}$ for $1\le i\le n-2$.} \ea\right.
\label{eq:symb} \eneq \Prop\label{prop:symalg} Assume that a family
of elements $\{a_i\}_{1\le i<n}$ of $Z(A)$ satisfies one of the
equivalent conditions \eqref{eq:syma} or \eqref{eq:symb}. Let $R$ be
the ring generated by $A$ and $\ts_i$ $(i=1,\ldots,n-1)$ with the
defining relations: \bna
\item $A\to R$ is a ring homomorphism,
\item $\ts_i\circ a=s_i(a)\circ \ts_i$ for any $a\in A$ and  $i=1,\ldots,n-1$,
\item the $\ts_i$'s satisfy the braid relations
\eq
&&\text{$\ts_i\ts_j=\ts_j\ts_i$ for $|i-j|>1$
and $\ts_{i+1}\ts_i\ts_{i+1}=\ts_i\ts_{i+1}\ts_i$.}
\label{eq:braid}
\eneq\label{cond:braid}
\item $\ts_i^2=a_i$ for any $i=1,\ldots,n-1$.
\ena For a reduced expression $w=s_{i_1}\cdots s_{i_\ell}$ of $w\in
\Sym_n$, set $\ts_w=\ts_{i_1}\cdots\ts_{i_\ell}$. Then $\ts_w$ is
independent of the choice of a reduced expression and there exists a
linear isomorphism
$$A\otimes\Z[\Sym_n]\To R$$
given by $a\otimes w\longmapsto a\ts_w$.
\enprop
\Proof
Let $$ \eps_i(w) = \begin{cases} 0 \ \ & \text{if} \ l(s_i w) >w, \\
1 \ \ & \text{if} \ l(s_i w) < w,
\end{cases}
$$
where $l(w)$ denotes the length of a reduced expression of $w$.  We
define the action of $\ts_i$ on $A\otimes\Z[\Sym_n]$ by \eqn
&&\ts_i(a\otimes w)= s_i(a)a_i^{\eps_i(w)}\otimes s_iw. \eneqn It is
enough to show that $\ts_i$ satisfies the relations (b)--(d). Since
(b) is obvious and (d) easily follows from
$\eps_i(w)+\eps_i(s_iw)=1$, we will show (c) only.

\smallskip
\noi
We shall first show $\ts_i\ts_j=\ts_j\ts_i$ for $|i-j|>1$.
By the relation $$\ts_i\ts_j(a\otimes w)=
\ts_i(s_j(a)a_j^{\eps_j(w)}\otimes s_jw)
=s_is_j(a)(s_ia_j)^{\eps_j(w)}a_i^{\eps_i(s_jw)}\otimes s_iw,$$
the assertion reduces to
$(s_ia_j)^{\eps_j(w)}a_i^{\eps_i(s_jw)}=(s_ja_i)^{\eps_i(w)}a_j^{\eps_j(s_iw)}$,
which easily follows from \eqref{eq:symb} and $\eps_i(s_jw)=\eps_i(w)$.

\medskip
\noindent
Let us show $\ts_{i+1}\ts_i\ts_{i+1}=\ts_{i}\ts_{i+1}\ts_{i}$.
We have
\eqn
&&\ts_{i+1}\ts_i\ts_{i+1}\bl(a\otimes w\br)
=\ts_{i+1}\ts_i\bl((s_{i+1}a)a_{i+1}^{\eps_{i+1}(w)}\otimes s_{i+1}w\br)\\
&&\hs{5ex}=
\ts_{i+1}\bl((s_is_{i+1}a)(s_ia_{i+1})^{\eps_{i+1}(w)}a_i^{\eps_i(s_{i+1}w)}
\otimes s_is_{i+1}w \br)\\
&&\hs{5ex}=(s_{i+1}s_is_{i+1}a)(s_{i+1}s_ia_{i+1})^{\eps_{i+1}(w)}
(s_{i+1}a_i)^{\eps_i(s_{i+1}w)}a_{i+1}^{\eps_{i+1}(s_is_{i+1}w)}
\otimes s_{i+1}s_is_{i+1}w\\
&&\hs{5ex}=(s_{i+1}s_is_{i+1}a)a_{i}^{\eps_{i+1}(w)}
(s_{i+1}a_i)^{\eps_i(s_{i+1}w)}a_{i+1}^{\eps_{i+1}(s_is_{i+1}w)}
\otimes s_{i+1}s_is_{i+1}w. \eneqn Similarly, we have \eqn
&&\ts_{i}\ts_{i+1}\ts_{i}\bl(a\otimes w\br)\\
&&\hs{5ex}=(s_{i}s_{i+1}s_{i}a)(s_{i}s_{i+1}a_{i})^{\eps_{i}(w)}
(s_{i}a_{i+1})^{\eps_{i+1}(s_{i}w)}a_{i}^{\eps_{i}(s_{i+1}s_{i}w)}
\otimes s_{i}s_{i+1}s_{i}w\\
&&\hs{5ex}=(s_{i}s_{i+1}s_{i}a)a_{i+1}^{\eps_{i}(w)}
(s_{i+1}a_{i})^{\eps_{i+1}(s_{i}w)}a_{i}^{\eps_{i}(s_{i+1}s_{i}w)}
\otimes s_{i}s_{i+1}s_{i}w.
\eneqn
Hence the desired result follows from
$$\eps_i(s_{i+1}w)=\eps_{i+1}(s_iw),$$
which can be easily checked by reducing it to the
$\Sym_3$-case. \QED

The proof of the following lemma is straightforward.
\Lemma\label{lem:braid}
Let $(c_{i,j})_{1\le i\not=j\le n}$ be a family of elements of $Z(A)$ such that
$w(c_{i,j})=c_{w(i),w(j)}$ for any $w\in \Sym_n$.
Then the $c_{i,i+1}\ts_i$'s satisfy the braid relations
\eqref{eq:braid}.
\enlemma

\subsection{The quiver Hecke-Clifford superalgebra $\RC_n$}\label{RC_subsection}

Let $J$ be a finite index set with an involution $c\cl J\to J$. We
denote by $J^c$ the set of fixed points $\set{j\in J}{c(j)=j}$ and
let $I$ denote the set of equivalence classes under the equivalence
relation given by $i\sim_c j\Leftrightarrow i=j$ or $i=c(j)$. We
denote by $\pr$ the canonical projection $J\to I$. The symmetric
group $\Sym_n$ acts on $J^n$ in a natural way. We define $c_p\cl
J\to J$ by
$$\CC_{p}\nu=(c^{\delta_{p\ell}}\nu_{\ell})_{1\leq \ell\leq n}
\quad\text{for $\nu=(\nu_1,\ldots,\nu_n)\in J^n$.}$$

\medskip
Let $\tQ=(\tQ_{i,j}(u,v))_{i,j\in J}$ be a family of polynomials in
$\cor[u,v]$ satisfying the following conditions: for all $i, j \in
J$, we have
\eq &&\ba{l}
\text{$\tQ_{j,i}(u,v)=\tQ_{i,j}(v,u)$,}\\[1ex]
\text{$\tQ_{ci,j}(-u,v)=\tQ_{i,cj}(u,-v)=\tQ_{i,j}(u,v)$,} \\[1ex]
\text{$\tQ_{i,j}(u,v)=0$ if $\pr(i)=\pr(j)$.} \ea\label{cond:Q}
\eneq
In particular, $\tQ_{i,j}(u,v)=\tQ_{i,j}(-u,v)$ for $i\in J^c$.

\Def \label{def_RC} The {\em quiver Hecke-Clifford superalgebra}
$\RC_n$ is the $\cor$-superalgebra generated by the even generators
$\{y_p\}_{1\le p\le n}$, $\{\sigma_a\}_{1\le a<n}$,
$\{e(\nu)\}_{\nu\in J^n}$ and the odd generators $\{\cc_p\}_{1\leq
p\leq n}$ with the following defining relations: for $\mu,\nu\in
J^n$, $1\leq p,\,q\leq n$, $1\leq a\leq n-1$, we have

\bnum
\item $e(\mu)e(\nu)=\delta_{\mu\nu}e(\mu)$, $1=\sum_{\nu\in J^n}e(\nu)$,
$y_pe(\nu)=e(\nu)y_p$,  $\cc_pe(\nu)=e(\CC_p\nu)\cc_p$,\label{MA1}
\item $y_py_q=y_qy_p$, $\cc_p\cc_q+\cc_q\cc_p=2\delta_{pq}$,\label{MA2}
\item $\cc_py_q=(-1)^{\delta_{p,q}}y_q\cc_p$, \label{MA3}
\item $\sigma_ae(\nu) = e(s_a\nu)\sigma_a, \sigma_a \cc_p = \cc_{s_a(p)}\sigma_a$,
\item $\sigma_a y_pe(\nu) = y_{p}\sigma_ae(\nu)$ if $p\not=a,a+1$,
\item \begin{align*}
(\sigma_ay_{a+1}-y_a\sigma_a)e(\nu) &=
\begin{cases}
e(\nu) & \nu_a=\nu_{a+1}\not\in J^c, \\
-\cc_a\cc_{a+1}e(\nu) & \nu_a=c\nu_{a+1}\not\in J^c, \\
(1-\cc_a\cc_{a+1})e(\nu) & \nu_a=\nu_{a+1}\in J^c,\\
0&\text{otherwise,}
\end{cases}
\end{align*}
or equivalently,
$$\sigma_ay_{a+1}-y_a\sigma_a=
\sum_{\nu_a=\nu_{a+1}}e(\nu)-\sum_{\nu_a=c\nu_{a+1}}\cc_a\cc_{a+1}e(\nu),$$
\item \begin{align*}
(y_{a+1}\sigma_a-\sigma_ay_a)e(\nu) &=
\begin{cases}
e(\nu) & \nu_a=\nu_{a+1}\not\in J^c, \\
\cc_a\cc_{a+1}e(\nu) & \nu_a=c\nu_{a+1}\not\in J^c, \\
(1+\cc_a\cc_{a+1})e(\nu) & \nu_a=\nu_{a+1}\in J^c,\\
0&\text{otherwise,}
\end{cases}
\end{align*}
or equivalently, $$y_{a+1}\sigma_a-\sigma_ay_a=
\sum_{\nu_a=\nu_{a+1}}e(\nu)+\sum_{\nu_a=c\nu_{a+1}}\cc_a\cc_{a+1}e(\nu),$$
\item $\sigma_a^2e(\nu) = \tQ_{\nu_a,\nu_{a+1}}(y_a,y_{a+1})e(\nu)$,
\item $\sigma_a\sigma_{b}=\sigma_{b}\sigma_{a}$ if $|a-b|>1$,
\item 
\begin{align*}
&(\sigma_{a+1}\sigma_a\sigma_{a+1}-\sigma_a\sigma_{a+1}\sigma_a)e(\nu) \\
&\hs{4ex}=\begin{cases}
\dfrac{\tQ_{\nu_a,\nu_{a+1}}(y_{a+2},y_{a+1})
-\tQ_{\nu_a,\nu_{a+1}}(y_{a},y_{a+1})}{y_{a+2}-y_{a}} e(\nu)
& \textrm{if $\nu_a=\nu_{a+2}\not\in J^c$}, \\[3ex]
\dfrac{\tQ_{\nu_a,\nu_{a+1}}(y_{a+2},y_{a+1})-\tQ_{\nu_a,\nu_{a+1}}(-y_{a},y_{a+1})}{y_{a+2}+y_{a}}
\cc_a\cc_{a+2}e(\nu)
& \textrm{if $\nu_a=c\nu_{a+2}\not\in J^c$}, \\[3ex]
\dfrac{\tQ_{\nu_a,\nu_{a+1}}(y_{a+2},y_{a+1})
-\tQ_{\nu_a,\nu_{a+1}}(y_{a},y_{a+1})}{y_{a+2}-y_{a}}e(\nu) &\\[1ex]
\hs{5ex}+\dfrac{\tQ_{\nu_a,\nu_{a+1}}(y_{a+2},y_{a+1})
-\tQ_{\nu_a,\nu_{a+1}}(y_{a},y_{a+1})}{y_{a+2}+y_{a}} \cc_a\cc_{a+2}e(\nu)
& \textrm{if $\nu_a=\nu_{a+2}\in J^c$}, \\[3ex]
0 & \textrm{otherwise},
\end{cases}
\end{align*}
or equivalently,
\eqn
\sigma_{a+1}\sigma_a\sigma_{a+1}-\sigma_a\sigma_{a+1}\sigma_a
&=&\kern-3ex\sum_{\nu_a=\nu_{a+2}}
\dfrac{\tQ_{\nu_a,\nu_{a+1}}(y_{a+2},y_{a+1})-\tQ_{\nu_a,\nu_{a+1}}(y_{a},y_{a+1})}%
{y_{a+2}-y_{a}}e(\nu) \\
&&\hs{-1ex}+\kern-2ex\sum_{\nu_a=c\nu_{a+2}}
\dfrac{\tQ_{\nu_a,\nu_{a+1}}(y_{a+2},y_{a+1})-\tQ_{\nu_a,\nu_{a+1}}(-y_{a},y_{a+1})}{y_{a+2}+y_{a}}
\cc_a\cc_{a+2}e(\nu).
\eneqn
\end{enumerate}
\edf

 We first calculate the commutation relation
between $\sigma_a$'s and the polynomials in $y_k$ ($k=1,\ldots,n$).
Set
\eqn&& \ba{l} \erm_{a,b}=\ssum_{\nu_a=\nu_{b}\in J}e(\nu), \quad
\ec_{a,b}=\ssum_{\nu_{a}=c\nu_b\in J}e(\nu), \\[1ex]
\erm_a=\erm_{a,a+1},\quad \ec_a=\ec_{a,a+1}.
\ea
\eneqn
Hence we have
\eq&&\ba{rcl}
\sigma_ay_{a+1}-y_a\sigma_a&=&\erm_a-\cc_a\cc_{a+1}\ec_a,\\[1ex]
y_{a+1}\sigma_a-\sigma_ay_a&=&\erm_a+\cc_a\cc_{a+1}\ec_a.
\ea
\eneq
\Lemma\label{lem:sigmacom} For $1\le a<n$ and
$f(y_1,\ldots,y_n)\in\cor[y_1,\ldots,y_n]$, we have
\eq&&
\sigma_a\circ
f=s_a(f)\circ\sigma_a+\dfrac{f-s_a(f)}{y_{a+1}-y_a}\erm_a
-\cc_a\cc_{a+1}\dfrac{f-\bs_a f}{y_{a+1}+y_a}\ec_a,
\label{com:sigmf} \eneq where \eq&& \ba{l}
(s_af)(y_1,\ldots,y_n)=f(y_1,\ldots,y_{a-1},y_{a+1},y_a,y_{a+2},\ldots,y_n),\\[1ex]
(\bs_af)(y_1,\ldots,y_n)=f(y_1,\ldots,y_{a-1},-y_{a+1},-y_a,y_{a+2},\ldots,y_n).
\ea \eneq \enlemma

\Proof We will use induction on the degree of $f$. Assume that
\eqref{com:sigmf} holds for $f$. We will show that it holds for $y_p
f$ as well.

It is evident for $p\not=a,a+1$. If $p=a$, then \eqn&&\hs{-5ex}
\ba{l}
\sigma_a\circ (y_af)=(y_{a+1}\sigma_a-\erm_a-\cc_a\cc_{a+1}\ec_a)\circ f\\[1ex]
\hs{5ex}=y_{a+1}\Bigl(s_a(f)\circ\sigma_a+\dfrac{f-s_a(f)}{y_{a+1}-y_a}\erm_a
-\cc_a\cc_{a+1}\dfrac{f-\bs_a f}{y_{a+1}+y_a}\ec_a\Bigr)
-f\erm_a-\cc_a\cc_{a+1}f\ec_a\\[2ex]
\hs{5ex}=s_a(y_af)\circ\sigma_a+\left(\dfrac{y_{a+1}(f-s_a(f))}{y_{a+1}-y_a}-f\right)\erm_a
+\cc_a\cc_{a+1}\left(\dfrac{y_{a+1}(f-\bs_a f)}{y_{a+1}+y_a}-f\right)\ec_a\\[2.5ex]
\hs{5ex}=s_a(y_af)\circ\sigma_a+\dfrac{y_af-s_a(y_af)}{y_{a+1}-y_a}\erm_a
+\cc_a\cc_{a+1}\dfrac{-y_af+\bs_a (y_af)}{y_{a+1}+y_a}\ec_a,
\ea
\eneqn
and if $p=a+1$, then
\eqn&&\hs{-5ex}
\ba{l}
\sigma_a\circ (y_{a+1}f)\\[1ex]
\hs{5ex}=(y_{a}\sigma_a+\erm_a-\cc_a\cc_{a+1}\ec_a)\circ f\\[2ex]
\hs{5ex}=y_{a}\Bigl(s_a(f)\circ\sigma_a+\dfrac{f-s_a(f)}{y_{a+1}-y_a}\erm_a
-\cc_a\cc_{a+1}\dfrac{f-\bs_a f}{y_{a+1}+y_a}\ec_a\Bigr)
+f\erm_a-\cc_a\cc_{a+1}f\ec_a\\[2ex]
\hs{5ex}=s_a(y_{a+1}f)\circ\sigma_a+\left(\dfrac{y_{a}(f-s_a(f))}{y_{a+1}-y_a}+f\right)\erm_a
+\cc_a\cc_{a+1}\left(\dfrac{y_{a}(f-\bs_a f)}{y_{a+1}+y_a}-f\right)\ec_a\\[2ex]
\hs{5ex}=s_a(y_{a+1}f)\circ\sigma_a
+\dfrac{y_{a+1}f-s_a(y_{a+1}f)}{y_{a+1}-y_a}\erm_a
+\cc_a\cc_{a+1}\dfrac{-y_{a+1}f+\bs_a (y_{a+1}f)}{y_{a+1}+y_a}\ec_a,
\ea \eneqn which proves our assertion.
\QED

\noi Note that the $\cor$-algebra $(\cor^{J})^{\otimes n}$ and
$\clif_n$ are canonically identified with $\soplus\nolimits_{\nu\in
J^n}\cor e(\nu)$ and $\langle \{\cc_p\}_{1\leq p\leq n} \mid
\cc_p\cc_q+\cc_q\cc_p=2\delta_{pq}\rangle$, respectively. We denote
by $\MA_n$ the superalgebra generated by $\{y_p, \cc_p,e(\nu)\mid
1\leq p\leq n,\;\nu\in J^n\}$ with the defining relations
\eqref{MA1}--\eqref{MA3} given above. Clearly, we have a
$\cor$-supermodule isomorphism
\begin{align*}
(\cor^{J})^{\otimes n}\otimes \cor[y_1,\ldots,y_n]\otimes \clif_n\isoto \MA_n,\quad
e({\nu})\otimes f\otimes c\longmapsto e(\nu)fc,
\end{align*}
where $\nu\in J^n,f\in \cor[y_1,\ldots,y_n]$ and $c\in\clif_n$.

\Rem
The algebras $\MA_n$ and $\RC_n$ have an anti-involution
that sends the generators $y_p$, $\cc_p$, $\sigma_a$, $e(\nu)$ to themselves.
\enrem

\subsection{Realization of quiver Hecke-Clifford superalgebras}

Let $\MK_n$ be the superalgebra
$\cor[y_1,\ldots,y_n][(y_a-y_{b})^{-1},(y_a+y_{b})^{-1}; 1\le a<b\le n]
\otimes_{\cor[y_1,\ldots,y_n]}\MA_n$.
The symmetric group $\Sym_n$ acts on $\MA_n$ and $\MK_n$.

For $1\le a, b\le n$ with $a\not=b$, set \eq&& \ba{rl}
\tQ_{a,b}&=\sum\limits_{\nu\in J^n}\tQ_{\nu_a,\nu_{b}}(y_a,y_{b})e(\nu),\\[2ex]
R_{a,b}&=\tQ_{a,b}-(y_a-y_{b})^{-2}\erm_{a,b}
-(y_a+y_{b})^{-2}\ec_{a,b}\\[1ex]
&=\sum\limits_{\pr(\nu_a)\not=\pr(\nu_b)}\tQ_{\nu_a,\nu_{b}}(y_a,y_{b})e(\nu)
-\sum\limits_{\nu_a=\nu_b\notin J^c}\dfrac{1}{(y_a-y_{b})^2}e(\nu)\\
&\hs{8ex}
-\sum\limits_{\nu_a=c\nu_b\notin J^c}\dfrac{1}{(y_a+y_{b})^2}e(\nu)
-\sum\limits_{\nu_a=\nu_b\in J^c}2\dfrac{y_a^2+y_{b}^2}{(y_a^2-y_{b}^2)^2}
e(\nu)\in\MK_n.\ea
\label{def R}
\eneq
Then the $R_{a,b}$'s belong to the center of $\MK_n$ and satisfy the properties
$$
R_{a,b}=R_{b,a}\quad\text{and}\quad w(R_{a,b})=R_{w(a),w(b)}
\quad\text{for all $w\in \Sym_n$.}$$

Let us denote by $\SK_n$ the $\cor$-superalgebra generated by the
$\cor$-superalgebra $\MK_n$ and $\ts_k$ ($k=1,\ldots,n-1$)
satisfying the defining relations:
\eq&&\left\{
\ba{l}\text{the
$\ts_k$'s satisfy the braid relations \eqref{eq:braid},}
\\[1ex]
\text{$\ts_k^2=R_{k,k+1}$ and,}\\[1ex]
\text{$\ts_k\circ a=(s_ka)\circ\ts_k$ for all $a\in \MK_n$.}
\ea\right. \eneq Then by Proposition \ref{prop:symalg}, we have
a $\cor$-linear isomorphism
\eq
&&\MK_n\otimes\cor[\Sym_n]\isoto \SK_n. \eneq

\Th\label{poly_rep} Assume that $(\tQ_{i,j})_{i,j\in J}$ satisfies
the conditions \eqref{cond:Q}. Set \eqn &&\sigma_ae(\nu)
=\begin{cases}
\ts_ae(\nu)&\text{if $\pr(\nu_a)\not=\pr(\nu_{a+1})$,}\\
(\ts_a-(y_a-y_{a+1})^{-1})e(\nu)&\text{if $\nu_a=\nu_{a+1}\not\in J^c$,}\\
(\ts_a+(y_a+y_{a+1})^{-1}\cc_a\cc_{a+1})e(\nu)&
\text{if $\nu_a=c\nu_{a+1}\not\in J^c$,}\\
\bigl(\ts_a-(y_a-y_{a+1})^{-1}
+(y_a+y_{a+1})^{-1}\cc_a\cc_{a+1}\bigr)e(\nu)&
\text{if $\nu_a=\nu_{a+1}\in J^c$.}
\end{cases}
\eneqn Then the $\sigma_a$'s satisfy the commutation relations in
{\rm Definition \ref{def_RC}}. Moreover, the superalgebra
homomorphism $\RC_n\To\SK_n$ thus obtained is injective.
\enth

\Proof For $1\le a, b\le n$ with $a\not=b$, set \eq
f_{a,b}=(y_b-y_a)^{-1}\erm_{a,b}+(y_b+y_a)^{-1}\cc_a\cc_b\ec_{a,b}\in\MK_n,
\label{def:f} \eneq and $f_a=f_{a,a+1}$. Then we have
$$\sigma_a=\ts_a+f_a$$
and $f_{a,b}$'s satisfy the relations
$$
\begin{aligned}
w(f_{a,b})& =f_{w(a),w(b)} \ \ \text{for any} \  w\in\Sym_n, \\
f_{a,b} & =-f_{b,a}.
\end{aligned}
$$
We can easily verify that $f_{a}\cc_b=\cc_{s_a(b)}f_a$.

Since $$f_a^2=(y_a-y_{a+1})^{-2}\erm_a +(y_a+y_{a+1})^{-2}\ec_a,$$ we
obtain \eq &&R_{a,b}=\tQ_{a,b}-f_{a,b}^2. \label{eq:QRf} \eneq

\bigskip
\noi Note that the commutation relations (i)--(v)  and (ix) are
obvious. Moreover, (vi) follows from
$f_{a,b}y_{b}-y_af_{a,b}=\erm_{a,b}-\cc_a\cc_b\ec_{a,b}$ and (vii)
can be verified in a similar manner.

\smallskip

For the relation (viii), since $\ts_a$ and $f_a$ anti-commute with
each other, we have
$$\sigma_a^2=\ts_a^2+f_a^2
=R_{a,a+1}+f_{a}^2 =\tQ_{a,a+1}$$ as desired, where the last
equality follows from \eqref{eq:QRf}.
\smallskip

It remains to verify the relation (x). By the definition, we have
\eqn &&\sigma_{a+1}\sigma_a\sigma_{a+1}
=\sigma_{a+1}\sigma_a(\ts_{a+1}+f_{a+1})\\
&&\hs{5ex}=\sigma_{a+1}(\ts_a\ts_{a+1}+(s_af_{a+1})\ts_a+f_a\ts_{a+1}+f_af_{a+1})\\
&&\hs{5ex}=\ts_{a+1}\ts_a\ts_{a+1}+(s_{a+1}s_af_{a+1})\ts_{a+1}\ts_a
+(s_{a+1}f_{a})\ts_{a+1}^2+(s_{a+1}f_a)(s_{a+1}f_{a+1})\ts_{a+1}\\
&&\hs{8ex}
+f_{a+1}\ts_a\ts_{a+1}+f_{a+1}(s_af_{a+1})\ts_a
+f_{a+1}f_a\ts_{a+1}+f_{a+1}f_af_{a+1}\\
&&\hs{5ex}=\ts_{a+1}\ts_a\ts_{a+1}+f_{a}\ts_{a+1}\ts_a
+f_{a+1}\ts_a\ts_{a+1}\\
&&\hs{8ex} +(f_{a+1}f_a-f_{a,a+2}f_{a+1})\ts_{a+1}
+f_{a+1}f_{a,a+2}\ts_a+f_{a+1}f_af_{a+1}+f_{a,a+2}R_{a+1,a+2}.
\eneqn
Similarly, we have \eqn
&&\sigma_{a}\sigma_{a+1}\sigma_{a}\\
&&\hs{5ex}=\ts_{a}\ts_{a+1}\ts_{a}+(s_{a}s_{a+1}f_{a})\ts_{a}\ts_{a+1}
+(s_{a}f_{a+1})\ts_{a}^2+(s_{a}f_{a+1})(s_{a}f_{a})\ts_{a}\\
&&\hs{8ex}
+f_{a}\ts_{a+1}\ts_{a}+f_{a}(s_{a+1}f_{a})\ts_{a+1}
+f_{a}f_{a+1}\ts_{a}+f_{a}f_{a+1}f_{a}\\
&&\hs{5ex}=\ts_{a}\ts_{a+1}\ts_{a}+f_{a+1}\ts_{a}\ts_{a+1}
+f_{a}\ts_{a+1}\ts_{a}\\
&&\hs{8ex} +(f_{a}f_{a+1}-f_{a,a+2}f_{a})\ts_{a}
+f_{a}f_{a,a+2}\ts_{a+1}+f_{a}f_{a+1}f_{a}+f_{a,a+2}R_{a,a+1}.
\eneqn Hence our assertion would follow from the following
equalities:
\eq\label{eq:fff} &&\ba{l}
f_{a}f_{a+1}=f_{a,a+2}f_{a}+f_{a+1}f_{a,a+2},\quad
f_{a+1}f_{a}=f_{a,a+2}f_{a+1}+f_{a}f_{a,a+2},\\[2ex]
f_{a+1}f_af_{a+1}-f_{a}f_{a+1}f_{a}
+f_{a,a+2}(R_{a+1,a+2}-R_{a,a+1})\\[1ex]
\hs{5ex}=\sum\limits_{\nu_a=\nu_{a+1}}
\dfrac{\tQ_{\nu_a,\nu_{a+1}}(y_{a+2},y_{a+1})-\tQ_{\nu_a,\nu_{a+1}}(y_{a},y_{a+1})}%
{y_{a+2}-y_{a}}e(\nu) \\[1ex]
\hs{15ex}+\kern-2ex\sum\limits_{\nu_a=c\nu_{a+1}}
\dfrac{\tQ_{\nu_a,\nu_{a+1}}(y_{a+2},y_{a+1})-\tQ_{\nu_a,\nu_{a+1}}(-y_{a},y_{a+1})}{y_{a+2}+y_{a}}
\cc_a\cc_{a+2}e(\nu).
\ea \eneq

In order to show the first and the second equalities, let us show
\eq&&f_{a,b}f_{b,c}+f_{b,c}f_{c,a}+f_{c,a}f_{a,b}=0 \quad\text{for
distinct elements $a,b,c\in\{1,\ldots,n\}$.} \label{eq:abc} \eneq
Set $h_{a,b}=(y_b-y_a)^{-1}\erm_{a,b}$. Then we have
$$f_{a,b}=h_{a,b}+\cc_ah_{a,b}\cc_b=h_{a,b}+\cc_bh_{a,b}\cc_a,$$
which yields
\eqn
f_{a,b}f_{b,c}&=&(h_{a,b}+\cc_bh_{a,b}\cc_a)(h_{b,c}+\cc_ch_{b,c}\cc_b)\\
&=&h_{ab}h_{bc}+\cc_bh_{a,b}h_{b,c}\cc_a+\cc_c h_{a,b}h_{b,c}\cc_b+
\cc_ah_{a,b}h_{b,c}\cc_c.
\eneqn
Hence \eqref{eq:abc} is a consequence of
$$h_{a,b}h_{b,c}+h_{b,c}h_{c,a}+h_{c,a}h_{a,b}=
h_{a,b}h_{b,c}h_{c,a}\bl((y_a-y_c)+(y_b-y_a)+(y_c-y_b)\br)=0.$$

Let us prove the last equality.
Since $f_{a,b}^2$ belongs to the center of $\MK_n$, by the first
equality in \eqref{eq:fff}, we have \eqn
&&f_{a+1}f_af_{a+1}-f_{a}f_{a+1}f_{a}+f_{a,a+2}(f_{a}^2-f_{a+1}^2)\\
&&\hs{10ex}=f_{a+1}(f_af_{a+1}-f_{a+1}f_{a,a+2})+(-f_af_{a+1}+f_{a,a+2}f_{a})f_{a}\\
&&\hs{10ex}=f_{a+1}(f_{a,a+2}f_a)+(-f_{a+1}f_{a,a+2})f_{a} =0.
\eneqn
It follows that
$$f_{a+1}f_af_{a+1}-f_{a}f_{a+1}f_{a}
+f_{a,a+2}(R_{a+1,a+2}-R_{a,a+1})
=f_{a,a+2}(\tQ_{a+1,a+2}-\tQ_{a,a+1}).$$

Since $\tQ_{a,b}$ is in $Z(\MA_n)$, we have
$$f_{a,a+2}(\tQ_{a+1,a+2}-\tQ_{a,a+1})
=h_{a,a+2}(\tQ_{a+1,a+2}-\tQ_{a,a+1})
+\cc_ah_{a,a+2}(\tQ_{a+1,a+2}-\tQ_{a,a+1})\cc_{a+2},$$ as desired.

\medskip
Thus we have constructed a superalgebra homomorphism $\RC_n\to
\SK_n$. For each $w\in \Sym_n$, we choose
a reduced expression $s_{i_1}\cdots s_{i_{\ell}}$ of $w$, and set
$\sigma_w=\sigma_{i_1}\cdots \sigma_{i_{\ell}}$. Then, by the
commutation relations, it is easy to see that $\RC_n$ is generated
by
$$
\set{
y^a\cc^{\eta}\sigma_{w}e(\nu)}%
{a\in\Z^n_{\geq 0},\, \eta\in(\FF)^n,\,w\in\Sym_n,\,\nu\in J^n}
$$
(see Corollary ~\ref{cor:pbwrc} below for the notation)
as a $\cor$-module. It is straightforward to verify that its image is
linearly independent in $\SK_n$. Hence $\RC_n\to\SK_n$ is injective.
\QED

Theorem \ref{poly_rep} immediately implies the following corollary.

\Cor\label{cor:pbwrc} For each $w\in \Sym_n$, we choose
a reduced expression $s_{i_1}\cdots s_{i_{\ell}}$
of $w$, and set $\sigma_w=\sigma_{i_1}\cdots
\sigma_{i_{\ell}}$. Then
\begin{align*}
\set{
y^a\cc^{\eta}\sigma_{w}e(\nu)}%
{a\in\Z^n_{\geq 0},\, \eta\in(\FF)^n,\,w\in\Sym_n,\,\nu\in J^n}
\end{align*}
is a basis of $\RC_n$. Here,
$y^a=y_1^{a_1}\cdots y_n^{a_n}$ for
$a=(a_1,\ldots,a_n)\in \Z_{\ge0}$ and
$\cc^{\eta}=\cc_1^{\eta_1}\cdots \cc_n^{\eta_n}$
for $\eta=(\eta_1,\cdots,\eta_n)\in(\FF)^n$.
\encor

\subsection{Weak Morita superequivalence between $\RC_n$ and $\R_n$}
Let $J$, $c\cl J\to J$, $J^c$, $\pr\cl J\to I$ and $\tQ=(\tQ_{i,j})_{i,j\in J}$
be as in the preceding subsection.

Set $I_\odd\seteq\pr(J^c)$ and $I_\even\seteq I\setminus I_\odd$.
Then $\#\pr^{-1}(i)=1$ or $2$ according as $i\in I_\odd$ or $i\in
I_\even$.
Choose $J^\pl\subset J$ such that the projection $\pr$ induces a
bijection $J^\pl\to I$ and let $\qr\cl I\to J^\pl$ be its inverse.

\Def\label{def_of_RCp} We define $\RCP_n=\EP\RC_n\EP$, where
$\EP=\sum_{\nu\in J^\pl{}^{n}}e(\nu)$. \edf

Since $\cc_{a_1}\cdots\cc_{a_k}e(\nu)\cc_{a_k}\cdots\cc_{a_1}=e(c_{a_1}\cdots c_{a_k}\nu)$,
we have $\RC_n\EP\RC_n=\RC_n$. Hence
Remark~\ref{rem:Morita}~\eqref{fulleven} implies that $\RC_n$
and $\RCP_n$ are Morita superequivalent.

It is easy to see that $\RCP_n$ is generated by the even generators
$\{y_p\EP\}_{1\leq p\leq n}$, $\{\sigma_a\EP\}_{1\leq a<n}$,
$\{e(\nu)\}_{\nu\in J^\pl{}^n}$ and the odd generators
$\{\EP\cc_p\EP\}_{1\leq p \leq n}$. For simplicity, we will write
$y_p$ for $y_p\EP$, $\sigma_a$ for $\sigma_a\EP$,
$\{e(\nu)\}_{\nu\in I^n}$ for $\{e(\nu)\}_{\nu\in J^\pl{}^n}$, and
$\cc_p$ for $\cc_p\EP$. Then the defining relations  for these
generators are give as follows:

\bnum
\item $e(\mu)e(\nu)=\delta_{\mu\nu}e(\mu)$ for all $\mu,\nu\in I^n$,
 $1=\sum_{\nu\in I^n}e(\nu)$,
\item $y_py_q=y_qy_p$, $\cc_p\cc_q=-\cc_q\cc_p$
for all $1\leq p<q\leq n$,
\item $y_pe(\nu)=e(\nu)y_p$ and $\cc_pe(\nu)=e(\nu)\cc_p$ for any $\nu\in I^n$,
\item
$\cc_pe(\nu)=0$ if $\nu_p\in I_\even$, and $\cc_p^2e(\nu)=e(\nu)$ if
$\nu_p\in I_\odd$,
\item $\cc_py_q=(-1)^{\delta_{pq}}y_q\cc_p$,
\item $\sigma_ae(\nu) = e(s_a\nu)\sigma_a, \sigma_a \cc_p = \cc_{s_a(p)}\sigma_a$,
\item $\sigma_a y_pe(\nu) = y_{p}\sigma_ae(\nu)$ if $p\ne a,a+1$,
\item \begin{align*}
(\sigma_ay_{a+1}-y_a\sigma_a)e(\nu) &=
\begin{cases}
e(\nu) & \text{if $\nu_a=\nu_{a+1}\in I_\even$,} \\
(1-\cc_a\cc_{a+1})e(\nu) & \text{if $\nu_a=\nu_{a+1}\in\IODD$,}\\
0&\text{if $\nu_a\ne\nu_{a+1}$.}
\end{cases}
\end{align*}
\item \begin{align*}
(y_{a+1}\sigma_a-\sigma_ay_a)e(\nu) &=
\begin{cases}
e(\nu) & \text{if $\nu_a=\nu_{a+1}\in\IEVEN$,} \\
(1+\cc_a\cc_{a+1})e(\nu) &\text{if $\nu_a=\nu_{a+1}\in\IODD$,}\\
0&\text{if $\nu_a\ne\nu_{a+1}$.}
\end{cases}
\end{align*}
\item $\sigma_a^2e(\nu) = \QTL_{\nu_a,\nu_{a+1}}(y_a,y_{a+1})e(\nu)$,
\item $\sigma_a\sigma_{b}=\sigma_{b}\sigma_{a}$ if $|a-b|>1$,
\item
\eqn
&&(\sigma_{a+1}\sigma_a\sigma_{a+1}-\sigma_a\sigma_{a+1}\sigma_a)e(\nu)\\
&&\hs{5ex}=
\begin{cases}
\dfrac{\QTL_{\nu_a,\nu_{a+1}}(y_{a+2},y_{a+1})-\QTL_{\nu_a,\nu_{a+1}}(y_{a},y_{a+1})}{y_{a+2}-y_{a}} e(\nu)
& \textrm{if $\nu_a=\nu_{a+2}\in \IEVEN$,} \\[3ex]
\dfrac{\QTL_{\nu_a,\nu_{a+1}}(y_{a+2},y_{a+1})-\QTL_{\nu_a,\nu_{a+1}}(y_{a},y_{a+1})}{y_{a+2}-y_{a}}e(\nu)&\\[0.5ex]
\hs{2ex}- \cc_a\cc_{a+2}\dfrac{\QTL_{\nu_a,\nu_{a+1}}(y_{a+2},y_{a+1})
-\QTL_{\nu_a,\nu_{a+1}}(y_{a},y_{a+1})}{y_{a+2}+y_{a}}e(\nu)
& \textrm{if $\nu_a=\nu_{a+2}\in \IODD$,} \\[3ex]
0 & \textrm{otherwise.}
\end{cases}
\eneqn
\end{enumerate}

Corollary~\ref{cor:pbwrc} implies that the set
\eq &&\ba{rl} \{ y^a
\cc^{\eta}e(\nu)\sigma_{w}\mid& a \in\Z^n_{\geq 0}
,\,\eta\in(\FF)^n,\,\nu\in
I^n, \,w\in\Sym_{n}\\
&\text{such that $\eta_p=0$ as soon as $\nu_p\in I_\even$}\}
\ea\label{eq:pbwCl}
\eneq
is a basis of $\RCP_n$.

\Def\label{def:At} For $i,j\in I$, recall that $\At =\cor\lan
w,z\ran/\lan zw-(-1)^{\prt{i}\,\prt{j}}wz\ran$. We
define $\tAt$ to be the subalgebra of $\At$ generated by
$w^{1+\prt{i}}$ and $z^{1+\prt{j}}$.
Similarly, we define $\Bt$ to be the commutative polynomial ring
$\cor[u,v]$ and denote by $\tBt$ its subalgebra generated by
$u^{1+\prt{i}}$ and $v^{1+\prt{j}}$.
\edf
Then $\tAt$ is
commutative, and $\tBt$ and $\tAt$ are isomorphic by the correspondence
\eq
&&u^{1+\prt{i}}=(-1)^{\prt{i}}w^{1+\prt{i}}, \quad
v^{1+\prt{j}}=(-1)^{\prt{j}}z^{1+\prt{j}}. \label{eq:AB} \eneq

Let $\tQ=(\tQ_{i,j}(u,v))_{i,j\in J}$ be as in \eqref{cond:Q}. For
$i,j\in I$, we denote by $Q_{i,j}$ the element of $\tAt\subset\At$
corresponding to $\tQ_{\qr(i),\qr(j)}\in\tBt[{\qr(i),\qr(j)}]$. Then
we can easily see that $Q\seteq(Q_{i,j})_{i,j\in I}$ satisfies the
condition \eqref{eq:Q}, and hence it defines the quiver Hecke
superalgebra $\R_n=\R_n(Q)$ (see \ref{DefRKLR}).

\Rem A different choice of $J^\pl$ would yield a different matrix
$Q'$ which is within the re-scaling of $Q$ given in
Remark~\ref{rem:scale}. Hence the quiver Hecke superalgebra $\R_n(Q)$
does not depend on the choice of $J^\pl$ up to isomorphism.
\enrem

For each $\beta=\sum_{i\in I}m_i\alpha_i\in
Q_+\seteq\sum_{i\in I}\Z_{\ge0}\alpha_i$ with $\hgt(\beta)\seteq\sum_{i\in I}m_i=n$,
we define
\eq\label{eq:Rbeta}
&&\R_{\beta}=e_{\beta}\R_{n}e_{\beta},\quad \RC_\beta=e^J(\beta)\RC_ne^J(\beta)
\quad\text{and}\quad
\RCP_{\beta}=e_{\beta}\RCP_{n}e_{\beta},
\eneq where
\eqn\ba{rlrl}
I^\beta&= \set{ (i_1,\ldots,i_{n})\in I^{n}}{\sum_{a=1}^{n}\alpha_{i_a}=\beta },
&\quad
e_{\beta}&=\sum_{\nu\in I^\beta}e(\nu),\\[2ex]
J^\beta&= \set{ (j_1,\ldots,j_{n})\in J^{n}}%
{\sum_{a=1}^{n}\alpha_{\pr(j_a)}=\beta },&\quad
e^J_{\beta}&=\sum_{\nu\in J^\beta}e(\nu).\ea
\eneqn
It is easy to see that $e_\beta$ is a central even idempotent and
$\R_n = \soplus_{\hgt(\beta)=n} \R_{\beta}$.

\Th\label{MEQRR} For each $\beta=\sum_{i\in I}m_i\alpha_i\in Q_+$, we have a
$\cor$-superalgebra isomorphism
\begin{align*}
\R_{\beta}\otimes \clif_{\ell}\isoto \RCP_{\beta},
\end{align*}
where $\ell=\sum_{i\in I_\odd}m_i$.
In particular, if $\cor$ contains $\sqrt{-1}$, then $\R_\beta$ and
$\RC_\beta$ are weakly Morita superequivalent.
\enth
 For the notion of weak Morita superequivalence,
 see \S~\ref{weak superequivalence}.
\begin{proof}
To prove our theorem, we use the following strategy.

\bna
\item We first construct the elements $\{x_p\}_{1\leq p\leq n}$, $\{\tau_a\}_{1\leq a<n}$, $\{e(\nu)\}_{\nu\in I^\beta}$
in $\RCP_{\beta}$ which satisfy the defining relations for
$\R_{\beta}$.
\item We then construct the elements $\{C_i\}_{1\leq i\leq \ell}$
in $\RCP_{\beta}$ satisfying the relations
$C_iC_j+C_jC_i=2\delta_{i,j}$ $(1 \le i, j \le \ell)$ for the
Clifford superalgebra.
\item  We show that the $C_i$'s supercommute with the $x_pe(\nu)$'s and the
$\sigma_ae(\nu)$'s.
\item Finally, we prove that the resulting superalgebra homomorphism $\R_{\beta}\otimes \clif_{\ell}\isoto \RCP_{\beta}$
is an isomorphism.
\ee

\bigskip
\noi
(a)\quad
Let us choose $\gamma_{i,j}\in\cor^\times$ ($i,j\in I$)
such that
\eq&&
\parbox{65ex}{\be[{\rm (i)}]
\item $\gamma_{i,j}=1$ if $i\in I_\even$ or $j\in I_\even$,\label{cond:gamma a}
\item $\gamma_{i,j}\gamma_{j,i}=-1/2$ if $i,j\in \IODD$ and $i\neq j$,\label{cond:gamma b}
\item $\gamma_{i,i}=1/2$ if $i\in \IODD$.\label{cond:gamma c}
\ee}\label{cond:gamma} \eneq
It is obvious that such $\bl(\gamma_{i,j}\br)_{i,j\in I}$ exists. We
define \eqn x_pe(\nu) &=&
\begin{cases}
\cc_py_pe(\nu) & \textrm{if $\nu_p\in \IODD$}, \\
y_pe(\nu) & \textrm{if $\nu_p\in \IEVEN$}, \\
\end{cases}
\\[1ex]
\tau_ae(\nu) &=&
\gamma_{\nu_a,\nu_{a+1}}(\cc_a-\cc_{a+1})^{\prt{\nu_a}\prt{\nu_{a+1}}}\sigma_ae(\nu).
\eneqn
Then we can easily check the commutation relations \eqn
&&x_ax_be(\nu)=(-1)^{\prt{\nu_a}\prt{\nu_b}}x_bx_ae(\nu)\hs{5.5ex}\text{for $a\ne b$,}\\
&&\tau_ax_pe(\nu)=(-1)^{\prt{\nu_p}\prt{\nu_a}\prt{\nu_{a+1}}}
x_p\tau_ae(\nu)\quad\text{for $p\ne a, a+1$,}\\
&&(\tau_ax_{a+1}-x_a\tau_a)e(\nu)=(x_{a+1}\tau_a-\tau_ax_a)e(\nu)
=\delta_{\nu_a,\nu_{a+1}}e(\nu)\quad\\
&&\hs{34ex}\text{if $\nu_a\in\IEVEN$ or $\nu_{a+1}\in \IEVEN$.}
\eneqn Here, we have used (\ref{cond:gamma}\;\ref{cond:gamma a}). If
$\nu_a\in\IODD$ and $\nu_{a+1}\in \IODD$, then we also have the
relations
$$(\tau_ax_{a+1}+x_a\tau_a)e(\nu)=(x_{a+1}\tau_a+\tau_ax_a)e(\nu)=\delta_{\nu_a,\nu_{a+1}}e(\nu).$$
Indeed, using  (\ref{cond:gamma}\;\ref{cond:gamma b}), we have
\begin{align*}
(\tau_ax_{a+1}+x_a\tau_a)e(\nu) &=
\gamma_{\nu_a,\nu_{a+1}}((\cc_a-\cc_{a+1})\sigma_a\cc_{a+1}y_{a+1}+\cc_ay_a(\cc_a-\cc_{a+1})\sigma_a)e(\nu) \\
&=\gamma_{\nu_a,\nu_{a+1}}
(1+\cc_a\cc_{a+1})(\sigma_ay_{a+1}-y_a\sigma_a)e(\nu) \\
&=\gamma_{\nu_a,\nu_{a+1}}
\delta_{\nu_a,\nu_{a+1}}(1+\cc_a\cc_{a+1})(1-\cc_a\cc_{a+1})e(\nu) \\
&=\delta_{\nu_a,\nu_{a+1}}e(\nu),
\end{align*}
and
\begin{align*}
(x_{a+1}\tau_a+\tau_ax_a)e(\nu) &=
\gamma_{\nu_a,\nu_{a+1}}(\cc_{a+1}y_{a+1}(\cc_a-\cc_{a+1})\sigma_a+(\cc_a-\cc_{a+1})\sigma_a\cc_ay_a)e(\nu) \\
&=\gamma_{\nu_a,\nu_{a+1}}
(1-\cc_a\cc_{a+1})(y_{a+1}\sigma_a-\sigma_ay_a)e(\nu) \\
&=\gamma_{\nu_a,\nu_{a+1}}
\delta_{\nu_a,\nu_{a+1}}(1-\cc_a\cc_{a+1})(1+\cc_a\cc_{a+1})e(\nu) \\
&=\delta_{\nu_a,\nu_{a+1}}e(\nu)
\end{align*}
as expected.

Let us show the relation
\eq&&\tau^2_ae(\nu)=Q_{\nu_a,\nu_{a+1}}(x_a,x_{a+1})e(\nu).
\label{eq:tausquare}
\eneq
Note that
\eq
&&y_k^2e(\nu)=-x_k^2e(\nu)\quad\text{if $\nu_k\in\IODD$,}
\eneq
implies
$$\tQ_{\nu_a,\nu_{a+1}}(y_a,y_{a+1})e(\nu)=Q_{\nu_a,\nu_{a+1}}(x_a,x_{a+1})e(\nu).$$
Hence \eqref{eq:tausquare} is obvious unless $\nu_a,\nu_{a+1}\in
\IODD$. If $\nu_a,\nu_{a+1}\in \IODD$, by
(\ref{cond:gamma}\;\ref{cond:gamma b}), we have
$$\tau^2_ae(\nu)=-\gamma_{\nu_a,\nu_{a+1}}\gamma_{\nu_{a+1},\nu_{a}}(\cc_a-\cc_{a+1})^2\sigma_a^2e(\nu)
=\tQ_{\nu_a,\nu_{a+1}}(y_a,y_{a+1})e(\nu).$$

\medskip

Let us show the commutation relation \eqref{eq:ttt} in
Definition~\ref{DefRKLR}. We first note that \eq
&&(\cc_a-\cc_{b})\cc_ce(\nu)=-\cc_{s_{ab}(c)}(\cc_a-\cc_{b})e(\nu)\quad\text{if
$a\not=b$ and $\nu_a,\nu_b\in \IODD$,}
\label{eq:cperm}\\
&&(\cc_a-\cc_{b})^2e(\nu)=(\prt{\nu_a}+\prt{\nu_b})e(\nu)\quad\text{if
$a\not=b$,} \eneq where $s_{ab}$ is the transposition of $a$ and $b$.

Set
$$\Delta\seteq(\tau_{a+1}\tau_a\tau_{a+1}-\tau_{a}\tau_{a+1}\tau_{a})e(\nu).$$
Then we have
\begin{align*}
\tau_{a+1}\tau_a\tau_{a+1}e(\nu) &=
\gamma_{\nu_{a},\nu_{a+1}}(\cc_{a+1}-\cc_{a+2})^{\prt{\nu_a}\prt{\nu_{a+1}}}\sigma_{a+1}\\
&\hs{10ex}\gamma_{\nu_{a},\nu_{a+2}}(\cc_a-\cc_{a+1})^{\prt{\nu_a}\prt{\nu_{a+2}}}\sigma_a\\
&\hs{15ex}\gamma_{\nu_{a+1},\nu_{a+2}}(\cc_{a+1}-\cc_{a+2})^{\prt{\nu_{a+1}}\prt{\nu_{a+2}}}\sigma_{a+1}e(\nu) \\
&=\gamma_{\nu_{a},\nu_{a+1}}\gamma_{\nu_{a},\nu_{a+2}}\gamma_{\nu_{a+1},\nu_{a+2}}\\
&\hs{10ex}(\cc_{a+1}-\cc_{a+2})^{\prt{\nu_a}\prt{\nu_{a+1}}}(\cc_a-\cc_{a+2})^{\prt{\nu_a}\prt{\nu_{a+2}}}\\
&\hs{15ex}(\cc_{a}-\cc_{a+1})^{\prt{\nu_{a+1}}\prt{\nu_{a+2}}}
\sigma_{a+1}\sigma_a\sigma_{a+1}e(\nu),
\end{align*}
and
\begin{align*}
\tau_{a}\tau_{a+1}\tau_{a}e(\nu) &=
\gamma_{\nu_{a+1},\nu_{a+2}}(\cc_a-\cc_{a+1})^{\prt{\nu_{a+1}}\prt{\nu_{a+2}}}\sigma_{a}\\
&\hs{10ex}\gamma_{\nu_{a},\nu_{a+2}}(\cc_{a+1}-\cc_{a+2})^{\prt{\nu_{a}}\prt{\nu_{a+2}}}\sigma_{a+1}\\
&\hs{15ex}\gamma_{\nu_{a},\nu_{a+1}}(\cc_a-\cc_{a+1})^{\prt{\nu_{a}}\prt{\nu_{a+1}}}\sigma_{a}e(\nu) \\
&=\gamma_{\nu_{a},\nu_{a+1}}\gamma_{\nu_{a},\nu_{a+2}}\gamma_{\nu_{a+1},\nu_{a+2}}\\
&\hs{10ex}(\cc_{a}-\cc_{a+1})^{\prt{\nu_{a+1}}\prt{\nu_{a+2}}}
(\cc_a-\cc_{a+2})^{\prt{\nu_{a}}\prt{\nu_{a+2}}}\\
&\hs{15ex}(\cc_{a+1}-\cc_{a+2})^{\prt{\nu_{a}}\prt{\nu_{a+1}}}\sigma_{a}\sigma_{a+1}\sigma_{a}e(\nu). \\
\end{align*}

We can check easily
\eqn&&
(\cc_{a+1}-\cc_{a+2})^{s_2}(\cc_a-\cc_{a+2})^{s_0}(\cc_{a}-\cc_{a+1})^{s_1}e(s_as_{a+1}s_a\nu)\\
&&\hs{15ex}=(\cc_{a}-\cc_{a+1})^{s_1}(\cc_a-\cc_{a+2})^{s_0}(\cc_{a+1}-\cc_{a+2})^{s_2}e(s_as_{a+1}s_a\nu)
\eneqn in the case $s_0=s_1s_2=0$ and the case $s_0-1=s_1-s_2=0$.
Indeed,  by \eqref{eq:cperm}, we have
\eq
&&\ba{l}
(\cc_{a+1}-\cc_{a+2})^{s}(\cc_a-\cc_{a+2})(\cc_{a}-\cc_{a+1})^{s}e(\nu)\\
\hs{5ex}=(\cc_{a+1}-\cc_{a+2})^{s}(-\cc_{a+2}+\cc_{a+1})^{s}(\cc_a-\cc_{a+2})e(\nu)\\
\hs{5ex}=(1+\prt{\nu_{a+1}})^s(\cc_a-\cc_{a+2})e(\nu).\label{eq:ccc}
\ea \eneq

\medskip
Hence we have
\eqn
\Delta&=&\gamma_{\nu_{a},\nu_{a+1}}\gamma_{\nu_{a},\nu_{a+2}}\gamma_{\nu_{a+1},\nu_{a+2}}\\
&&\hs{10ex}(\cc_{a+1}-\cc_{a+2})^{\prt{\nu_a}\prt{\nu_{a+1}}}(\cc_a-\cc_{a+2})^{\prt{\nu_a}\prt{\nu_{a+2}}}\\
&&\hs{15ex}(\cc_{a}-\cc_{a+1})^{\prt{\nu_{a+1}}\prt{\nu_{a+2}}}
(\sigma_{a+1}\sigma_a\sigma_{a+1}-\sigma_{a}\sigma_{a+1}\sigma_{a})e(\nu),
\eneqn
which vanishes when $\nu_a\not=\nu_{a+2}$.

Assume $\nu_a=\nu_{a+2}=i$ and $\nu_{a+1}=j$. Then \eqn
\Delta&=&\gamma_{i,i}\gamma_{i,j}\gamma_{j,i}
(\cc_{a+1}-\cc_{a+2})^{\prt{i}\prt{j}}(\cc_a-\cc_{a+2})^{\prt{i}}\\
&&\hs{15ex}(\cc_{a}-\cc_{a+1})^{\prt{i}\prt{j}}
(\sigma_{a+1}\sigma_a\sigma_{a+1}-\sigma_{a}\sigma_{a+1}\sigma_{a})e(\nu).
\eneqn If $i\in I_\even$, we have
\eqn
&& \Delta=\dfrac{\tQ_{\nu_a,\nu_{a+1}}(y_{a+2},y_{a+1})-\tQ_{\nu_a,\nu_{a+1}}(y_{a},y_{a+1})}{y_{a+2}-y_a}e(\nu)\\
&&\hs{15ex}=\dfrac{Q_{\nu_a,\nu_{a+1}}(x_{a+2},x_{a+1})-Q_{\nu_a,\nu_{a+1}}(x_{a},x_{a+1})}{x_{a+2}-x_a}e(\nu).
\eneqn
If $i\in I_\odd$, then
\eqn
\Delta&=&\frac{1}{2}(-\frac{1}{2})^{\prt{j}}
(\cc_{a+1}-\cc_{a+2})^{\prt{j}}(\cc_a-\cc_{a+2})
(\cc_{a}-\cc_{a+1})^{\prt{j}}\\
&&\hs{15ex}(\sigma_{a+1}\sigma_a\sigma_{a+1}-\sigma_{a}\sigma_{a+1}\sigma_{a})e(\nu).
\eneqn
By \eqref{eq:ccc}, we have
$
(\cc_{a+1}-\cc_{a+2})^{\prt{j}}(\cc_a-\cc_{a+2})
(\cc_{a}-\cc_{a+1})^{\prt{j}}
=2^{\prt{j}}(\cc_a-\cc_{a+2})
$.
Hence
\begin{align*}
\Delta &= \frac{(-1)^{\prt{j}}}{2}
(\cc_a-\cc_{a+2})\Bigl(\frac{\tQ_{\nu_a,\nu_{a+1}}(y_{a+2},y_{a+1})
-\tQ_{\nu_a,\nu_{a+1}}(y_{a},y_{a+1})}{y_{a+2}-y_{a}} \\
&\hs{28ex}
-\cc_a\cc_{a+2}\frac{\tQ_{\nu_a,\nu_{a+1}}(y_{a+2},y_{a+1})
-\tQ_{\nu_a,\nu_{a+1}}(y_{a},y_{a+1})}{y_{a+2}+y_{a}}\Bigr)e(\nu)\\
&= \frac{(-1)^{\prt{j}}}{2}\bl((\cc_a-\cc_{a+2})(y_{a+2}+y_{a})
-(\cc_a-\cc_{a+2})\cc_a\cc_{a+2}(y_{a+2}-y_{a})\br)\\
&\hs{20ex}\dfrac{\tQ_{\nu_a,\nu_{a+1}}(y_{a+2},y_{a+1})-
\tQ_{\nu_a,\nu_{a+1}}(y_{a},y_{a+1})}{y_{a+2}^2-y_a^2}e(\nu).
\end{align*}

Since
$$
\begin{aligned}
&
(\cc_a-\cc_{a+2})(y_{a+2}+y_{a})-(\cc_a-\cc_{a+2})\cc_a\cc_{a+2}(y_{a+2}-y_{a})
\\
& \hs{5ex}=2\bl(\cc_a y_a -\cc_{a+2}y_{a+2}\br) =2(x_a-x_{a+2}),
\end{aligned}
$$
we finally obtain
$$\Delta=(-1)^{\prt{j}}(x_a-x_{a+2})\dfrac{Q_{\nu_a,\nu_{a+1}}(x_{a+2},x_{a+1})
-Q_{\nu_a,\nu_{a+1}}(x_{a},x_{a+1})}{-x_{a+2}^2+x_a^2}e(\nu).$$

\smallskip
The other commutation relations are obvious.

\bigskip
\noi (b)\quad For $\nu\in I^\beta$, we define $1\leq
p_{1}(\nu)<\cdots <p_{\ell}(\nu)\leq n$ as a unique sequence of
integers such that $\set{a}{\text{$1\le a\le n$ and $\nu_a\in
I_\odd$}}=\set{p_{k}(\nu)}{1\le k\le \ell}$. Set $C_k=\sum_{\nu\in
I^\beta}\cc_{p_{k}(\nu)}e(\nu)$. It is easy to see that
$C_kC_j+C_jC_k=2\delta_{kj}$ and $C_ke(\nu)=e(\nu)C_k$.

\medskip
\noi (c)\ Let us show $C_kx_ae(\nu)=(-1)^{\prt{\nu_a}}x_ae(\nu)C_k$.
Since it is obvious if $\nu_a\in I_\even$, we assume that $\nu_a\in
I_\odd$. By the definition, we have
$C_kx_ae(\nu)=\cc_{p_k(\nu)}\cc_ay_ae(\nu)$. If $p_k(\nu)\not=a$,
then $\cc_{p_k(\nu)}$ anticommutes with $ \cc_a$ and commutes with
$y_ae(\nu)$. If $p_k(\nu)=a$, then $\cc_{p_k(\nu)}$ commutes with $
\cc_a$ and anticommutes with $y_ae(\nu)$. Hence in both cases,
$\cc_{p_k(\nu)}$ anticommutes with $x_ae(\nu)=\cc_ay_ae(\nu)$.

Let us show
\eq
&&C_k\tau_ae(\nu)=(-1)^{\prt{\nu_a}\prt{\nu_{a+1}}}\tau_ae(\nu)C_k.
\label{eq:ctau}
\eneq
It is obvious if $\nu_a\in I_\even$ or $\nu_{a+1}\in I_\even$.
Assume that $\nu_a,\nu_{a+1}\in I_\odd$.
Then we have
$$C_k\tau_ae(\nu)=C_ke(s_a\nu)\tau_ae(\nu)=\cc_{p_k(s_a\nu)}(\cc_a-\cc_{a+1})\gamma_{\nu_a,\nu_{a+1}}\sigma_ae(\nu).$$
Since $p_k(s_a\nu)=p_k(\nu)$, by \eqref{eq:cperm}, we have
$$\cc_{p_k(s_a\nu)}(\cc_a-\cc_{a+1})\sigma_a
=-(\cc_a-\cc_{a+1})\cc_{s_ap_k(\nu)}\sigma_a=-(\cc_a-\cc_{a+1})\sigma_a\cc_{p_k(\nu)},$$
which implies \eqref{eq:ctau}.

Now assume that one of $\nu_a$ or $\nu_{a+1}$ belongs to $I_\odd$
and the other belongs to $I_\even$. In this case
$p_k(s_a\nu)=s_a(p_k(\nu))$ holds, and hence we have
$$C_k\tau_ae(\nu)=\cc_{p_k(s_a\nu)}\sigma_ae(\nu)
=\cc_{s_a(p_k(\nu))}\sigma_ae(\nu)=\sigma_a\cc_{p_k(\nu)}e(\nu)=\sigma_ae(\nu)C_k.$$

Thus we have constructed a superalgebra homomorphism $\R_n\otimes\clif_\ell
\to\RCP_n$.

\medskip
\noi (d) It is obvious that $\R_n\otimes\clif_\ell$ is generated by
$$\set{
x^a\tau_{w}e(\nu)\otimes C^\eta}%
{a\in\Z^n_{\geq 0},\,w\in\Sym_n, \,\nu\in I^n,\,\eta\in (\FF)^\ell}$$
as a $\cor$-module.
Its image by the homomorphism $\R_n\otimes\clif_\ell\to\RCP_n$
forms a basis of $\RCP_n$ as a $\cor$-module by \eqref{eq:pbwCl}.
Hence $\R_n\otimes\clif_\ell\to\RCP_n$ is bijective.
\end{proof}

\Rem\label{rem:QtQ} We have constructed $I$ and
$Q=\bl(Q_{i,j}(w,z)\br)_{i,j\in I}$ starting from $J$ and
$\bl(\tQ_{i,j}(u,v)\br)_{i,j\in J}$. Conversely, we can construct
$J$ and $\bl(\tQ_{i,j}(u,v)\br)_{i,j\in J}$ starting from $I$ and
$Q=\bl(Q_{i,j}(w,z)\br)_{i,j\in I}$. Assume that $I=I_\even\sqcup
I_\odd$ and $Q=\bl(Q_{i,j}(w,z)\br)_{i,j\in I}$ are given so that
\eqref{eq:Q} is satisfied. Then set $J=(I_\even\times \{0,1\})\sqcup
(I_\odd\times \{0\})$. Let $\pr\cl J\to I$ be the canonical map, and
let $c$  be the involution given by $c(i,\eps)=(i,1-\eps)$ for $i\in
I_\even$ and $c(i,0)=(i,0)$ for $i\in I_\odd$. For $i,j\in I$, let
$\tAt\subset\At$ and $\tBt\subset \Bt$ be the algebras as in
Definition~\ref{def:At}, and let $\tAt\simeq\tBt$ be the isomorphism
given in \eqref{eq:AB}. Let $Q'_{i,j}\in\tBt$ be the element
corresponding to $Q_{i,j}\in\tAt$. Define
$$\tQ_{(i,\eps),(i',\eps')}(u,v)=Q'_{i,i'}((-1)^{\eps}u,(-1)^{\eps'}v)
\ \ \text{for} \  (i,\eps), (i',\eps')\in J.$$ Then
$(\tQ_{j,j'})_{j,j'\in J}$ satisfies the condition \eqref{cond:Q}.
\enrem

\medskip
Combining Theorem~\ref{MEQRR} with Remark~\ref{rem:QtQ} and
Corollary~\ref{cor:pbwrc}, we immediately obtain the following
corollary.

\Cor\label{cor:pbw} Let $I=I_\even\sqcup I_\odd$ and let
$Q=\bl(Q_{i,j}(w,z)\br)_{i,j\in I}$ be a family of skew polynomials
satisfying the conditions \eqref{eq:Q}. Let $\R_n(Q)$ be the
associated quiver Hecke $\cor$-superalgebra. For each $w\in \Sym_n$,
we choose a reduced expression $s_{i_1}\cdots s_{i_{\ell}}$ of $w$
and write $\tau_w=\tau_{i_1}\cdots \tau_{i_{\ell}}$. Then
\begin{align*}
\set{
x^a\tau_{w}e(\nu)}%
{a\in\Z^n_{\geq 0},\,w\in\Sym_n,\,\nu\in I^n}
\end{align*}
is a basis of $\R_n(Q)$.
\encor

\bigskip

\subsection{Quiver Hecke superalgebra associated with generalized Cartan matrices}\label{specRC}

Let $\IRED$ be a finite index set with a decomposition $\IRED =
\IODD \sqcup \IEVEN$ and let $A=(a_{i,j})_{i,j\in \IRED}$ be a
symmetrizable \GCM\ indexed by $\IRED$.
We assume that
\eq &&\text{$a_{i,j}\in 2\Z$ if $i\in \IODD$.} \eneq
For $i,j\in I$, let $S_{i,j}$ be the set of $(r,s)$ where $r$ and
$s$ are integers satisfying the following conditions: \eq&&\left\{
\parbox{70ex}{\bnum
\item $0\leq r\le -a_{i,j}$, $0\le s\le -a_{j,i}$,
\item $a_{j,i}r+a_{i,j}s=-a_{i,j}a_{j,i}$,
\item $r\in2\Z$ if $i\in I_\odd$,
\item $s\in2\Z$ if $j\in I_\odd$.
\ee} \right. \eneq Let $\{t_{i,j;r,s}\}_{i\ne j,\;(r,s)\in S_{i,j}}$
be  a family of indeterminates such that $t_{i,j;r,s}=t_{j,i;s,r}$
and $t_{i,j;-a_{i,j},0}$ is invertible, and let $\cor_A$ be the
algebra $\Z[1/2][\{t_{i,j;r,s}\}][\{(t_{i,j,-a_{i,j},0})^{-1}\}]$. We
take $Q=(Q_{i,j})_{i,j\in \IRED}$, where
$$Q_{i,j}(w,z)=\sum_{(r,s)\in S_{i,j}}t_{i,j;r,s}w^rz^s
\in\cor_{A} \langle w,z \rangle/\langle
zw-(-1)^{\prt{i}\,\prt{j}}wz\rangle$$ for $i\not=j$
and $Q_{i,j}=0$ for $i=j$. We write $\R_n=\R_n(A)=\R_n(Q)$ for the
associated quiver Hecke $\cor_A$-superalgebra.

Let $J=I\times \{0\}\sqcup I_\even\times\{1\}$ and
$\tQ=(\tQ_{i,j})_{i,j\in J}$ associated with $I$ and $Q$ as in
Remark~\ref{rem:QtQ}. Then we can define the quiver Hecke-Clifford
$\cor_A$-superalgebra $\RC_n=\RC_n(A)\seteq \RC_n(\QTL)$.
Note that $\R_n$ and $\RC_n$ are $(\Z\times
(\FF))$-graded via the following assignment. \eqn&& \ba{l} \DEG
e(\nu) = (0;0),\quad
\DEG x_pe(\nu)=((\alpha_{\nu_p},\alpha_{\nu_p});\prt{\nu_a}),\\
\DEG\tau_ae(\nu)=(-(\alpha_{\nu_a},\alpha_{\nu_{a+1}});\prt{\nu_a}\prt{\nu_{a+1}}),\\[1ex]
\text{and}\\[1ex]
\DEG y_pe(\nu)=((\alpha_{\nu_p},\alpha_{\nu_p});{0}),\quad
\DEG \cc_pe(\nu) = (0;{1}),\quad \\
\DEG\sigma_ae(\nu)=(-(\alpha_{\nu_a},\alpha_{\nu_{a+1}});{0}).
\ea
\eneqn

If $\cor_A\to \cor$ is a ring homomorphism, we can consider the
superalgebra $\cor\otimes_{\cor_A}\R_n$, which will be called the
quiver Hecke superalgebra associated with a generalized Cartan
matrix $A$.

\subsection{Cyclotomic quiver Hecke superalgebras}\label{sec:cyclq}
Let $Q^\vee\seteq\oplus_{i\in I}\Z h_i$ be the coroot lattice with
the bilinear form $Q^\vee\times Q\to \Z$ by $\lan h_i,
\alpha_j\ran=a_{i,j}$. Let $P\seteq\Hom(Q^\vee,\Z)$ be the weight
lattice, and $\bl(\Lambda_i\br)_{i\in I}$ the dual basis of
$\bl(h_i\br)_{i\in I}$. Set $$P_+\seteq\set{\la\in
P}{\text{$\la(h_i)\ge0$ for any $i\in I$}}.$$ For $\Lambda\in P_+$,
set $x_1^\Lambda=\sum_{\nu\in I^n} x_1^{\Lambda(h_{\nu_1})}e(\nu)\in
\R_n$. We define
$$\R_n^\Lambda=\R_n/\R_nx_1^\Lambda\R_n =\R_n/\bl(\ssum_{\nu\in
I^n}x_1^{\Lambda(h_{\nu_1})}e(\nu)\R_n\br),$$
$$R_\beta^\Lambda=\R_\beta/\R_\beta x_1^\Lambda\R_\beta
\quad\text{for $\beta\in Q_+$},$$ and call them the {\em cyclotomic
quiver Hecke superalgebras}. Similarly, we can define the {\em
cyclotomic quiver Hecke-Clifford superalgebras} $\RC_n^\Lambda$ and
$\RC_\beta^\Lambda$. By Theorem~\ref{MEQRR} the superalgebras $\R_\beta^\Lambda$ and
$\RC_\beta^\Lambda$ are weakly Morita superequivalent. Note that
$\RC_n^\Lambda$ and $\R_n^\Lambda$ are finitely generated
$\cor$-modules.

\subsection{Completion of $\RC_n$}
For the data $J$, $c$, $\bl(\tQ_{i,j}\br)_{i,j\in J}$ as in
\S\,\ref{RC_subsection}, let $\RC_n$ be the associated quiver
Hecke-Clifford superalgebra. Let $\ga$ be the ideal of
$\cor[y_1,\ldots,y_n]$ generated by $y_1,\ldots,y_n$. Then we have
\eqn &&\text{for any $s\in\RC_n$, there exists $m$ such that
$s\ga^{k+m}\subset\ga^k\RC_n$.} \eneqn Indeed, it is enough to
verify this for the generators, and it is obvious for $y_p$ and
$\cc_p$. The case for $s=\sigma_a$ follows from
Lemma~\ref{lem:sigmacom}.

Hence we can  easily see that the superalgebra structure of $\RC_n$
induces a superalgebra structure
on $\prolim[k]\RC_n/\ga^k\RC_n$.
\Def\label{compRC}
 We define the completion $\HRC$ of $\RC_n$ to be
 $$\HRC=\prolim[k]\RC_n/\ga^k\RC_n\simeq\cor[[y_1,\ldots,y_n]]\
\otimes_{\cor[y_1,\ldots,y_n]}\RC_n.$$
 \edf
 Then the  formula \eqref{com:sigmf} holds in $\HRC$
for all $f\in\cor[[y_1,\ldots,y_n]] $. The algebra $\HRC$ contains
$\soplus_{\nu\in J}\cor[[y_1,\ldots,y_n]]e(\nu)$ as a subalgebra.

\section{Relation to the affine Hecke-Clifford superalgebras}

In the rest of this article, $\cor$ is an algebraically closed field
of characteristic\;$\not=2$. We fix a non-zero element $q\in
\cor^\times$ and set $\xi=q-q^{-1}$. We assume \eq (q^2)^2\not=1.
\eneq The purpose of this section is to prove that the affine
Hecke-Clifford superalgebras are isomorphic to quiver Hecke-Clifford
superalgebras after a completion.

\subsection{Affine Hecke-Clifford superalgebra}
For a more comprehensive treatment of  the affine  Hecke-Clifford
superalgebra and its cyclotomic quotients, see ~\cite{BK3} and the
references therein.

\Def \label{def:AHC}Let $n\geq 0$ be an integer. The {\em affine
Hecke-Clifford superalgebra} $\AHC$ of degree $n$ is the
$\cor$-superalgebra generated by the even generators $X_1^{\pm1},
\ldots, X_n^{\pm1}, T_1,\ldots, T_{n-1}$ and the odd generators
$C_1,\ldots,C_n$ with the following relations:
\bnum
\item $X_{i}X_{i}^{-1}=X_{i}^{-1}X_i=1$, $X_iX_j=X_jX_i$ for all $1\leq i,j\leq n$,
\item $C_iC_j+C_jC_i=0$ $( i\ne j)$ and $C_i^2=1$,
\label{clifford_rel}
\item $T_i^2=\xi T_i+1$, $T_iT_j=T_jT_i$ $(|i-j|>1)$, $T_iT_{i+1}T_i=T_{i+1}T_iT_{i+1}$,
\label{iwahori_rel}
\item $T_iC_j=C_jT_i$ $(j\not=i,i+1)$, $T_iC_i=C_{i+1}T_i$ for all $1\leq i\leq n-1$,
\label{non_triv_eq}
\item  $C_iX_j=X_jC_i$ for all $i\ne j$ and $C_iX_i^{\pm1}=X_i^{\mp1}C_i$,\label{def:CX}
\item $ T_iX_j=X_jT_i$  $(j\ne i, i+1)$,
$(T_i+\xi C_iC_{i+1})X_iT_i=X_{i+1}$.
\end{enumerate}
\label{def_aff} \edf

It is known that $\AHC$ has a PBW-type basis.

\Prop\label{prop:pbw} For each $w\in \Sym_n$, let $w=s_{i_1}\cdots
s_{i_{\ell}}$ be a reduced expression and set $T_w =T_{i_1}\cdots
T_{i_{\ell}}$ \ro which is independent of the choice of a reduced
expression\rf.
Then $\AHC$ is a free left
$\cor[X_1^{\pm1},\ldots,X_n^{\pm1}]$-module with a basis
\begin{align*}
\set{
C^{\eta}T_{w}e(\nu)}%
{\eta\in(\FF)^n,\,w\in\Sym_n,\,\nu\in J^n},
\end{align*}
where $C^{\eta}=C_1^{\eta_1}\cdots C_n^{\eta_n}$ for
$\eta=(\eta_1,\dots,\eta_n)\in(\FF)^n$.
\enprop
\subsection{Intertwiners}\label{Intertwin}

We recall the definition of intertwiners of $\AHC$
(\cite{Naz},\cite[\S14.8]{Kle} ??):
\eq&&\Phi_a=T_a+\xi\dfrac{X_{a+1}}{X_a-X_{a+1}}+\xi C_aC_{a+1}
\dfrac{X_aX_{a+1}}{X_aX_{a+1}-1}
\in\cor(X_1,\ldots,X_n)\tens_{\cor[X_1^{\pm1},\ldots,X_n^{\pm1}]}
\AHC.\eneq %
They satisfy the relations: %
\eq&&\Phi_aX_i=X_{s_a(i)}\Phi_a\quad\text{and}\quad
\Phi_aC_i=C_{s_a(i)}\Phi_a. \eneq %
Setting %
\eqn&& K(u,v)=u^2-(q^2+q^{-2})uv+v^2+4(q-q^{-1})^2,\quad\\[1ex]
&&F(X_a,X_{a+1})=\dfrac{X_{a}^2X_{a+1}^2\;K(X_a+X_{a}^{-1},\,X_{a+1}+X_{a+1}^{-1})}%
{(X_a-X_{a+1})^2(X_aX_{a+1}-1)^2} ,\eneqn%
we have \eq &&\Phi_a^2=F(X_a,X_{a+1}).\eneq%
 Note that by setting $w_k=X_k+X_k^{-1}=2\dfrac{\la_k+\la_k^{-1}}{q+q^{-1}}$ ($k=1,2$),
we have
\eq&&
\ba{rl}
K(w_1,w_2)&= \dfrac{4}{(q+q^{-1})^2\lambda_1^2\la_2^2}(\lambda_1-q^2\la_2)
(\lambda_1-q^{-2}\la_2)(\lambda_1\la_2-q^2)(\lambda_1\la_2-q^{-2}),\\[4ex]
F(X_1,X_2)&=\dfrac%
{(\la_2-q^2\la_1)(\la_2-q^{-2}\la_1)(\la_2-q^2\la_1^{-1})(\la_2-q^{-2}\la_1^{-1})}%
{(\la_2-\la_{1})^2(\la_2-\la_1^{-1})^2}. \ea \label{eq:KLam} \eneq
Moreover, the $\Phi_a$'s satisfy the braid relations:
$$\Phi_{a}\Phi_{c}=
\Phi_{c}\Phi_{a} \ \ (|a-c|>1), \ \ \ \ \Phi_{a}\Phi_{a+1}\Phi_{a} =
\Phi_{a+1}\Phi_{a}\Phi_{a+1}.$$

\medskip
\subsection{Localization of $\AHC$}
For a finite-dimensional representation of $\AHC$,
 we consider the simultaneous eigenspaces of the
$X_i$'s.
 Let us denote by $\AAA^1$ the one-dimensional affine space and by $\T$
 the algebraic torus $\AAA^1\setminus \{0\}$.
We denote by $\T_X$ the variety $\T$ with $X$ as a coordinate.
It plays a role of the set of eigenvalues of
$X_k$.
Let $c\cl \T_X\to\T_X$ be the involution of $\T_X$ given by
$c(X)=X^{-1}$ (corresponding to Definition \ref{def:AHC}
\eqref{def:CX}).

For an algebraic variety $X$ over $\cor$ and its $\cor$-valued point
$x$, let us denote by ${{\CO}}_{X,x}$  the germ of the structure
sheaf $\CO_X$ at $x$, by ${\HAT{\CO}}_{X,x}$ its completion, and by
$\FRAC({\HAT{\CO}}_{X,x})$ its fraction field. For $p$ such that
$1\le p\le n$, let us denote by $c_p\cl \T^n\to\T^n$ the involution
$$(X_1,\ldots,X_n)\mapsto
(X_1,\ldots,X_{p-1},X_p^{-1},X_{p+1},\ldots,X_n),$$ and by the same
letter $c_p$ the induced isomorphism
$\FRAC({\HAT{\CO}}_{{\T}^n,q})\isoto\FRAC({\HAT{\CO}}_{{\T}^n,c_p(q)})$
($q\in \T^n$). For $1\le p<n$, we denote by $s_p\cl \T^n\to\T^n$ the
involution
$$(X_1,\ldots,X_n)\mapsto
(X_1,\ldots,X_{p-1},X_{p+1},X_{p},X_{p+2},\ldots,X_n),$$ and by
the same letter $s_p$ the induced isomorphism
$\FRAC({\HAT{\CO}}_{{\T}^n,q})\isoto\FRAC({\HAT{\CO}}_{{\T}^n,s_p(q)})$.
Similarly, we denote by $\bs_p\cl \T^n\to\T^n$ the
involution
$$(X_1,\ldots, X_n)\mapsto
(X_1,\ldots,X_{p-1},X_{p+1}^{-1},X_{p}^{-1},X_{p+2},\ldots,X_n),$$
and the induced isomorphism
$\FRAC({\HAT{\CO}}_{{\T}^n,q})\isoto\FRAC({\HAT{\CO}}_{{\T}^n,\bs_p(q)})$.
Clearly, we have $\bs_p=c_pc_{p+1}s_p$.

\medskip
Let us denote by $\HC$ the $\cor$-superalgebra generated by $C_p$
($1\le p\le n$), $T_a$ ($1\le a<n$)
 with  \eqref{clifford_rel}--\eqref{non_triv_eq} in
 Definition~\ref{def:AHC} as the  defining relations.
 The superalgebra $\HC$ is called the {\em Hecke-Clifford superalgebra}
 and can be regarded as a subsuperalgebra of $\AHC$.

 \Def\label{compHC} Let $J$ be a finite subset of $\T_X$ invariant under $c$, and
 let $X\cl J\to\T_X$ be the inclusion map.
 We define the commutative $\cor$-algebras
\begin{align*}
\CO_{n} = \bigoplus_{\nu\in J^n}
{\HAT{\CO}}_{\T^n,X(\nu)}e(\nu),\quad \CK_n = \bigoplus_{\nu\in J^n}
\FRAC({\HAT{\CO}}_{{\T}^n,X(\nu)})e(\nu),
\end{align*}
where $X(\nu)=(X(\nu_1),\dots,X(\nu_n))\in \T^n$ for
$\nu=(\nu_1,\ldots,\nu_n)$. We define the algebra structure on
\begin{align*}
\KHC=\CK_n\otimes_{\cor}\HC
\end{align*}
by
\begin{equation}\label{def:tformal}
\begin{aligned}
C_pe(\nu)f & =e(\CC_p\nu)c_p(f)C_p  \ \ (1\leq p\leq n), \\
T_ae(\nu)f  & = e(s_a\nu)s_a(f)T_a -\xi
\frac{e(\nu)f-e(s_a\nu)s_af}{X_aX_{a+1}^{-1}-1}  -\xi
C_aC_{a+1}\dfrac{e(\nu)f-e(\bs_a\nu)\bs_a(f)}{X_a^{-1}X_{a+1}^{-1}-1}
\end{aligned}
\end{equation}
for $\nu\in J^n$, $1\leq a<n$, $f\in
\FRAC({\HAT{\CO}}_{{\T}^n,X(\nu)})$.

We define $\OHC$ to be the subsuperalgebra $\CO_n\otimes \HC$ of
$\KHC$.
\edf

Let us denote by $\ga_J$ the ideal of $\cor[X_1^{\pm1},\ldots,X_n^{\pm1}]$
consisting of functions which vanish on $\set{X(\nu)}{\nu\in J^n}$.
Then we have
\eqn
\CO_n&\simeq&\prolim[k]\cor[X_1^{\pm1},\ldots,X_n^{\pm1}]/\ga_J^k,\\
\OHC&\simeq&\prolim[k]\AHC/\ga_J^k\AHC,\\
\KHC&\simeq &\CK_n\otimes_{\CO_n}\OHC. \eneqn
 Thus there exists a
$\cor$-superalgebra homomorphism $\AHC\to\OHC$, which is injective
if $J$ is not empty.

\medskip
Let $\MOD{\AHC}$ be the category of $\AHC$-modules
finite-dimensional over $\cor$. Let $\mathcal{C}(J)$ be the full
subcategory of $\MOD{\AHC}$ consisting of $M\in \MOD{\AHC}$ whose
$X_k$-eigenvalues lie in $\set{X(i)}{i\in J}$ for all $1\leq k\leq
n$. For such an $M$ and $\nu\in J^n$, let $e(\nu)$ be the projection
operator of $M$ onto the simultaneous generalized eigenspace of
$\bl(X_k\br)_{1\le k\le n}$ with $\bl(X(\nu_k)\br)_{1\le k\le n}$ as
eigenvalues. Then $M$ has a natural structure of $\OHC$-module. Thus
$\mathcal{C}(J)$ is equivalent to the category $\MOD{\OHC}$ of
$\OHC$-modules finite-dimensional over $\cor$.
\medskip
We have
\eq
&&\text{$\Phi_a e(\nu)=e(s_a\nu)\Phi_a$ for any $1\le a<n$.}
\label{eq:comphi}
\eneq
Indeed, this formula must hold since $\Phi_a$ is an intertwiner,
and the definition \eqref{def:tformal} is
equivalent to $\Phi_a e(\nu)f=e(s_a\nu)s_a(f)\Phi_a$ for $f\in
\FRAC({\HAT{\CO}}_{{\T}^n,X(\nu)})$.
\bigskip
\subsection{Dynkin diagram}\label{sec:Dynahc}
Let us denote by $\AAA^1_w$ the one-dimensional affine space $\AAA^1$ with
the coordinate $w$. Let $\pr\cl \T_X\to\AAA^1_w$ be the map
$\pr(X)=X+X^{-1}$.
Hence the fiber of $\pr$ is of the form $\{j,c(j)\}$ for some $j\in \T_X$.
Let $\T_\la$ be the variety $\T$ with the
coordinate $\la$. Let $g\cl\T_\la\to\AAA^1_w$ be the map
$g(\la)=2\dfrac{\la+\la^{-1}}{q+q^{-1}}$.
Thus these three varieties are related by
$$\xymatrix@C=6ex{{\T_X}\ar[r]^-{\pr}&{\AAA^1_w}
&{\T_\la}\ar[l]_-g.}$$

\medskip
We will regard $\AAA^1_w$ as the set of vertices of a Dynkin diagram. We consider
a $\Q$-vector space with $\{\alpha_w\}_{w\in\AAA^1_w}$ as a basis, and
define an inner product on it as follows:
\eq &&\left\{
\ba{l}
(\alpha_w,\alpha_w)=
\begin{cases}
1&\text{if $\#\pr^{-1}(w)=1$, }\\
4&\text{if $\#g^{-1}(w)=1$,}\\
2&\text{otherwise,}
\end{cases}\\[5ex]
\text{if $w_1\not=w_2$ and $K(w_1,w_2)=0$, then}\\[1ex]
\hs{3ex}(\alpha_{w_1},\alpha_{w_2})=
\begin{cases}
-2&\text{if $(\alpha_{w_1},\alpha_{w_1})=4$ or $(\alpha_{w_2},\alpha_{w_2})=4$,}\\[2ex]
-1&\text{if $(\alpha_{w_1},\alpha_{w_1})$, $(\alpha_{w_2},\alpha_{w_2})\not=4$,}
\end{cases}\\[4ex]
\text{if $w_1\not=w_2$ and $K(w_1,w_2)\not=0$, then}\\[2ex]
\hs{3ex}(\alpha_{w_1},\alpha_{w_2})=0.
\ea
\right.\label{def:inner}
\eneq

Note that by \eqref{eq:KLam}, $K(w_1,w_2)=0$ if and only if
$\la_1=q^{\pm2}\la_2$ or $\la_1=q^{\pm2}\la_2^{-1}$ when writing
$w_k=g(\la_k)=2\dfrac{\la_k+\la_k^{-1}}{q+q^{-1}}$ ($k=1,2$). The
generalized Cartan matrix $\bl(a_{w,w'}\br)_{w,w'\in\AAA^1_w}$ is
given by $a_{w,w'}=2(\alpha_w,\alpha_{w'})/(\alpha_w,\alpha_w)$.

Set $(\AAA^1_w)_\odd=\{\pr(\pm1)\}=\{g(\pm q)\}$ and
$(\AAA^1_w)_\even=\AAA^1_w\setminus\AAA^1_\odd$. Then $\AAA^1_w$
decomposes into the even part and the odd part.

The connected components of the Dynkin diagram $\AAA^1_w$ are
classified as follows. (The odd vertices are marked by $\times$.)

\bnum \item when $q^2$ is not a root of unity, there are three types of
connected components.

\vs{1ex}
\begin{enumerate}[(a)]
\item $\set{g(\zeta q^{2k})}{k\in\Z}$ for some
$\zeta\not\in \pm q^\Z$, where $\pm q^\Z=\set{\pm q^k}{k\in\Z}$.
The Dynkin diagram is of type $A_{\infty}$.

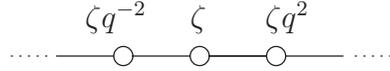
\vs{2ex}
\begin{center}
\unitlength 0.1in
\begin{picture}( 20.0000,  3.0000)(  0.0000, -2.5000)
%
{\color[named]{Black}{%
\special{pn 8}%
\special{ar 1000 200 50 50  0.0000000 6.2831853}%
}}%
%
{\color[named]{Black}{%
\special{pn 8}%
\special{ar 1400 200 50 50  0.0000000 6.2831853}%
}}%
%
{\color[named]{Black}{%
\special{pn 8}%
\special{ar 600 200 50 50  0.0000000 6.2831853}%
}}%
%
{\color[named]{Black}{%
\special{pn 8}%
\special{pa 650 200}%
\special{pa 950 200}%
\special{fp}%
}}%
%
{\color[named]{Black}{%
\special{pn 8}%
\special{pa 1050 200}%
\special{pa 1350 200}%
\special{fp}%
}}%
%
{\color[named]{Black}{%
\special{pn 8}%
\special{pa 1450 200}%
\special{pa 1750 200}%
\special{fp}%
}}%
%
{\color[named]{Black}{%
\special{pn 8}%
\special{pa 1050 200}%
\special{pa 1350 200}%
\special{fp}%
}}%
%
{\color[named]{Black}{%
\special{pn 8}%
\special{pa 250 200}%
\special{pa 550 200}%
\special{fp}%
}}%
\put(9.5000,-0.8000){\makebox(0,0)[lb]{$\zeta$}}%
\put(12.0000,-0.8000){\makebox(0,0)[lb]{\kern2ex$\zeta q^2$}}%
\put(4.0000,-0.8000){\makebox(0,0)[lb]{$\zeta q^{-2}$}}%
%
{\color[named]{Black}{%
\special{pn 8}%
\special{pa 200 200}%
\special{pa 0 200}%
\special{dt 0.045}%
}}%
%
{\color[named]{Black}{%
\special{pn 8}%
\special{pa 2000 200}%
\special{pa 1800 200}%
\special{dt 0.045}%
}}%
\end{picture}%

\end{center}

\vs{3ex}

\item $\set{g(\eps q^{2k+1})}{k\in\Z_{\ge 0}}$ ($\eps=\pm1$).
The Dynkin diagram is of type $B_\infty$.

\medskip
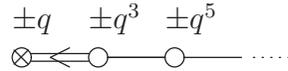
\begin{center}
\unitlength 0.1in
\begin{picture}( 14.5000,  3.0000)(  1.5000, -2.5000)

{\color[named]{Black}{%
\special{pn 8}%
\special{ar 600 200 50 50  0.0000000 6.2831853}%
}}%

{\color[named]{Black}{%
\special{pn 8}%
\special{ar 1000 200 50 50  0.0000000 6.2831853}%
}}%

{\color[named]{Black}{%
\special{pn 8}%
\special{ar 200 200 50 50  0.0000000 6.2831853}%
}}%

{\color[named]{Black}{%
\special{pn 8}%
\special{pa 650 200}%
\special{pa 950 200}%
\special{fp}%
}}%

{\color[named]{Black}{%
\special{pn 8}%
\special{pa 1050 200}%
\special{pa 1350 200}%
\special{fp}%
}}%

{\color[named]{Black}{%
\special{pn 8}%
\special{pa 650 200}%
\special{pa 950 200}%
\special{fp}%
}}%
\put(9.5000,-0.8000){\makebox(0,0)[lb]{$\pm q^5$}}%

{\color[named]{Black}{%
\special{pn 8}%
\special{pa 1600 200}%
\special{pa 1400 200}%
\special{dt 0.045}%
}}%
\put(5.5000,-0.8000){\makebox(0,0)[lb]{$\pm q^3$}}%
\put(1.5000,-0.8000){\makebox(0,0)[lb]{$ \pm q$}}%

{\color[named]{Black}{%
\special{pn 8}%
\special{pa 250 180}%
\special{pa 550 180}%
\special{fp}%
}}%
\put(2.0500,-2.55000){\makebox(0,0)[lb]{\kern-1ex$\times$}}
{\color[named]{Black}{%
\special{pn 8}%
\special{pa 250 220}%
\special{pa 550 220}%
\special{fp}%
}}%

{\color[named]{Black}{%
\special{pn 8}%
\special{pa 350 200}%
\special{pa 450 150}%
\special{fp}%
}}%

{\color[named]{Black}{%
\special{pn 8}%
\special{pa 450 250}%
\special{pa 350 200}%
\special{fp}%
}}%
\end{picture}%

\end{center}

\vs{3ex}
\item $\set{g(\eps q^{2k})}{k\in\Z_{\geq 0}}$ ($\eps=\pm1$).
The Dynkin diagram is of type $C_\infty$.

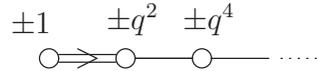
\vs{3ex}
\begin{center}
\unitlength 0.1in
\begin{picture}( 16.0000,  3.0000)(  0.0000, -2.5000)
%
{\color[named]{Black}{%
\special{pn 8}%
\special{ar 600 200 50 50  0.0000000 6.2831853}%
}}%
%
{\color[named]{Black}{%
\special{pn 8}%
\special{ar 1000 200 50 50  0.0000000 6.2831853}%
}}%
%
{\color[named]{Black}{%
\special{pn 8}%
\special{ar 200 200 50 50  0.0000000 6.2831853}%
}}%
%
{\color[named]{Black}{%
\special{pn 8}%
\special{pa 650 200}%
\special{pa 950 200}%
\special{fp}%
}}%
%
{\color[named]{Black}{%
\special{pn 8}%
\special{pa 1050 200}%
\special{pa 1350 200}%
\special{fp}%
}}%
%
{\color[named]{Black}{%
\special{pn 8}%
\special{pa 650 200}%
\special{pa 950 200}%
\special{fp}%
}}%

\put(9.0000,-0.8000){\makebox(0,0)[lb]{$ \pm q^4$}}%

{\color[named]{Black}{%
\special{pn 8}%
\special{pa 1600 200}%
\special{pa 1400 200}%
\special{dt 0.045}%
}}%

\put(5.0000,-0.8000){\makebox(0,0)[lb]{$\pm q^2$}}%

\put(0.0000,-0.8000){\makebox(0,0)[lb]{$\pm 1$}}%
%
{\color[named]{Black}{%
\special{pn 8}%
\special{pa 250 180}%
\special{pa 550 180}%
\special{fp}%
}}%

{\color[named]{Black}{%
\special{pn 8}%
\special{pa 250 220}%
\special{pa 550 220}%
\special{fp}%
}}%

{\color[named]{Black}{%
\special{pn 8}%
\special{pa 350 150}%
\special{pa 450 200}%
\special{fp}%
}}%

{\color[named]{Black}{%
\special{pn 8}%
\special{pa 450 200}%
\special{pa 350 250}%
\special{fp}%
}}%
\end{picture}%

\end{center}
\end{enumerate}

\vs{3ex}
\item  When $q^2$ is a primitive $\ell$-th root of unity, there are
three types of connected components.

\vs{1ex}
\begin{enumerate}
\item $\set{g(\zeta q^{2k})}{k\in\Z/\ell\Z}$
for some $\zeta\not\in\pm q^\Z$. The Dynkin diagram is of type
$A^{(1)}_{\ell-1}$.

\medskip
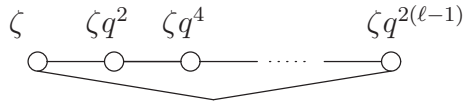
\begin{center}
\unitlength 0.1in
\begin{picture}( 20.5000,  4.5000)(  0.5000, -4.0000)
%
{\color[named]{Black}{%
\special{pn 8}%
\special{ar 600 200 50 50  0.0000000 6.2831853}%
}}%
%
{\color[named]{Black}{%
\special{pn 8}%
\special{ar 1000 200 50 50  0.0000000 6.2831853}%
}}%
%
{\color[named]{Black}{%
\special{pn 8}%
\special{ar 200 200 50 50  0.0000000 6.2831853}%
}}%
%
{\color[named]{Black}{%
\special{pn 8}%
\special{pa 250 200}%
\special{pa 550 200}%
\special{fp}%
}}%
%
{\color[named]{Black}{%
\special{pn 8}%
\special{pa 650 200}%
\special{pa 950 200}%
\special{fp}%
}}%
%
{\color[named]{Black}{%
\special{pn 8}%
\special{pa 1050 200}%
\special{pa 1350 200}%
\special{fp}%
}}%
%
{\color[named]{Black}{%
\special{pn 8}%
\special{pa 650 200}%
\special{pa 950 200}%
\special{fp}%
}}%
\put(8.5000,-0.8000){\makebox(0,0)[lb]{$\zeta q^4$}}%
%
{\color[named]{Black}{%
\special{pn 8}%
\special{pa 1600 200}%
\special{pa 1400 200}%
\special{dt 0.045}%
}}%
\put(4.5000,-0.8000){\makebox(0,0)[lb]{$\zeta q^2$}}%
\put(0.5000,-0.8000){\makebox(0,0)[lb]{$\zeta$}}%
%
{\color[named]{Black}{%
\special{pn 8}%
\special{pa 1700 200}%
\special{pa 2000 200}%
\special{fp}%
}}%
%
{\color[named]{Black}{%
\special{pn 8}%
\special{ar 2050 200 50 50  0.0000000 6.2831853}%
}}%
\put(18.5000,-0.8000){\makebox(0,0)[lb]{\kern1ex$\zeta q^{2(\ell-1)}$}}%
%
{\color[named]{Black}{%
\special{pn 8}%
\special{pa 200 250}%
\special{pa 1120 400}%
\special{fp}%
\special{pa 1120 400}%
\special{pa 2050 250}%
\special{fp}%
}}%
\end{picture}%

\end{center}
\label{case:Loop}

\medskip
\item $\{g(\pm q),\ldots,g(\pm q^{2s+1})\}$
when $\ell$ is odd ($\ell=2s+1$ with $s\ge1$).
In this case $q^{2s+1}=\pm1$.
The Dynkin diagram is of type $A^{(2)}_{2s}$.

\vs{3ex}
\begin{center}
\hs{10ex}
\unitlength 0.1in
\begin{picture}(  5.0000,  3.0000)(  1.5000, -2.5000)
{\color[named]{Black}{%
\special{pn 8}%
\special{ar 600 200 50 50  0.0000000 6.2831853}%
}}%
{\color[named]{Black}{%
\special{pn 8}%
\special{ar 200 200 50 50  0.0000000 6.2831853}%
}}%
{\color[named]{Black}{%
\special{pn 8}%
\special{pa 250 170}%
\special{pa 550 170}%
\special{fp}%
}}%
\put(5.5000,-0.8000){\makebox(0,0)[lb]{\kern-.8ex$\pm q^3$}}%
\put(1.5000,-0.8000){\makebox(0,0)[lb]{\kern-1ex$\pm q$}}%
\put(2.0500,-2.55000){\makebox(0,0)[lb]{\kern-1ex$\times$}}%
{\color[named]{Black}{%
\special{pn 8}%
\special{pa 250 230}%
\special{pa 550 230}%
\special{fp}%
}}%
{\color[named]{Black}{%
\special{pn 8}%
\special{pa 350 200}%
\special{pa 450 150}%
\special{fp}%
}}%
{\color[named]{Black}{%
\special{pn 8}%
\special{pa 450 250}%
\special{pa 350 200}%
\special{fp}%
}}%
{\color[named]{Black}{%
\special{pn 8}%
\special{pa 250 190}%
\special{pa 550 190}%
\special{fp}%
}}%
{\color[named]{Black}{%
\special{pn 8}%
\special{pa 250 210}%
\special{pa 550 210}%
\special{fp}%
}}%
\end{picture}%
\hs{15ex}\text{($(q^2)^3=1$)}
\end{center}

\vs{3ex}
\begin{center}
\unitlength 0.1in
\begin{picture}( 19.6000,  3.0000)(  1.5000, -2.5000)
%
{\color[named]{Black}{%
\special{pn 8}%
\special{ar 600 200 50 50  0.0000000 6.2831853}%
}}%
%
{\color[named]{Black}{%
\special{pn 8}%
\special{ar 200 200 50 50  0.0000000 6.2831853}%
}}%
%
{\color[named]{Black}{%
\special{pn 8}%
\special{pa 250 180}%
\special{pa 550 180}%
\special{fp}%
}}%
%
{\color[named]{Black}{%
\special{pn 8}%
\special{pa 650 200}%
\special{pa 950 200}%
\special{fp}%
}}%
%
{\color[named]{Black}{%
\special{pn 8}%
\special{pa 650 200}%
\special{pa 950 200}%
\special{fp}%
}}%
%
{\color[named]{Black}{%
\special{pn 8}%
\special{pa 1200 200}%
\special{pa 1000 200}%
\special{dt 0.045}%
}}%
\put(5.5000,-0.8000){\makebox(0,0)[lb]{$\pm q^3$}}%
\put(1.5000,-0.8000){\makebox(0,0)[lb]{$\pm q$}}%
\put(2.0500,-2.55000){\makebox(0,0)[lb]{\kern-1ex$\times$}}%
%
{\color[named]{Black}{%
\special{pn 8}%
\special{pa 1300 200}%
\special{pa 1600 200}%
\special{fp}%
}}%
%
{\color[named]{Black}{%
\special{pn 8}%
\special{ar 1650 200 50 50  0.0000000 6.2831853}%
}}%
\put(15.0000,-0.8000){\makebox(0,0)[lb]{$\pm q^{2s-1}$}}%
%
{\color[named]{Black}{%
\special{pn 8}%
\special{pa 250 220}%
\special{pa 550 220}%
\special{fp}%
}}%
%
{\color[named]{Black}{%
\special{pn 8}%
\special{pa 350 200}%
\special{pa 450 150}%
\special{fp}%
}}%
%
{\color[named]{Black}{%
\special{pn 8}%
\special{pa 450 250}%
\special{pa 350 200}%
\special{fp}%
}}%
%
{\color[named]{Black}{%
\special{pn 8}%
\special{ar 2060 200 50 50  0.0000000 6.2831853}%
}}%
%
{\color[named]{Black}{%
\special{pn 8}%
\special{pa 1710 180}%
\special{pa 2010 180}%
\special{fp}%
}}%
%
{\color[named]{Black}{%
\special{pn 8}%
\special{pa 1710 220}%
\special{pa 2010 220}%
\special{fp}%
}}%
%
{\color[named]{Black}{%
\special{pn 8}%
\special{pa 1810 200}%
\special{pa 1910 150}%
\special{fp}%
}}%
%
{\color[named]{Black}{%
\special{pn 8}%
\special{pa 1910 250}%
\special{pa 1810 200}%
\special{fp}%
}}%
\put(20.0000,-0.8000){\makebox(0,0)[lb]{$\pm q^{2s+1}$}}%
\end{picture}%
\hs{5ex}\text{($s>1$)}
\end{center}

\vs{3ex}
\item
$\{g(\pm 1), g(\pm q^2),\ldots, g(\pm q^{2s)})\}$ when $\ell$ is
even ($\ell=2s$ with $s\ge2$). In this case $q^{2s}=-1$. The Dynkin
diagram is of type $C^{(1)}_s$.

\vs{3ex}
\begin{center}
\unitlength 0.1in
\begin{picture}( 19.6000,  3.0000)(  1.5000, -2.5000)
{\color[named]{Black}{%
\special{pn 8}%
\special{ar 600 200 50 50  0.0000000 6.2831853}%
}}%
{\color[named]{Black}{%
\special{pn 8}%
\special{ar 200 200 50 50  0.0000000 6.2831853}%
}}%
{\color[named]{Black}{%
\special{pn 8}%
\special{pa 1710 180}%
\special{pa 2010 180}%
\special{fp}%
}}%
{\color[named]{Black}{%
\special{pn 8}%
\special{pa 650 200}%
\special{pa 950 200}%
\special{fp}%
}}%
{\color[named]{Black}{%
\special{pn 8}%
\special{pa 650 200}%
\special{pa 950 200}%
\special{fp}%
}}%
{\color[named]{Black}{%
\special{pn 8}%
\special{pa 1200 200}%
\special{pa 1000 200}%
\special{dt 0.045}%
}}%

\put(5.5000,-0.8000){\makebox(0,0)[lb]{$\pm q^2$}}%

\put(1.5000,-0.8000){\makebox(0,0)[lb]{$\pm 1$}}%

{\color[named]{Black}{%
\special{pn 8}%
\special{pa 1300 200}%
\special{pa 1600 200}%
\special{fp}%
}}%

{\color[named]{Black}{%
\special{pn 8}%
\special{ar 1650 200 50 50  0.0000000 6.2831853}%
}}%

\put(15.0000,-0.8000){\makebox(0,0)[lb]{$\kern-2ex\pm q^{2(s-1)}$}}%

{\color[named]{Black}{%
\special{pn 8}%
\special{pa 1710 220}%
\special{pa 2010 220}%
\special{fp}%
}}%

{\color[named]{Black}{%
\special{pn 8}%
\special{pa 1810 200}%
\special{pa 1910 150}%
\special{fp}%
}}%

{\color[named]{Black}{%
\special{pn 8}%
\special{pa 1910 250}%
\special{pa 1810 200}%
\special{fp}%
}}%

{\color[named]{Black}{%
\special{pn 8}%
\special{ar 2060 200 50 50  0.0000000 6.2831853}%
}}%

{\color[named]{Black}{%
\special{pn 8}%
\special{pa 250 180}%
\special{pa 550 180}%
\special{fp}%
}}%

{\color[named]{Black}{%
\special{pn 8}%
\special{pa 250 220}%
\special{pa 550 220}%
\special{fp}%
}}%

{\color[named]{Black}{%
\special{pn 8}%
\special{pa 350 250}%
\special{pa 450 200}%
\special{fp}%
}}%

{\color[named]{Black}{%
\special{pn 8}%
\special{pa 450 200}%
\special{pa 350 150}%
\special{fp}%
}}%
\put(20.0000,-0.8000){\makebox(0,0)[lb]{\kern-0.1ex$\pm q^{2s}=\mp1$}}%
\end{picture}%

\end{center}

\vs{3ex}

\item
$\{g(\pm q),g(\pm q^3),\ldots,g(\pm q^{2s-3}), g(\pm q^{2s-1})\}$
where $\ell$ is even ($\ell=2s$ with $s\ge2$).
In this case, $q^{2s}=-1$.
The Dynkin diagram is of type $D^{(2)}_s$.

\vs{3ex}
\begin{center}
\hs{30ex}

\unitlength 0.1in
\begin{picture}(  5.0000,  3.0000)(  1.5000, -2.5000)

{\color[named]{Black}{%
\special{pn 8}%
\special{ar 600 200 50 50  0.0000000 6.2831853}%
}}%

{\color[named]{Black}{%
\special{pn 8}%
\special{ar 200 200 50 50  0.0000000 6.2831853}%
}}%

{\color[named]{Black}{%
\special{pn 8}%
\special{pa 250 180}%
\special{pa 550 180}%
\special{fp}%
}}%

\put(5.5000,-0.8000){\makebox(0,0)[lb]{$\pm q^3=(\mp q)^{-1}$}}%
\put(5.37000,-2.55000){\makebox(0,0)[lb]{$\times$}}%
\put(1.5000,-0.8000){\makebox(0,0)[lb]{\kern-.5ex$\pm q$}}%
\put(1.34500,-2.55000){\makebox(0,0)[lb]{$\times$}}%
{\color[named]{Black}{%
\special{pn 8}%
\special{pa 250 220}%
\special{pa 550 220}%
\special{fp}%
}}%

{\color[named]{Black}{%
\special{pn 8}%
\special{pa 300 200}%
\special{pa 400 150}%
\special{fp}%
}}%

{\color[named]{Black}{%
\special{pn 8}%
\special{pa 400 250}%
\special{pa 300 200}%
\special{fp}%
}}%

{\color[named]{Black}{%
\special{pn 8}%
\special{pa 400 250}%
\special{pa 500 200}%
\special{fp}%
}}%

{\color[named]{Black}{%
\special{pn 8}%
\special{pa 500 200}%
\special{pa 400 150}%
\special{fp}%
}}%
\end{picture}%
\hs{23ex}($s=2$, $(q^2)^2=-1$)
\end{center}

\vs{3ex}
\begin{center}

\unitlength 0.1in
\begin{picture}( 19.6000,  3.0000)(  1.5000, -2.5000)
%
{\color[named]{Black}{%
\special{pn 8}%
\special{ar 600 200 50 50  0.0000000 6.2831853}%
}}%
%
{\color[named]{Black}{%
\special{pn 8}%
\special{ar 200 200 50 50  0.0000000 6.2831853}%
}}%
%
{\color[named]{Black}{%
\special{pn 8}%
\special{pa 250 180}%
\special{pa 550 180}%
\special{fp}%
}}%
%
{\color[named]{Black}{%
\special{pn 8}%
\special{pa 650 200}%
\special{pa 950 200}%
\special{fp}%
}}%
%
{\color[named]{Black}{%
\special{pn 8}%
\special{pa 650 200}%
\special{pa 950 200}%
\special{fp}%
}}%
%
{\color[named]{Black}{%
\special{pn 8}%
\special{pa 1200 200}%
\special{pa 1000 200}%
\special{dt 0.045}%
}}%
\put(5.5000,-0.8000){\makebox(0,0)[lb]{$\pm q^3$}}%
\put(1.5000,-0.8000){\makebox(0,0)[lb]{$\pm q$}}%
\put(2.0500,-2.55000){\makebox(0,0)[lb]{\kern-1ex$\times$}}%
%
{\color[named]{Black}{%
\special{pn 8}%
\special{pa 1300 200}%
\special{pa 1600 200}%
\special{fp}%
}}%
%
{\color[named]{Black}{%
\special{pn 8}%
\special{ar 1650 200 50 50  0.0000000 6.2831853}%
}}%
\put(15.0000,-0.8000){\makebox(0,0)[lb]{$\kern-2ex\pm q^{2s-3}$}}%
%
{\color[named]{Black}{%
\special{pn 8}%
\special{pa 250 220}%
\special{pa 550 220}%
\special{fp}%
}}%
%
{\color[named]{Black}{%
\special{pn 8}%
\special{pa 350 200}%
\special{pa 450 150}%
\special{fp}%
}}%
%
{\color[named]{Black}{%
\special{pn 8}%
\special{pa 450 250}%
\special{pa 350 200}%
\special{fp}%
}}%
%
{\color[named]{Black}{%
\special{pn 8}%
\special{ar 2060 200 50 50  0.0000000 6.2831853}%
}}%
%
{\color[named]{Black}{%
\special{pn 8}%
\special{pa 1710 180}%
\special{pa 2010 180}%
\special{fp}%
}}%
%
{\color[named]{Black}{%
\special{pn 8}%
\special{pa 1710 220}%
\special{pa 2010 220}%
\special{fp}%
}}%
%
{\color[named]{Black}{%
\special{pn 8}%
\special{pa 1810 250}%
\special{pa 1910 200}%
\special{fp}%
}}%
%
{\color[named]{Black}{%
\special{pn 8}%
\special{pa 1910 200}%
\special{pa 1810 150}%
\special{fp}%
}}%
\put(20.0000,-0.8000){\makebox(0,0)[lb]{$\pm q^{2s-1}=(\mp q)^{-1}$}}%
\put(19.980000,-2.55000){\makebox(0,0)[lb]{$\times$}}%
\end{picture}%
\hs{15ex}($s>2$)
\end{center}

\end{enumerate}
\ee

\bigskip
Observe that any vertex of $\AAA^1_w$
has one or two edges, and $\{g(q),g(-q),g(1),g(-1)\}$ is the set of
points with one edge. Note that the odd vertices have the minimal
length.

\bigskip
Let $J$ be a finite subset of $\T_X$ invariant under $c$,
and let $X\cl J\to \T_X$ be the inclusion map.
 Let $I$ be the image of the composition $J\To[X]\T_X\To[\pr]\AAA^1_w$,
 and we also denote by $\pr\cl J\to I$ the surjection.
Set $J^c\seteq\set{j\in J}{c(j)=j}$ and
$I_\odd=\pr(J^c)=I\cap\{g(\pm q)\}$, $I_\even=I\setminus I_\odd$.

 We choose a function $h\cl I\to \T_\la$
 such that
 \eq
 &&\left\{\parbox{55ex}{\be[{\small$\bullet$}]
 \item $g(h(i))=i$,

\vs{1ex}
 \item if $(\alpha_i,\alpha_{i'})<0$, 
 then $h(i)=q^{\pm2}h(i')$ for any $i,i'\in I$,
 \ee}\right.
 \label{cond:h}\eneq
Such a function $h$ exists. Indeed, it is enough to show it for a
connected component of $\AAA_w$. Since any vertices has at most two
edges, the assertion is obvious when the connected component has no
loop,
 and the remaining case \eqref{case:Loop} can
be checked directly.

\medskip
 We define the map
$\lambda \cl J\to \T_\lambda$ by $\lambda(j)=h(\pr(j))$. Hence we
have $X(j)+X(j)^{-1}=2\dfrac{\la(j)+\la(j)^{-1}}{q+q^{-1}}$. Since
$h\cl I \to\T_\la$ is injective, we
have \eq && \ba{l}\text{for $j,j'\in J$, we have $\la(j)=\la(j')$ if
and only if $j=j'$ or $j=c(j')$.} \ea \eneq

\medskip
We choose a function $\eps\cl J\to\{0,1\}$ such that $\eps^{-1}(0)\to I$
is bijective. Hence we have
$$\text{$\eps(cj)=1-\eps(j)$ for $j\in J\setminus J^c$, and $\eps(j)=0$ for $j\in J^c$.}$$

For $i,j\in J$, we write $(i,j)$ for $(\alpha_{\pr(i)},\alpha_{\pr(j)})$
and $a_{i,j}\seteq2(i,j)/(i,i)$.
For $i,j\in J$ such that $i\not=j$ and $\eps(i)=\eps(j)=0$, we choose
\eq &&
\tQ_{i,j}(u,v)=
\begin{cases}
\pm(u^{-a_{i,j}}-v^{-a_{j,i}} )&\text{if $(i,j)<0$,}\\
1&\text{if $(i,j)=0$,}
\end{cases}\label{def:tQ}
\eneq
where we choose $\pm$ such that $\tQ_{i,j}(u,v)=\tQ_{j,i}(v,u)$.

We can extend uniquely the definition of $\tQ_{i,j}$ for all $i,j\in
J$ such that the conditions in \eqref{cond:Q} hold; namely,
\eq&&\tQ_{i,j}(u,v)
=\tQ_{c^{\eps(i)}(i),\,c^{\eps(j)}(j)}((-1)^{\eps(i)}u,\,(-1)^{\eps(j)}v).
\label{eq:extQ}
\eneq

Then we can associate to the data $\{\tQ_{i,j}\}_{i,j\in J}$ the quiver Hecke-Clifford superalgebra $\RC_n$
and its completion $\HRC$.

The following theorem is a main result of this section.

 \Th\label{th:isom_aDHC} We have $\HRC \simeq\OHC$. \enth

The rest of this subsection is devoted to the proof of this theorem.

\bigskip
\subsection{Proof of Theorem~\ref{th:isom_aDHC}}
\subsubsection{Strategy of the proof}
We  will construct the elements $y_k\in\CO_n\subset\OHC$
($k=1,\ldots,n$) and $\ts_a\in\KHC$ ($a=1,\ldots,n-1$)
 such that
\eq&&\left\{
\parbox{70ex}{
\bnum
\item $y_pe(\nu)\in
(\HCO_{\T^n,X(\nu)})^{\times}\,\bl(X_p-X(\nu_p)\br)e(\nu)$,
\item $C_ay_pe(\nu)=(-1)^{\delta_{a,p}}y_pe(c_a\nu)C_a$,

\vs{.5ex}
\item $\ts_a^2=R_{a,a+1}$ (see \eqref{def R}),\label{eq:tsR}
\item $\{\ts_a\}_{1\le a<n}$ satisfies the braid relations (see \eqref{eq:braid}),
\label{eq:tsbraid}

\vs{.5ex}
\item $\ts_aC_k=C_{s_a(k)}\ts_a$, $\ts_ay_k=y_{s_a(k)}\ts_a$,
\label{eq:tsperm}
\vs{0ex}
\item setting $\sigma_a\seteq\ts_a+f_{a,a+1}$ (see \eqref{def:f}),
the element $\sigma_ae(\nu)$ belongs to
$T_ae(\nu)({\HAT{\CO}}_{\T^n,X(\nu)})^{\times}
+\CO_n\langle C_1,\cdots,C_n\rangle$.\label{eq:Taureg}
\ee}\right.
\label{eq:hcrc}
\eneq
Setting $\cc_k=C_k$, the elements
 $y_p\,$, $\cc_p$ ($p=1,\ldots,n$), $\sigma_a$ ($a=1,\ldots,n-1$)
satisfy the defining relations of $\RC_n$ by Theorem~\ref{poly_rep}.
Hence we obtain a homomorphism $F\cl\HRC\to \OHC$. By
Corollary~\ref{cor:pbwrc},
$\{\cc^\eta\sigma_w\}_{\mbox{{\scriptsize$\eta\in(\FF)^n$,
$w\in\Sym_n$}}}$ is a basis of $\HRC$ as a left $\CO_n$-module. On
the other hand, it is easy to see that its image by $F$ forms a
basis of $\OHC$ as a left $\CO_n$-module by
Proposition~\ref{prop:pbw}. Hence we conclude that $F\cl\HRC\to
\OHC$ is an isomorphism.

\smallskip
Now we will construct $y_k\in\CO_n$ ($k=1,\ldots,n$) and
$\ts_a\in\KHC$ ($a=1,\ldots,n-1$) which satisfy the
condition~\eqref{eq:hcrc}.

\subsubsection{Definition of $y_p\,e(\nu)$}\quad

We will construct $y_p\,e(\nu)$ such that
$$
\begin{aligned}
& y_pe(\nu)\in
({\HCO}_{\T^n,X(\nu)})^{\times}\,(X_p-X(\nu_p))e(\nu), \\
& C_a y_pe(\nu)=(-1)^{\delta_{ap}}y_pe(c_a\nu)C_a.
\end{aligned}
$$

Let $C$ be the subscheme of $\T_X\times\T_\la$ defined by the equation
$$X+X^{-1}=2\dfrac{\la+\la^{-1}}{q+q^{-1}}.$$
Then $C$ is a smooth curve and we have a Cartesian product \eqn
&&\xymatrix{
C\ar[r]\ar[d]&\T_\la\ar[d]^g\\
\T_X\ar[r]^\pr&\AAA_w.
}
\eneqn
Then for every $i\in J$, it induces
the injective homomorphisms
\eq
&&\HCO_{\T_X, X(i)}\monoto\HCO_{C,(X(i),\la(i))}\monofrom\HCO_{\T_\la, \la(i)}.
\eneq
Let us set $\CO_i\seteq\HCO_{\T_X,X(i)}=\cor[[X-X(i)]]$
and $\tCO_i\seteq\HCO_{C,\;(X(i),\,\la(i))}$,
and we regard $\CO_i$ as a subalgebra of $\tCO_i$.
The algebra $\CO_i$ is a discrete valuation ring with the maximal ideal $\m_i\seteq\CO_i(X-X(i))$.
The indeterminates $X$ and  $\la$ are considered as elements of $\tCO_i$ and
they satisfy $X+X^{-1}=2\dfrac{\la+\la^{-1}}{q+q^{-1}}$.
Let $c_i\cl\CO_i\to\CO_{c(i)}$ be the homomorphism induced by the map $X\mapsto X^{-1}$.
For $1\le a\le n$, let us denote by $q_{a,\nu}\cl\CO_{\nu_a}\To\HCO_{\T^n,X(\nu)}e(\nu)$
the homomorphism induced by the $a$-th projection $(X_1,\ldots,X_n)\mapsto X_a$.

\medskip
Recall that $\eps\cl J\to\{0,1\}$ is a map such that
$\eps^{-1}(0)\to I$  is bijective. For all $i\in J$ with
$\eps(i)=0$, we will construct $y_i\in\CO_i$ such that \eq
&&\left\{\parbox{60ex}{ \bnum
\item
$y_i\in \CO_i^\times(X-X(i))$,

\vs{.5ex}
\item if $c(i)=i$, then
$c_i(y_i)=-y_i$, \ee }\right. \label{cond:y} \eneq and define
$$y_ae(\nu)=\begin{cases}
q_{a,\nu}(y_{\nu_a})&\text{if $\eps(\nu_a)=0$,}\\
-c_aq_{a,c_a\nu}(y_{c(\nu_a}))&\text{if $\eps(\nu_a)=1$.}\end{cases}$$
Then we have $y_pe(\nu)\in
({\HCO}_{\T^n,X(\nu)})^{\times}\,(X_p-X(\nu_p))e(\nu)$ and $C_a
y_pe(\nu)=(-1)^{\delta_{ap}}y_pe(c_a\nu)C_a$.

\medskip
Now, we will construct $y_i$ for $i\in J$ such that $\eps(i)=0$. We
choose $\psi(\la)\in\HCO_{\T_{\la},1}$ such that
\eq&&\ba{l}
\psi\in(\HCO_{\T,1})^\times(\la-1),\\[1ex]
\psi(\la^{-1})=-\psi(\la). \ea \eneq
For example, we can take
$(\la-1)/(\la+1)$ or $\la-\la^{-1}$ as $\psi(\la)$.

\bigskip
\noindent
$\bullet$ Case $\lambda(i)\ne \pm1,\pm q, \pm q^{-1}$.

The projections $C\to \T_{X}$ and $C\to \T_{\la}$
are \'etale at $(X(i),\la(i))$, and hence
we have isomorphisms
\begin{align*}
\cor[[X-X(i)]]=\CO_i
 \isoto\tCO_i=\HAT{\CO}_{C,\,(X(i),\la(i))}
\isofrom \HAT{\CO}_{\T_{\lambda},\lambda(i)}
=\cor[[\lambda-\lambda(i)]].
\end{align*}
We define $y_i\in\CO_{i}$ by $y_i=\psi(\la(i)^{-1}\la)$.

\bigskip
\noindent
$\bullet$ Case $\lambda(i)=\pm q, \pm q^{-1}$.

\smallskip
In this case, we have $X(i)=\pm1$ and
\begin{align*}
\cor[[X-X(i)]]=\CO_i=\HAT{\CO}_{\T,X(i)}
 \isoto
\tCO_i=\HAT{\CO}_{C,\,(X(i),\la(i))}
\monofrom\HAT{\CO}_{\T_{\lambda},\la(i)}=\cor[[\lambda-\la(i)]] .
\end{align*}
Then $\HAT{\CO}_{\T_{\lambda},\la(i)}$ is identified with
$\set{f(X)\in \CO_i}{c(f)=f}$. (Recall that $c(f)(X)=f(X^{-1})$.)
Since $\CO_i\psi(\la(i)^{-1}\la)=\m_i^2$, there exists
$y_i\in \CO_i$ such that $y_i^2=\psi(\la(i)^{-1}\la)$.
Then $y_i$ generates the maximal ideal $\m_i$.
Since $c(y_i)^2=y_i^2$ and $c(y_i)$ cannot be equal to $y_i$, we obtain $c(y_i)=-y_i$.

\bigskip
\noindent
$\bullet$ Case $\lambda(i)=\pm 1$.

\smallskip
In this case, we have
\begin{align*}
\HCO_{\T_{X},X(i)}
\monoto \HCO_{C,\,(X(i),\pm1)}
\isofrom\HCO_{\T_{\lambda},\,\pm1}.
\end{align*}
Then $\CO_i$ is identified with $\set{f(\la)\in\HCO_{\T_{\lambda},\,\pm1}}{f(\la^{-1})=f(\la)}$.
 Hence $\psi(\la(i)^{-1}\la)^2$ belongs to $\CO_i$.
We define $y_i\in \CO_i$ by
$y_i=\psi(\la(i)^{-1}\la)^2$.
Then we have $\m_i=\CO_i y_i$.

\subsubsection{Definition of $\ts_ae(\nu)$}
\quad We will define $\ts_ae(\nu)$ for $1\le a<n$ and $\nu\in J^n$.

In the preceding subsection, we have constructed $y_i\in
\CO_i$ such that $y_ae(\nu)=q_{a,\nu}(y_{\nu_a})$ for any $\nu\in
J^n$. For $i,j\in J$, let us denote
$$\text{$\CO_{i,j}=\HCO_{\T^2,(X(i),X(j))}$
and $\tCO_{i,j}=\HCO_{C\times C,\;((X(i),\la(i)),(X(j),\la(j)))}$.}$$
We regard $\CO_{i,j}$ as a subalgebra of $\tCO_{i,j}$. Let
$r^1_{i,j}\cl\tCO_i\to \tCO_{i,j}$ and $r^2_{i,j}\cl\tCO_j\to
\tCO_{i,j}$ be the algebra homomorphisms induced by the first and
second projections from $C\times C$ to $C$, respectively. We write
$y_1$ for $r^1_{i,j}(y_i)\in\CO_{i,j}$ and $y_2$ for
$r^2_{i,j}(y_j)\in\CO_{i,j}$. Similarly, we define
$\la_1\in\tCO_{i,j}$ and $\la_2\in\tCO_{i,j}$ as $r^1_{i,j}(\la)$ and
$r^2_{i,j}(\la)$. Note that $X_1,X_2\in\CO_{i,j}$ may be regarded as
$r^1_{i,j}(X)$ and $r^2_{i,j}(X)$.

Let $c_1\cl\CO_{i,j}\to\CO_{c(i),j}$ and $c_2\cl\CO_{i,j}\to\CO_{i,c(j)}$
be the isomorphisms induced by $(X_1,X_2)\mapsto (X_1^{-1},X_2)$ and $(X_1,X_2)\mapsto (X_1,X_2^{-1})$,
respectively.
We denote by $s_{12}\cl \CO_{i,j}\to\CO_{j,i}$ the
isomorphism induced by $\T^2\ni(X_1,X_2)\mapsto(X_2,X_1)\in\T^2$.

Let $\psi_{a,\nu}\cl\CO_{\nu_a,\nu_{a+1}}\To\HCO_{\T^n,X(\nu)}e(\nu)$
be the algebra homomorphism induced by
the projection $(X_1,\ldots,X_n)\longmapsto(X_a,X_{a+1})$.
Let us define $R_{i,j}\in\CO_{i,j}$ 
by
\eq
R_{i,j}&=&\begin{cases}
\tQ_{i,j}(y_1,y_2)&\text{if $j\not=i,c(i)$,}\\
-(y_1-y_2)^{-2}&\text{if $j=i\not\in J^c$,}\\
-(y_1+y_2)^{-2}&\text{if $j=c(i)\not\in J^c$,}\\
-(y_1-y_2)^{-2}-(y_1+y_2)^{-2}=-2\dfrac{y_1^2+y_2^2}{(y_1^2-y_2^2)^{2}}&\text{if $j=i\in J^c$.}
\end{cases}\label{eq:R}
\eneq
Then we have $R_{a,a+1}e(\nu)=\psi_{a,\nu}(R_{\nu_a,\nu_{a+1}})$.

Let us recall (see Remark~\ref{Intertwin})
\begin{align*}
\Phi_a^2
&=F(X_a,X_{a+1}).
\end{align*}

\Lemma\label{lem:FR} For $i,j\in J$ such that $\eps(i)=\eps(j)=0$,
the following statements hold.
\bnum
\item $F(X_1,X_2)^{-1}R_{i,j}$ belongs to $\CO_{i,j}{}^\times$.
\item
$F(X_1,X_2)^{-1}R_{i,i}-\left(\dfrac{X_1X_2^{-1}-1}{\xi(y_1-y_2)}\right)^2$
belongs to $\CO_{i,i}(X_1-X_2)$. \ee \enlemma
We will prove this
lemma later: (i) by case-by-case verification and (ii) as a
consequence of (i). Admitting this lemma for a while, we shall
construct $\ts_ae(\nu)$ and prove the theorem.

\Lemma For $i,j\in J$ such that $\eps(i)=\eps(j)=0$, there exists a
family of elements  $G_{i,j}\in\CO_{i,j}$ satisfying the following
conditions:
\eq&&\left\{
\parbox{70ex}{
\bna
\item
$G_{i,j}\cdot s_{12}G_{j,i} = F(X_1,X_2)^{-1}R_{i,j}$,\label{eq:Gsg}
\vs{1ex}
\item $G_{i,i}-\dfrac{X_1X_2^{-1}-1}{\xi(y_1-y_2)}\in\CO_{i,i}(X_1-X_2)$,
\label{cond:Gtau}
\vs{.5ex}
\item if $c(i)=i$, then $c_1\bl(G_{i,j}\br)=G_{i,j}$,
 and if $c(j)=j$, then $c_2\bl(G_{i,j}\br)=G_{i,j}$.\label{cond:cG}
\ee
}\right.\label{cond:G}
\eneq
\enlemma \Proof Assume first that $i\not=j$. In order to see
\eqref{eq:Gsg}, it is enough to take $G_{j,i}=1$ and
$G_{i,j}=F(X_1,X_2)^{-1}R_{i,j}$. Indeed, we have $G_{j,i}\cdot
s_{12}G_{i,j}=F(X_1,X_2)^{-1}R_{j,i}$ because $F(X_1,X_2)$ is
invariant under $s_{12}$ and $s_{12}R_{i,j}=R_{j,i}$.

If $i=j$, then by  Lemma~\ref{lem:FR}, there exists a unique
$G_{i,i}$ satisfying $G_{i,i}^2=F(X_1,X_2)^{-1}R_{i,i}$ and the
condition \eqref{cond:Gtau}. Since
$F(X_1,X_2)^{-1}R_{i,i}$ is invariant under $s_{12}$, we derive that
$(s_{12}G_{i,i})^2=G_{i,i}^2$. Since $G_{i,i}$ and $s_{12}G_{i,i}$ take
the same non-zero value at $(X(i),X(j))$, we have
$s_{12}G_{i,i}=G_{i,i}$. Hence \eqref{eq:Gsg} is satisfied.

Let us show \eqref{cond:cG}. Assume $c(i)=i$. Since
$F(X_1,X_2)^{-1}R_{i,j}$ is invariant under $c_1$, we have
$c_1(G_{i,j})=\pm G_{i,j}$. Since $G_{i,j}$ and $c_1G_{i,j}$ take
the same non-zero value at $(X(i),X(j))$, $c_1G_{i,j}=G_{i,j}$
holds. Similarly, we have $c_2G_{i,j}=G_{i,j}$ if
$c(j)=j$. \QED
\medskip
Thus we have proved that there exists $\bl(G_{i,j}\br)_{i,j}$
($i,j\in J$ with $\eps(i)=\eps(j)=0$) which satisfies the
conditions~\eqref{cond:G}. We extend the definition of $G_{i,j}$ for
all $i,j\in J$ by
$$G_{i,j}=c_1^{\eps(i)}c_2^{\eps(j)}G_{c^{\eps(i)}(i),\,c^{\eps(j)}(j)}.$$
Then $(G_{i,j})_{i,j\in J}$ satisfies
\eqref{cond:G}\,\eqref{eq:Gsg}, \eqref{cond:G}\,\eqref{cond:Gtau}
and \eq&&\text{$c_1(G_{i,j})=G_{c(i),j}$ and
$c_2(G_{i,j})=G_{i,c(j)}$}\label{eq:cG} \eneq as well.

Now we define $\ts_a\in\KHC$ by
\eq&&\ts_ae(\nu)=\Phi _ae(\nu)\psi_{a,\nu}(G_{\nu_{a},\nu_{a+1}}).\eneq

Let us verify the conditions in \eqref{eq:hcrc}. By
Lemma~\ref{lem:braid}, the $\ts_a$'s satisfy the braid relations. By
the construction, we have
$$\text{$\ts_a\cdot e(\nu)f=s_a(f)e(s_a\nu)\cdot\ts_a$ for any $f\in \FRAC({\HAT{\CO}}_{{\T}^n,X(\nu)})$.}$$
The condition \eqref{eq:cG} implies that
$$\ts_aC_b=C_{s_a(b)}\ts_a.$$
Hence we have proved the conditions \eqref{eq:tsbraid} and \eqref{eq:tsperm} in \eqref{cond:G}.

The hypothesis \eqref{cond:G}\;\eqref{eq:Gsg} implies
$$\text{$\ts_a^2e(\nu)=R_{a,a+1}e(\nu)$.}$$
Indeed,
$s_a\bl(\psi_{a,s_a\nu}(G_{\nu_{a+1},\nu_{a}})\br)=\psi_{a,\nu}(s_{12}G_{\nu_{a+1},\nu_{a}})$
and we have \eqn \ts_a^2e(\nu)&=&\ts_ae(s_a\nu)\cdot\ts_ae(\nu)=
\Phi_ae(s_a\nu)\psi_{a,s_a\nu}(G_{\nu_{a+1},\nu_{a}})\Phi_ae(\nu)\psi_{a,\nu}(G_{\nu_{a},\nu_{a+1}})
\\
&=&\Phi_a^2e(\nu)s_a\bl(\psi_{a,s_a\nu}(G_{\nu_{a+1},\nu_{a}})\br)\psi_{a,\nu}(G_{\nu_{a},\nu_{a+1}})\\
&=&\Phi_a^2e(\nu)\psi_{a,s_a\nu}\bl(s_{12}G_{\nu_{a+1},\nu_{a}}\cdot G_{\nu_{a},\nu_{a+1}}\br)\\
&=&F(X_a,X_{a+1})e(\nu)\bl(F(X_a,X_{a+1})^{-1}R_{a,a+1}e(\nu)\br)=R_{a,a+1}e(\nu).
\eneqn

\medskip

Using the condition \eqref{cond:G}\;\eqref{cond:Gtau}, let us show
 the condition \eqref{eq:hcrc}\;\eqref{eq:Taureg}:
$$\sigma_ae(\nu)\seteq\ts_ae(\nu)+f_{a,a+1}e(\nu)\in
T_ae(\nu)(\HCO_{\T^n,X(\nu)})^{\times}+\CO_n\langle C_1,\cdots,C_n\rangle.$$
Setting $i=\nu_a$ and $j=\nu_{a+1}$ and
 writing $G_{i,j}$ for $\psi_{a,\nu}(G_{\nu_{a},\nu_{a+1}})$, we have
 \eqn
\sigma_ae(\nu)&=&\Bigl(\ts_a-(y_a-y_{a+1})\delta_{i,j}+(y_a+y_{a+1})^{-1}C_aC_{a+1}\delta_{c(i),j}\Bigr)e(\nu)\\
&=&\Bigl(T_a+\xi(X_aX_{a+1}^{-1}-1)^{-1}-C_aC_{a+1}\xi(X_a^{-1}X_{a+1}^{-1}-1)^{-1}\Bigr)e(\nu)G_{i,j}\\
&&\hs{10ex}-(y_a-y_{a+1})^{-1}\delta_{i,j}e(\nu)-C_aC_{a+1}(y_a+y_{a+1})^{-1}\delta_{c(i),j}e(\nu)\\
&=&T_ae(\nu)G_{i,j}+\Bigl(\xi(X_aX_{a+1}^{-1}-1)^{-1}G_{i,j}-(y_a-y_{a+1})^{-1}\delta_{i,j}\Bigr)e(\nu)\\
&&\hs{5ex}+C_a\Bigl(-\xi(X_a^{-1}X_{a+1}-1)^{-1}(c_{a+1}G_{i,j})+(-y_a+y_{a+1})^{-1}\delta_{c(i),j}\Bigr)C_{a+1}e(\nu)\\
&=&T_ae(\nu)G_{i,j}+\Bigl(\xi(X_aX_{a+1}^{-1}-1)^{-1}G_{i,j}-(y_a-y_{a+1})^{-1}\delta_{i,j}\Bigr)e(\nu)\\
&&\hs{5ex}+C_{a}\Bigl(\xi\bl((X_aX_{a+1}^{-1}-1)^{-1}+1\br)G_{i,c(j)}
-(y_a-y_{a+1})^{-1}\delta_{i, c(j)}\Bigr)e(c_{a+1}\nu)C_{a+1}.
\eneqn
Here, the last equality follows from \eqref{eq:cG}. Hence it is
enough to show that
\eq&&\xi(X_aX_{a+1}^{-1}-1)^{-1}G_{i,j}-(y_a-y_{a+1})^{-1}\delta_{i,j}\in\CO_n.
\label{eq:reg} \eneq
 If $i\not=j$, then $(X_aX_{a+1}^{-1}-1)^{-1}e(\nu)\in\CO_n$ and hence
\eqref{eq:reg} holds.

If $i=j$, then
$$\xi(X_aX_{a+1}^{-1}-1)^{-1}G_{i,j}-(y_a-y_{a+1})^{-1}
=\xi(X_aX_{a+1}^{-1}-1)^{-1}\Bigl(G_{i,j}-(X_aX_{a+1}^{-1}-1)\xi^{-1}(y_a-y_{a+1})^{-1}\Bigr)$$
belongs to $\CO_n$ by \eqref{cond:G}\;\eqref{cond:Gtau}. Thus we
have verified all the conditions in \eqref{eq:hcrc}, and hence
$\widehat{\RC}_n $ and $\OHC$ are isomorphic.

\medskip
 \subsubsection{Proof of Lemma~\ref{lem:FR}}
 Now,
 we shall prove Lemma~\ref{lem:FR} for $i,j\in J$ such that $\eps(i)=\eps(j)=0$.

 \smallskip
 Let us first derive (ii) admitting (i).
 Let $\eta$ be a generic point of the irreducible subscheme $\{X_1=X_2\}$
 of $\on{Spec}(\CO_{i,i})$, and let $(\CO_{i,i})_\eta$ be
 the localization of $\CO_{i,i}$ at $\eta$.
 Then (i) implies that
$S\seteq F(X_1,X_2)^{-1}R_{i,i}-\Bigl(\dfrac{X_1X_2^{-1}-1}{\xi(y_1-y_2)}\Bigr)^{2}$
belongs to $\CO_{i,i}$.
Hence in order to see (ii), it is enough to show that
$S$ belongs to $(\CO_{i,i})_\eta(X_1-X_2)$.

We have, modulo $(\CO_{i,i})_\eta(X_1-X_2)$
\eqn
&&F(X_1,X_2)^{-1}(X_1X_2^{-1}-1)^{-2}\\
&&=\dfrac{(X_2-X_1^{-1})^2}%
{(X_1+X_1^{-1})^2-(q^2+q^{-2})(X_1+X_1^{-1})(X_2+X_2^{-1})+(X_2+X_2^{-1})^2+4(q-q^{-1})^2}\\
&&\equiv\dfrac{(X_1-X_1^{-1})^2}%
{(2-q^2-q^{-2})(X_1+X_1^{-1})^2+4(q-q^{-1})^2}\\
&&=\dfrac{(X_1-X_1^{-1})^2}%
{-(q-q^{-1})^2(X_1-X_1^{-1})^2}=-\xi^{-2}
\eneqn
Since $R_{i,j}\equiv -(y_1-y_2)^{-2} \mod (\CO_{i,j})_\eta$,
we obtain
$$F(X_1,X_2)^{-1}R_{i,j}\equiv \dfrac{(X_1X_2^{-1}-1)^2}{\xi^2(y_1-y_2)^2}
\mod (\CO_{i,j})_\eta(X_1-X_2).$$

\bigskip
It only remains to prove (i). Let us show it by case-by-case
verification. Let us recall that \eq&&\ba{rl}
F(X_1,X_2)^{-1}&=\dfrac{(\la_2-\la_{1})^2(\la_2-\la_1^{-1})^2}%
{(\la_2-q^2\la_1)(\la_2-q^{-2}\la_1)(\la_2-q^2\la_1^{-1})(\la_2-q^{-2}\la_1^{-1})}.
\ea\label{eq:RF} \eneq
Note that (see \eqref{cond:h}) for $i,j\in J$
such that $\eps(i)=\eps(j)=0$, we have
$$
\parbox{60ex}
{
\bna
\item $\la(j)=\la(i)^{-1}$ implies $i=j$,

\vs{0.5ex}
\item $(\alpha_{\pr(i)},\alpha_{\pr(j)})<0$
implies $\la(j)=q^{\pm2}\la(i)$,
\ee}
$$

\bigskip
\noindent$\bullet$ Case $\lambda(j)\ne \lambda(i),
q^{\pm2}\lambda(i)$.

\smallskip
In this case, $F(X_1,X_2)\in\CO^{\times}_{i,j}$ and $R_{i,j}=\tQ_{i,j}=1$.

\bigskip
\noindent
$\bullet$ Case $\lambda(j)=q^{\pm2}\lambda(i)$ and
$\lambda(i),\lambda(j)\ne \pm1,\pm q,\pm q^{-1}$.

\smallskip
We may assume  that $\lambda(j)=q^2\lambda(i)$. In this case,
$y_1=\psi(\la(i)^{-1}\la_1)$, $y_2=\psi(\la(j)^{-1}\la_2)$, and
$$R_{i,j}=\tQ_{i,j}(y_1,y_2)=\pm(y_1-y_2)\in
\CO_{i,j}{}^\times(\la(i)^{-1}\la_1-\la(j)^{-1}\la_2)
=\CO_{i,j}{}^\times(\la_2-q^{2}\la_1).$$ Therefore, since
\begin{align*}
F(X_1,X_2)^{-1}R_{i,j}
&\in\CO_{i,j}{}^\times
\dfrac{(\la_2-\la_{1})^2(\la_2-\la_1^{-1})^2}%
{(\la_2-q^{-2}\la_1)(\la_2-q^2\la_1^{-1})(\la_2-q^{-2}\la_1^{-1})},
\end{align*}
$F(X_1,X_2)^{-1}R_{i,j}$ belongs to $\CO_{i,j}^{\times}$.

\bigskip
\noindent$\bullet$ Case $(\lambda(i),\,\lambda(j))=\pm (1, q^{2c})$ or $\pm( q^{2c},1)$
for some $c=1,-1$

\smallskip
We may assume that $(\lambda(i),\,\lambda(j))=\pm (1, q^{2c})$. Thus
we have $y_1=\psi(\la(i)^{-1}\la_1)^2$. Moreover, we have
$\pr(j)=g(\pm q)$ if and only if $(q^2)^3=1$.

\medskip
Let us first assume $(q^2)^3=1$. In this case,
$y_{2}^2=\psi(\la(j)^{-1}\la_2)$,
and
\eqn
 R_{i,j}&=& \tQ_{i,j}(y_1,y_2)=
 \pm(y_1-y_2^4)
 =\pm\bl(\psi(\la(i)^{-1}\la_1)^2
 -\psi(\la(j)^{-1}\la_2)^2\br)\\
&=&\pm\bl(\psi(\la(i)^{-1}\la_1)-\psi(\la(j)^{-1}\la_2)\br)
\bl(\psi(\la(i)^{-1}\la_1)+\psi(\la(j)^{-1}\la_2)\br)\\
&=&\pm\bl(\psi(\la(i)^{-1}\la_1)-\psi(\la(j)^{-1}\la_2)\br)
\bl(-\psi(\la(i)\la_1^{-1})+\psi(\la(j)^{-1}(\la_2))\br)\\
&\in&\CO_{i,j}{}^\times(\la(j)^{-1}\la_2-\la(i)^{-1}\la_1)(\la(j)^{-1}\la_2-\la(i)\la_1^{-1})\\
&=&\CO_{i,j}{}^\times(\la_2-\la(j)\la(i)^{-1}\la_1)(\la_2-\la(i)\la(j)\la_1^{-1})\\
&=&\CO_{i,j}{}^\times(\la_2-q^{2c}\la_1)(\la_2-q^{2c}\la_1^{-1}).
\eneqn
Hence we have
\eqn
F(X_1,X_2)^{-1}R_{i,j}&\in&\CO_{i,j}{}^\times
\dfrac{(\la_2-\la_{1})^2(\la_2-\la_1^{-1})^2}%
{(\la_2-q^{-2c}\la_1)(\la_2-q^{-2c}\la_1^{-1})}=\CO_{i,j}^\times.
\eneqn

\medskip
Now assume that $(q^2)^3\not=1$.
 In this case,
$y_1=\psi(\la(i)^{-1}\la_1)^2$, $y_{2}=\psi(\la(j)^{-1}\la_2)$,
$R_{i,j}=\tQ(y_1,y_2)=y_2^2-y_1=\psi(\la(i)^{-1}\la_1)^2-\psi(\la(j)^{-1}\la_2)^2$.
Hence, similarly to the case of $(q^2)^3=1$, we have
$R_{i,j}\in\CO_{i,j}{}^\times(\la_2-q^{2c}\la_1)(\la_2-q^{2c}\la_1^{-1})$
and hence we have
\eqn
F(X_1,X_2)^{-1}R_{i,j}&\in&\CO_{i,j}{}^\times
\dfrac{(\la_2-\la_{1})^2(\la_2-\la_1^{-1})^2}%
{(\la_2-q^{-2c}\la_1)(\la_2-q^{-2c}\la_1^{-1})}=\CO_{i,j}^\times.
\eneqn

\bigskip
\noindent$\bullet$ Case $(\lambda(i),\lambda(j))=\pm(q^{c},q^{3c})$ or $\pm(q^{3c},q^{c})$
for some $c=1,-1$.

\smallskip
Assume that $(\lambda(i),\lambda(j))=\pm(q^c,q^{3c})$.
 First note that
$\la(j)=\pm1$ if and only if $(q^2)^3=1$,
and $\pr(j)=g(\pm q)$ if and only if $(q^2)^4=1$.

\noi The case $(q^2)^3=1$ was already treated. \noi Assume that
$(q^2)^4=1$. Then $q^4=-1$ and hence $\la(j)=\mp q^{-c}$. In this
case, $y_1^2=\psi(\la(i)^{-1}\la_1)$ and
$y_2^2=\psi(\la(j)^{-1}\la_2)$, and $R_{i,j}=\pm(y_1^2-y_2^2)$. Hence
$R_{i,j}\in\CO_{i,j}{}^\times(\la(i)^{-1}\la_1-\la(j)^{-1}\la_2)
=\CO_{i,j}{}^\times(\la_2-q^{2c}\la_1)$. Therefore, we have \eqn
F(X_1,X_2)^{-1}R_{i,j}&\in&\CO_{i,j}{}^\times
\dfrac{(\la_2-\la_{1})^2(\la_2-\la_1^{-1})^2}%
{(\la_2-q^{-2c}\la_1)(\la_2-q^2\la_1^{-1})(\la_2-q^{-2}\la_1^{-1})}=\CO_{i,j}^\times.
\eneqn

Assume now $(q^2)^k\not=1$ for $k=3,4$. In this case,
$y_1^2=\psi(\la(i)^{-1}\la_1)$, $y_{2}=\psi(\la(j)^{-1}\la_2)$, and
$R_{i,j}=\tQ_{i,j}(y_1,y_2)=\pm(y_1^2-y_2)\in\CO_{i,j}{}^\times(\la_2-q^{2c}\la_1)$.
Therefore, since
\eqn F(X_1,X_2)^{-1}R_{i,j}&\in &\CO_{i,j}{}^\times
\dfrac{(\la_2-\la_{1})^2(\la_2-\la_1^{-1})^2}%
{(\la_2-q^{-2c}\la_1)(\la_2-q^2\la_1^{-1})(\la_2-q^{-2}\la_1^{-1})},
\eneqn $F(X_1,X_2)^{-1}R_{i,j}$ belongs to $\CO_{i,j}^\times$.

\bigskip
\noindent
$\bullet$ Case $i=j$ and
$\lambda(i)\ne \pm1,\pm q^{\pm1}$.

\smallskip
In this case, $y_k=\psi(\la(i)^{-1}\la_k)$ ($k=1,2$)
and $R_{i,j}=-(y_1-y_2)^{-2}\in \CO_{i,j}{}^\times (\la_1-\la_2)^{-2}$.
Then by \eqref{eq:RF}
\eqn
F(X_1,X_2)^{-1}R_{i,i}&\in&\CO_{i,j}{}^\times \dfrac{(\la_2-\la_1^{-1})^2}%
{(\la_2-q^2\la_1)(\la_2-q^{-2}\la_1)(\la_2-q^2\la_1^{-1})(\la_2-q^{-2}\la_1^{-1})}\\
&=&\CO_{i,j}{}^\times.
\eneqn

\bigskip
\noindent$\bullet$ Case $i=j$ and $\lambda(i)=\pm q^c$ for some $c=1,-1$.

In this case, $y_k^2=\psi(\la(i)^{-1}\la_k)$ ($k=1,2$), and
 \eqn
 R_{i,j}&=&-2(y_1^2+y_2^2)(y_1^2-y_2^2)^{-2}\\
 &=&-2\bl(\psi(\la(i)^{-1}\la_2)-\psi(\la(i)\la_1^{-1})\br)
 \bl(\psi(\la(i)^{-1}\la_1)-\psi(\la(i)^{-1}\la_2)\br)^{-2}\\
 &\in&\CO_{i,j}{}^\times(\la_2-q^{2c}\la_1^{-1})(\la_2-\la_1)^{-2}.
 \eneqn
Hence by \eqref{eq:RF}, we have \eqn
F(X_1, X_2)^{-1}R_{i,j}&\in&\CO_{i,j}{}^\times\dfrac{(\la_2-\la_1^{-1})^2}%
{(\la_2-q^2\la_1)(\la_2-q^{-2}\la_1)(\la_2-q^{-2c}\la_1^{-1})}
\eneqn
and it belongs to $\CO_{i,j}{}^\times$.

\bigskip
\noindent$\bullet$ Case $i=j$ and $\lambda(i)=\pm1$.

\smallskip
In this case,  $y_k=\psi(\la(i)^{-1}\la_k)^2$ ($k=1,2$),
and
\eqn
R_{i,j}=-(y_1-y_2)^{-2}&\in&\CO_{i,j}{}^\times
\bl(\psi(\la(i)^{-1}\la_2)-\psi(\la(i)^{-1}\la_1)\br)^{-2}
\bl(\psi(\la(i)^{-1}\la_2)-\psi(\la(i)\la_1^{-1})\br)^{-2}\\
&=&\CO_{i,j}{}^\times(\la_2-\la_1)^{-2}(\la_2-\la_1^{-1})^{-2}.
\eneqn
 Hence
\eqn
F(X_1,X_2)^{-1}R_{i,j}&\in&
\CO_{i,j}{}^\times\dfrac{1}%
{(\la_2-q^2\la_1)(\la_2-q^{-2}\la_1)(\la_2-q^2\la_1^{-1})(\la_2-q^{-2}\la_1^{-1})}
\eneqn and it belongs to $\CO_{i,j}^\times$.

\bigskip
Thus we have proved Lemma~\ref{lem:FR}, which completes the proof of
Theorem~\ref{th:isom_aDHC}. \qed

\bigskip
\subsection{The case $q^2=-1$}
Let us briefly explain what happens in the case $q^2=-1$. In this
case also, we can prove that  $\OHC$ is isomorphic to $\HRC$ by a
similar argument. We have  a factorization $K(u,v)=(u+v-4)(u+v+4)$.
Therefore $\pr\cl\T_X\to\AAA^1_w$ is the same but $g$ is the
identity. The Dynkin diagram structure of $\AAA^1_w$ is given by the
same formula \eqref{def:inner} as in the $(q^2)^2\not=1$ case. The
odd vertices are $w=\pm2$. We will explain which types
of Dynkin diagrams appear.

Let $s_1(w)=4-w$ and $s_2(w)=-4-w$ be the automorphisms of
$\AAA^1_w$.
Let $G$ be the group generated by $s_1$ and $s_2$.
Then it is an affine reflection group.
The connected components of $\AAA^1_w$ are nothing but the $G$-orbits.
They are described as follows
\bnum
\item
When $\CHAR\cor=0$, there are two types of  connected components.
\bna
\item $G\zeta$ ($\zeta\not\in 2\Z$). The Dynkin diagram is of type $A_\infty$.
\item $G\gamma$ with $\gamma=\pm2$. The Dynkin diagram is of type $B_\infty$.
\ee
\item
When $\CHAR\cor=p>2$, there are two types of  connected components.
Note that $(s_1s_2)^p=\id$.
\bna
\item $G\zeta$ ($\zeta\not\in \F_p$). The Dynkin diagram is of type $A^{(1)}_{2p-1}$.
\item $\F_p$. The Dynkin diagram is of type $D^{(2)}_{p}$.
\ee
\ee

\subsection{Cyclotomic Hecke-Clifford superalgebras}
For a subset $I_1$ of $\AAA^1_w$, set $J_1=\pr^{-1}(I_1)$.
For $\Lambda=\sum_{i\in I_1}m_i\Lambda_i$ with $m_i\in\Z_{\ge0}$
(see \S\;\ref{sec:cyclq}), we set
$$f_\Lambda(X_1)=\prod_{i\in J_1}(X_1-X(i))^{m_{\pr(i)}},$$
and
$$\AHCC_n^\Lambda=\AHC/\bl(\AHC f_\Lambda(X_1)\AHC\br).$$
We call $\AHCC_n^\Lambda$
the {\em cyclotomic Hecke-Clifford superalgebra}.
Then we have the following lemma.
\Lemma Let $n\ge1$ and $\Lambda$ as above.
\bnum
\item
$\AHCC_n^\Lambda$ is a finite-dimensional $\cor$-module.
\item
Set $I_n=\set{g(q^{2k}\la)}{\text{$g(\la)\in I_1$ and $-n<k<n$}}$
and $J_n=\pr^{-1}(I_n)$. Then for any finite-dimensional
$\AHCC_n^\Lambda$-module $M$, the eigenvalues of $X_k\vert_M$ belongs
to $\set{X(i)}{i\in J_n}$.
\ee
\enlemma

Hence by Theorem~\ref{th:isom_aDHC},we have

\Cor If $J$ is a subset of $\T_X$
invariant under $c$ and containing $J_n$, then the cyclotomic
Hecke-Clifford superalgebra $\AHCC_n^\Lambda$ is isomorphic to the
cyclotomic quiver Hecke-Clifford superalgebra $\RC_n^\Lambda$.
\encor

Note that the Hecke-Clifford algebra $\HC$ is isomorphic to
$\AHCC_n^{\Lambda_{i_0}}$, where $i_0\in  I$ is an odd vertex ($i_0=\pr(\pm1)=g(\pm q)$).
\bigskip

\section{Relations to affine Sergeev superalgebras}
In this section, we will prove that the affine Sergeev superalgebra,
a degenerate version of the affine Hecke-Clifford superalgebra, is
also isomorphic to a quiver Hecke-Clifford superalgebra. Since the
proof is parallel to the affine Hecke-Clifford superalgebra
case, the discussion will be brief. We still assume that
$\cor$ is an algebraically closed field of characteristic $\not=2$.

\subsection{Affine Sergeev superalgebras}

 \Def[{\cite{JN}}]\label{def_aDHC} For an integer $n\geq 0$,
the {\em affine Sergeev superalgebra} $\dAHC{n}$ is the
$\cor$-superalgebra generated by the even generators $x_1, \ldots,
x_n, t_1,\ldots, t_{n-1}$ and the odd generators $C_1,\ldots,C_n$
with the following defining relations. \bnum
\item $x_ix_j=x_jx_i$ for all $1\leq i,j\leq n$,
\item $C_i^2=1$, $C_iC_j+C_jC_i=0$ for all $1\leq i\ne j\leq n$,
\label{clifford_rel1}
\item $t_i^2=1$, $t_it_{i+1}t_i=t_{i+1}t_it_{i+1}$, $t_it_j=t_jt_i$ $(|i-j|\geq 2)$,
\item $t_iC_j=C_{s_i(j)}t_i$, \label{clifford_rel2}
\item $C_ix_j=x_jC_i$ for all $1\leq i\ne j\leq n$,
\item $C_ix_i=-x_iC_i$ for all $1\le i\le n$,\label{cond:C-}
\item $t_ix_i=x_{i+1}t_i-1-C_iC_{i+1}$, $t_ix_{i+1}=x_it_i+1-C_iC_{i+1}$ for all $1\leq i\leq n-1$,
\label{non_triv_eq d}
\item  $t_ix_j=x_jt_i$ if $j\ne i, i+1$.
\end{enumerate}
\label{def:deg_aff} \edf

\Rem\label{intertwin}(\cite{Naz},\cite[\S14.8]{Kle}) The
intertwiners of $\dAHC{n}$ are defined by
\eq&&\vphi_a=t_a+(x_a-x_{a+1})^{-1}-(x_a+x_{a+1})^{-1}C_aC_{a+1}
\in\cor(x_1,\ldots,x_n)\otimes_{\cor[x_1,\ldots,x_n]}\dAHC{n}\eneq %
and they satisfy the relations
\eq&&\vphi_ax_i=x_{s_a(i)}\vphi_a\quad\text{and}\quad
\vphi_aC_i=C_{s_a(i)}\vphi_a. \eneq %
Setting %
\eqn&& K(u,v)=(u-v)^2-2(u+v)\quad\text{and}\quad
F(x_a,x_{a+1})=K(x_a^2,x_{a+1}^2)/(x_a^2-x_{a+1}^2)^2 ,\eneqn%
we have \eq &&\Phi_a^2=F(x_a,x_{a+1}).\eneq%
 By setting $u=\lambda^2-1/4$ and $v=\mu^2-1/4$, we
 have
\begin{align}
K(u,v)= (\lambda-\mu-1)(\lambda-\mu+1)(\lambda+\mu-1)(\lambda+\mu+1).
\label{eq:Klam}
\end{align}

\enrem

\medskip
 Let $\AAA^1_x$ be
the one-dimensional affine space $\AAA^1$ with the coordinate $x$
and let $c\cl \AAA^1_x\to\AAA^1_x$ be the
involution of $\AAA^1_x$ given by $c(x)=-x$ (corresponding to
Definition \ref{def_aDHC} \eqref{cond:C-}).

 For $1\le p\le n$, let us denote by $c_p\cl \AAA^n\to\AAA^n$ the
involution $$(x_1,\ldots,x_n)\mapsto
(x_1,\ldots,x_{p-1},-x_p,x_{p+1},\ldots,x_n),$$ and by
the same letter the induced isomorphism
$\FRAC({\HAT{\CO}}_{{\AAA}^n,q})\To[{c_p}]\FRAC({\HAT{\CO}}_{{\AAA}^n,c_p(q)})$.
For $1\le p<n$, we denote by $s_p\cl \AAA^n\to\AAA^n$ the involution
$$(x_1,\ldots,x_n)\mapsto
(x_1,\ldots,x_{p-1},x_{p+1},x_{p},x_{p+2},\ldots,x_n),$$ and by
the same letter $s_p$ the induced isomorphism
$\FRAC({\HAT{\CO}}_{{\AAA}^n,q})\to\FRAC({\HAT{\CO}}_{{\AAA}^n,s_p(q)})$.
Similarly, we denote by
$\bs_p\cl \AAA^n\to\AAA^n$ the
involution
$$(x_1,\ldots,x_n)\mapsto
(x_1,\ldots,x_{p-1},-x_{p+1},-x_{p},x_{p+2},\ldots,x_n),$$ and by
the same letter $\bs_p$ the induced isomorphism
$\FRAC({\HAT{\CO}}_{{\AAA}^n,q})\to\FRAC({\HAT{\CO}}_{{\AAA}^n,\bs_p(q)})$.

\medskip
Let us denote by $\hc$ the $\cor$-superalgebra generated by $C_p$
($1\le p\le n$), $t_a$ ($1\le a<n$)
 with the defining relations \eqref{clifford_rel1}--\eqref{clifford_rel2} in
 Definition~\ref{def_aDHC}.
The superalgebra $\hc$ is called the {\em Sergeev superalgebra} and
can be regarded as a subsuperalgebra of $\dAHC{n}$.

 \Def\label{comphc}Let $J$ be a finite subset of $\AAA^1_x$ invariant by $c$,
 and let $x\cl J\to\T_X$ be the inclusion map.
 We define
\begin{align*}
\CO_{n} = \bigoplus_{\nu\in J^n}
{\HAT{\CO}}_{\AAA^n,x(\nu)}e(\nu),\quad \CK_n = \bigoplus_{\nu\in
J^n} \FRAC({\HAT{\CO}}_{{\AAA}^n,x(\nu)})e(\nu),
\end{align*}
where $x(\nu)=(x(\nu_1),\cdots,x(\nu_n))\in \AAA^n$. We define the
algebra structure on
\begin{align*}
\Khc=\CK_n\otimes_{\cor}\hc
\end{align*}
by
\begin{equation}\label{def:tformald}
\begin{aligned}
C_pe(\nu)f & =e(\CC_p\nu)c_p(f)C_p \ \ (1\leq p\leq n), \\
t_ae(\nu)f & = e(s_a\nu)s_a(f)t_a +
\frac{e(\nu)f-e(s_a\nu)s_af}{x_{a+1}-x_a}
+C_aC_{a+1}\dfrac{e(\nu)f-e(\bs_a\nu)\bs_a(f)}{x_{a+1}+x_a}
\end{aligned}
\end{equation}
for $\nu\in J^n$, $1\leq a<n$ and $f\in
\FRAC({\HAT{\CO}}_{{\AAA}^n,x(\nu)})$.

We define $\Ohc$ to be the subsuperalgebra of
$\Khc$ generated by $\CO_n$ and $\hc$. \edf

Thus there exists a $\cor$-superalgebra homomorphism
$\dAHC{n}\to\Ohc$.
\medskip
We have
\eq
&&\text{$\vphi_a e(\nu)=e(s_a\nu)\vphi_a$ for any $1\le a<n$.}
\eneq

\medskip

Let us denote by $\AAA^1_w$ the space $\AAA^1$ with
the coordinate $w$. Let $\pr\cl \AAA^1_x\to\AAA^1_w$ be the map
$\pr(x)=x^2$. Let $\AAA^1_\la$ be the affine space $\AAA^1$ with the
coordinate $\la$. Let $g\cl\AAA^1_\la\to\AAA^1_w$ be the map
$g(\la)=\la^2-1/4$.

\medskip
\subsection{Dynkin diagram}\label{sec:Dyndahc}
We  regard $\AAA^1_w$ as a Dynkin diagram
by the same formula \eqref{def:inner}.
Note that $g(\la_1)$ and $g(\la_2)$ ($\la_1\not=\pm\la_2$) are connected
if and only if $\la_1=\la_2\pm1$ or $\la_1=-\la_2\pm1$.

The connected components of the Dynkin diagram $\AAA^1_w$ are
classified as follows. Here the odd vertices are marked
by $\times$.

\bnum \item when $\CHAR \cor=0$, there are three types of
connected components.

\vs{1ex}
\begin{enumerate}[(a)]
\item $\set{g(\zeta+k)}{k\in\Z}$ for some
$\zeta\not\in\Z/2$. The Dynkin diagram is of type $A_{\infty}$.

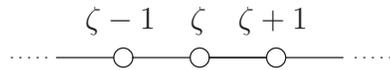
\vs{2ex}
\begin{center}
\unitlength 0.1in
\begin{picture}( 20.0000,  3.0000)(  0.0000, -2.5000)
%
{\color[named]{Black}{%
\special{pn 8}%
\special{ar 1000 200 50 50  0.0000000 6.2831853}%
}}%
%
{\color[named]{Black}{%
\special{pn 8}%
\special{ar 1400 200 50 50  0.0000000 6.2831853}%
}}%
%
{\color[named]{Black}{%
\special{pn 8}%
\special{ar 600 200 50 50  0.0000000 6.2831853}%
}}%
%
{\color[named]{Black}{%
\special{pn 8}%
\special{pa 650 200}%
\special{pa 950 200}%
\special{fp}%
}}%
%
{\color[named]{Black}{%
\special{pn 8}%
\special{pa 1050 200}%
\special{pa 1350 200}%
\special{fp}%
}}%
%
{\color[named]{Black}{%
\special{pn 8}%
\special{pa 1450 200}%
\special{pa 1750 200}%
\special{fp}%
}}%
%
{\color[named]{Black}{%
\special{pn 8}%
\special{pa 1050 200}%
\special{pa 1350 200}%
\special{fp}%
}}%
%
{\color[named]{Black}{%
\special{pn 8}%
\special{pa 250 200}%
\special{pa 550 200}%
\special{fp}%
}}%
\put(9.5000,-0.8000){\makebox(0,0)[lb]{$\zeta$}}%
\put(12.0000,-0.8000){\makebox(0,0)[lb]{$\zeta+1$}}%
\put(4.0000,-0.8000){\makebox(0,0)[lb]{$\zeta-1$}}%
%
{\color[named]{Black}{%
\special{pn 8}%
\special{pa 200 200}%
\special{pa 0 200}%
\special{dt 0.045}%
}}%
%
{\color[named]{Black}{%
\special{pn 8}%
\special{pa 2000 200}%
\special{pa 1800 200}%
\special{dt 0.045}%
}}%
\end{picture}%

\end{center}

\vs{3ex}

\item $\set{g(\frac{1}{2}+k)}{k\in\Z_{\geq 0}}$.
The Dynkin diagram is of type $B_\infty$.

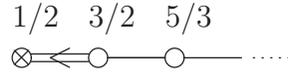
\medskip
\begin{center}
\unitlength 0.1in
\begin{picture}( 14.5000,  3.0000)(  1.5000, -2.5000)
%
{\color[named]{Black}{%
\special{pn 8}%
\special{ar 600 200 50 50  0.0000000 6.2831853}%
}}%
%
{\color[named]{Black}{%
\special{pn 8}%
\special{ar 1000 200 50 50  0.0000000 6.2831853}%
}}%
%
{\color[named]{Black}{%
\special{pn 8}%
\special{ar 200 200 50 50  0.0000000 6.2831853}%
}}%
%
{\color[named]{Black}{%
\special{pn 8}%
\special{pa 650 200}%
\special{pa 950 200}%
\special{fp}%
}}%
%
{\color[named]{Black}{%
\special{pn 8}%
\special{pa 1050 200}%
\special{pa 1350 200}%
\special{fp}%
}}%
%
{\color[named]{Black}{%
\special{pn 8}%
\special{pa 650 200}%
\special{pa 950 200}%
\special{fp}%
}}%
\put(9.5000,-0.8000){\makebox(0,0)[lb]{$5/3$}}%
%
{\color[named]{Black}{%
\special{pn 8}%
\special{pa 1600 200}%
\special{pa 1400 200}%
\special{dt 0.045}%
}}%
\put(5.5000,-0.8000){\makebox(0,0)[lb]{$3/2$}}%
\put(1.5000,-0.8000){\makebox(0,0)[lb]{$1/2$}}%
\put(1.350,-2.54000){\makebox(0,0)[lb]{$\times$}}%
%
{\color[named]{Black}{%
\special{pn 8}%
\special{pa 250 180}%
\special{pa 550 180}%
\special{fp}%
}}%
%
{\color[named]{Black}{%
\special{pn 8}%
\special{pa 250 220}%
\special{pa 550 220}%
\special{fp}%
}}%
%
{\color[named]{Black}{%
\special{pn 8}%
\special{pa 350 200}%
\special{pa 450 150}%
\special{fp}%
}}%
%
{\color[named]{Black}{%
\special{pn 8}%
\special{pa 450 250}%
\special{pa 350 200}%
\special{fp}%
}}%
\end{picture}%

\end{center}

\vs{3ex}
\item $\set{g(k)}{k\in\Z_{\geq 0}}$.
The Dynkin diagram is of type $C_\infty$.

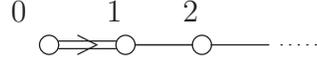
\vs{3ex}
\begin{center}
\unitlength 0.1in
\begin{picture}( 16.0000,  3.0000)(  0.0000, -2.5000)
%
{\color[named]{Black}{%
\special{pn 8}%
\special{ar 600 200 50 50  0.0000000 6.2831853}%
}}%
%
{\color[named]{Black}{%
\special{pn 8}%
\special{ar 1000 200 50 50  0.0000000 6.2831853}%
}}%
%
{\color[named]{Black}{%
\special{pn 8}%
\special{ar 200 200 50 50  0.0000000 6.2831853}%
}}%
%
{\color[named]{Black}{%
\special{pn 8}%
\special{pa 650 200}%
\special{pa 950 200}%
\special{fp}%
}}%
%
{\color[named]{Black}{%
\special{pn 8}%
\special{pa 1050 200}%
\special{pa 1350 200}%
\special{fp}%
}}%
%
{\color[named]{Black}{%
\special{pn 8}%
\special{pa 650 200}%
\special{pa 950 200}%
\special{fp}%
}}%
\put(9.0000,-0.8000){\makebox(0,0)[lb]{$2$}}%
%
{\color[named]{Black}{%
\special{pn 8}%
\special{pa 1600 200}%
\special{pa 1400 200}%
\special{dt 0.045}%
}}%
\put(5.0000,-0.8000){\makebox(0,0)[lb]{$1$}}%
\put(0.0000,-0.8000){\makebox(0,0)[lb]{$0$}}%
%
{\color[named]{Black}{%
\special{pn 8}%
\special{pa 250 180}%
\special{pa 550 180}%
\special{fp}%
}}%
%
{\color[named]{Black}{%
\special{pn 8}%
\special{pa 250 220}%
\special{pa 550 220}%
\special{fp}%
}}%
%
{\color[named]{Black}{%
\special{pn 8}%
\special{pa 350 150}%
\special{pa 450 200}%
\special{fp}%
}}%
%
{\color[named]{Black}{%
\special{pn 8}%
\special{pa 450 200}%
\special{pa 350 250}%
\special{fp}%
}}%
\end{picture}%

\end{center}
\end{enumerate}

\vs{3ex}
\item  When $\CHAR \cor=p>2$, there are
three types of connected components.

\vs{1ex}
\begin{enumerate}
\item $\set{g(\zeta+k)}{k\in\F_{p}}$ for some $\zeta\not\in\F_{p}$. The Dynkin diagram is of type
$A^{(1)}_{p-1}$.

\medskip
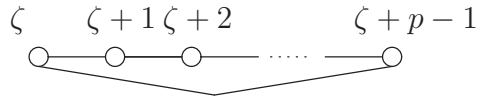
\begin{center}
\unitlength 0.1in
\begin{picture}( 20.5000,  4.5000)(  0.5000, -4.0000)
%
{\color[named]{Black}{%
\special{pn 8}%
\special{ar 600 200 50 50  0.0000000 6.2831853}%
}}%
%
{\color[named]{Black}{%
\special{pn 8}%
\special{ar 1000 200 50 50  0.0000000 6.2831853}%
}}%
%
{\color[named]{Black}{%
\special{pn 8}%
\special{ar 200 200 50 50  0.0000000 6.2831853}%
}}%
%
{\color[named]{Black}{%
\special{pn 8}%
\special{pa 250 200}%
\special{pa 550 200}%
\special{fp}%
}}%
%
{\color[named]{Black}{%
\special{pn 8}%
\special{pa 650 200}%
\special{pa 950 200}%
\special{fp}%
}}%
%
{\color[named]{Black}{%
\special{pn 8}%
\special{pa 1050 200}%
\special{pa 1350 200}%
\special{fp}%
}}%
%
{\color[named]{Black}{%
\special{pn 8}%
\special{pa 650 200}%
\special{pa 950 200}%
\special{fp}%
}}%
\put(8.5000,-0.8000){\makebox(0,0)[lb]{$\zeta+2$}}%
%
{\color[named]{Black}{%
\special{pn 8}%
\special{pa 1600 200}%
\special{pa 1400 200}%
\special{dt 0.045}%
}}%
\put(4.5000,-0.8000){\makebox(0,0)[lb]{$\zeta+1$}}%
\put(0.5000,-0.8000){\makebox(0,0)[lb]{$\zeta$}}%
%
{\color[named]{Black}{%
\special{pn 8}%
\special{pa 1700 200}%
\special{pa 2000 200}%
\special{fp}%
}}%
%
{\color[named]{Black}{%
\special{pn 8}%
\special{ar 2050 200 50 50  0.0000000 6.2831853}%
}}%
\put(18.5000,-0.8000){\makebox(0,0)[lb]{$\zeta+p-1$}}%
%
{\color[named]{Black}{%
\special{pn 8}%
\special{pa 200 250}%
\special{pa 1120 400}%
\special{fp}%
\special{pa 1120 400}%
\special{pa 2050 250}%
\special{fp}%
}}%
\end{picture}%

\end{center}
\label{case:loop}

\medskip
\item $\{g(\frac{1}{2}),\ldots,g(\frac{p}{2})\}$ of type $A^{(2)}_{p-1}$.

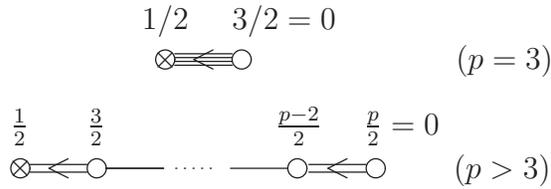
\vs{2ex}
\begin{center}
\hs{10ex}
\unitlength 0.1in
\begin{picture}(  5.0000,  3.0000)(  1.5000, -2.5000)
{\color[named]{Black}{%
\special{pn 8}%
\special{ar 600 200 50 50  0.0000000 6.2831853}%
}}%
{\color[named]{Black}{%
\special{pn 8}%
\special{ar 200 200 50 50  0.0000000 6.2831853}%
}}%
{\color[named]{Black}{%
\special{pn 8}%
\special{pa 250 170}%
\special{pa 550 170}%
\special{fp}%
}}%
\put(5.5000,-0.8000){\makebox(0,0)[lb]{$3/2=0$}}%
\put(1.5000,-0.8000){\makebox(0,0)[lb]{\kern-1ex$1/2$}}%
\put(1.3500,-2.54000){\makebox(0,0)[lb]{$\times$}}
{\color[named]{Black}{%
\special{pn 8}%
\special{pa 250 230}%
\special{pa 550 230}%
\special{fp}%
}}%
{\color[named]{Black}{%
\special{pn 8}%
\special{pa 350 200}%
\special{pa 450 150}%
\special{fp}%
}}%
{\color[named]{Black}{%
\special{pn 8}%
\special{pa 450 250}%
\special{pa 350 200}%
\special{fp}%
}}%
{\color[named]{Black}{%
\special{pn 8}%
\special{pa 250 190}%
\special{pa 550 190}%
\special{fp}%
}}%
{\color[named]{Black}{%
\special{pn 8}%
\special{pa 250 210}%
\special{pa 550 210}%
\special{fp}%
}}%
\end{picture}%
\hs{15ex}\text{($p=3$)}
\end{center}

\vs{3ex}
\begin{center}
\unitlength 0.1in
\begin{picture}( 19.6000,  3.0000)(  1.5000, -2.5000)
%
{\color[named]{Black}{%
\special{pn 8}%
\special{ar 600 200 50 50  0.0000000 6.2831853}%
}}%
%
{\color[named]{Black}{%
\special{pn 8}%
\special{ar 200 200 50 50  0.0000000 6.2831853}%
}}%
%
{\color[named]{Black}{%
\special{pn 8}%
\special{pa 250 180}%
\special{pa 550 180}%
\special{fp}%
}}%
%
{\color[named]{Black}{%
\special{pn 8}%
\special{pa 650 200}%
\special{pa 950 200}%
\special{fp}%
}}%
%
{\color[named]{Black}{%
\special{pn 8}%
\special{pa 650 200}%
\special{pa 950 200}%
\special{fp}%
}}%
%
{\color[named]{Black}{%
\special{pn 8}%
\special{pa 1200 200}%
\special{pa 1000 200}%
\special{dt 0.045}%
}}%
\put(5.5000,-0.8000){\makebox(0,0)[lb]{$\frac{3}{2}$}}%
\put(1.5000,-0.8000){\makebox(0,0)[lb]{$\frac{1}{2}$}}%
\put(1.3500,-2.54000){\makebox(0,0)[lb]{$\times$}}%
%
{\color[named]{Black}{%
\special{pn 8}%
\special{pa 1300 200}%
\special{pa 1600 200}%
\special{fp}%
}}%
%
{\color[named]{Black}{%
\special{pn 8}%
\special{ar 1650 200 50 50  0.0000000 6.2831853}%
}}%
\put(15.0000,-0.8000){\makebox(0,0)[lb]{\kern0.5ex$\frac{p-2}{2}$}}%
%
{\color[named]{Black}{%
\special{pn 8}%
\special{pa 250 220}%
\special{pa 550 220}%
\special{fp}%
}}%
%
{\color[named]{Black}{%
\special{pn 8}%
\special{pa 350 200}%
\special{pa 450 150}%
\special{fp}%
}}%
%
{\color[named]{Black}{%
\special{pn 8}%
\special{pa 450 250}%
\special{pa 350 200}%
\special{fp}%
}}%
%
{\color[named]{Black}{%
\special{pn 8}%
\special{ar 2060 200 50 50  0.0000000 6.2831853}%
}}%
%
{\color[named]{Black}{%
\special{pn 8}%
\special{pa 1710 180}%
\special{pa 2010 180}%
\special{fp}%
}}%
%
{\color[named]{Black}{%
\special{pn 8}%
\special{pa 1710 220}%
\special{pa 2010 220}%
\special{fp}%
}}%

{\color[named]{Black}{%
\special{pn 8}%
\special{pa 1810 200}%
\special{pa 1910 150}%
\special{fp}%
}}%

{\color[named]{Black}{%
\special{pn 8}%
\special{pa 1910 250}%
\special{pa 1810 200}%
\special{fp}%
}}%

\put(20.0000,-0.8000){\makebox(0,0)[lb]{$\frac{p}{2}=0$}}%
\end{picture}%
\hs{5ex}\text{($p>3$)}
\end{center}
\end{enumerate}
\ee

\bigskip
Let us take a finite subset $J\subset \AAA^1_x$ invariant under the
involution $c$ and set $I=\pr(J)\subset\AAA^1_w$. Set
$J^c\seteq\set{j\in J}{c(j)=j}$, $I_\odd=\pr(J^c)$ and
$I_\even=I\setminus\{0\}$.
 We choose a function $h\cl I\to \AAA^1_\la$
 such that
 \eq
 &&\left\{\parbox{55ex}{\be[{\small$\bullet$}]
 \item $g(h(i))\seteq h(i)^2-1/4=i$,

\vs{.0ex}
 \item if $(\alpha_i,\alpha_{i'})<0$, then $h(i')=h(i)\pm1$,

\vs{.0ex}
 \item $h(g(1/2))=1/2$.
 \ee}\right.
 \label{cond:hd}\eneq

\medskip
 We define
$\lambda\cl J\to \AAA^1_\lambda $ by $\lambda(j)=h(\pr(j))$, and
$x \cl J\to\AAA^1_x$ by $x(j)=j$.
We choose $\eps\cl J\to\{0,1\}$ such that $\eps^{-1}(0)\to I$ is
bijective. For $i, j \in I$, we write $(i,j)$ for $(\alpha_{\pr(i)},
\alpha_{\pr(j)})$ and define $\tQ_{i,j}$ by the same formula
\eqref{def:tQ} and \eqref{eq:extQ}.

Then we associate to the data $\{\tQ_{i,j}\}_{i,j\in J}$ the quiver Hecke-Clifford superalgebra $\RC_n$
and its completion $\widehat\RC_n$.

 \Th\label{th:isom_aDhc} We have $\widehat{\RC}_n \simeq\Ohc$. \enth

\bigskip
\subsection{Proof of Theorem~\ref{th:isom_aDhc}}
\subsubsection{Strategy of the proof}
By Theorem~\ref{poly_rep}, it is sufficient to construct the
elements $y_k\in\CO_n$ ($k=1,\ldots,n$) and $\ts_a\in\Khc$
($a=1,\ldots,n-1$)
 such that
\eq&&\left\{
\parbox{73ex}{
\bnum
\item $y_pe(\nu)\in
({\HCO}_{\AAA^n,x(\nu)})^{\times}\,(x_p-x(\nu_p))e(\nu)$,

\vs{.3ex}
\item $C_ay_pe(\nu)=(-1)^{\delta{a,p}}y_pe(c_a\nu)C_a$,

\vs{1ex}
\item $\ts_a^2=R_{a,a+1}$ (see \eqref{def R}),

\vs{1ex}
\item $\{\ts_a\}_{1\le a<n}$ satisfies the braid relations (see \eqref{eq:braid}),

\vs{1ex}
\item $\ts_aC_k=C_{s_a(k)}\ts_a$, $\ts_ay_k=y_{s_a(k)}\ts_a$,

\vs{.3ex}

\item setting $\sigma_a\seteq\ts_a+f_{a,a+1}$ (see \eqref{def:f}),
the element $\sigma_ae(\nu)$ belongs to
$t_ae(\nu)({\HCO}_{\AAA^n,x(\nu)})^{\times}
+\CO_n\langle C_1,\ldots,C_n\rangle$.\label{eq:taureg}
\ee}\right.
\label{eq:hcrcd}
\eneq

\subsubsection{Definition of $y_p\,e(\nu)$}\quad

We will construct $y_pe(\nu)\in
({\HCO}_{\AAA^n,x(\nu)})^{\times}\,(x_p-x(\nu_p))e(\nu)$.
Let $C$ be the curve defined by the equation
$x^2=\lambda^2-1/4$. 
We have two projections:
\eqn
&&\AAA^1_x\from C\To \AAA^1_\la.
\eneqn
For every $i\in J$, they induce
the injective homomorphisms
\eq
&&\HCO_{\AAA^1, x(i)}\monoto\HCO_{C,(x(i),\la(i))}\monofrom\HCO_{\AAA^1, \la(i)}.
\eneq
Let us set $\CO_i=\HCO_{\AAA^1,x(i)}$
and $\tCO_i=\HCO_{C,(x(i),\la(i))}$,
and we regard $\CO_i$ as a subalgebra of $\tCO_i$.
Let $c_i\cl\CO_i\to\CO_{c(i)}$ be the homomorphism induced by the map $x\mapsto-x$.
For $1\le a\le n$, let us denote by $q_{a,\nu}\cl\CO_{\nu_a}\To\HCO_{\AAA^n,x(\nu)}e(\nu)$
the homomorphism induced by the $a$-th projection $(x_1,\ldots,x_n)\mapsto x_a$.

\medskip
For all $i\in J$ such that $\eps(i)=0$, we will construct
$y_i\in\CO_i$ such that
\eq
&&\left\{\parbox{60ex}{
\bnum
\item
$y_i\in \CO_i^\times(x-x(i))$,

\vs{.5ex}
\item if $c(i)=i$, then
$c_i(y_i)=-y_i$, \ee }\right. \label{cond:yd} \eneq and we define
$y_ae(\nu)=q_{a,\nu}((-c_{\nu_a})^{-\eps(\nu_a)}y_{c^{\eps(\nu_a)}\nu_a})$.

Then we have
$$
\begin{aligned}
& y_pe(\nu)\in
({\HCO}_{\AAA^n,x(\nu)})^{\times}\,(x_p-x(\nu_p))e(\nu), \\
& C_a y_pe(\nu)=(-1)^{\delta_{ap}}y_pe(c_a\nu)C_a.
\end{aligned}
$$

\bigskip
\noindent
$\bullet$ Case $\lambda(i)\ne 0,1/2$.

In this case, we have isomorphisms
\begin{align*}
\cor[[x-x(i)]]=\CO_i\isoto\tCO_i=\HAT{\CO}_{C,\,(\lambda(i),x(i))}
\isofrom
\HAT{\CO}_{\AAA^1_{\lambda},\lambda(i)}=\cor[[\lambda-\lambda(i)]].
\end{align*}
We define $y_i\in\CO_{i}$ by $y_i=\lambda-\lambda(i)$.

\bigskip
\noindent
$\bullet$ Case $\lambda(i)=1/2$.

\smallskip
In this case, we have $x(i)=0$,
\begin{align*}
\cor[[x]]=\CO_i=\HAT{\CO}_{\AAA^1_{x},0}
 \isoto
\tCO_i=\HAT{\CO}_{C,\,(0,1/2)}
\monofrom\HAT{\CO}_{\AAA^1_{\lambda},1/2}=\cor[[\lambda-1/2]],
\end{align*}
and
$\lambda-\frac{1}{2}=\frac{1}{2}(\sqrt{1+4x^2}-1)=x^2+\cdots\in \cor[[x^2]]$.
Take $y_i=x+\cdots\in \cor[[x]]=\CO_i$ such that $y_i^2=\lambda$.
Since $y_i\in x\cor[[x^2]]$, we have
$c_i(y_i)=-y_i$.

\bigskip
\noindent
$\bullet$ Case $\lambda(i)=0$.

\smallskip
In this case, we have $x(i)^2=-1/4$ and
\begin{align*}
\cor[[x-x(i)]]=\HAT{\CO}_{\AAA^1_{x},x(i)}
\monoto \HAT{\CO}_{C,\,(x(i),0)}
\isofrom\HAT{\CO}_{\AAA^1_{\lambda},\,0}=\cor[[\lambda]].
\end{align*}
Since $x^2-x(i)^2=\la^2$, we have
$\CO_i=\cor[[x-x(i)]]=
\cor[[\la^2]]$.
We define $y_i\in \cor[[x-x(i)]]$ by
$y_i=\la^2$.

\subsubsection{Definition of $\ts_ae(\nu)$}
\quad We will define $\ts_ae(\nu)$ for $1\le a<n$ and $\nu\in J^n$.

For $i,j\in J$, let us denote $\CO_{i,j}=\HCO_{\AAA^2,(x(i),x(j))}$
and $\tCO_{i,j}=\HCO_{C\times C,((x(i),\la(i)),(x(j),\la(j)))}$. We
regard $\CO_{i,j}$ as a subalgebra of $\tCO_{i,j}$. Let
$r^1_{i,j}\cl\tCO_i\to \tCO_{i,j}$ and $r^2_{i,j}\cl\tCO_j\to
\tCO_{i,j}$ be the algebra homomorphism induced by the first and
second projections from $C\times C$ to $C$, respectively. We write
$y_1$ for $r^1_{i,j}(y_i)\in\CO_{i,j}$ and $y_2$ for
$r^2_{i,j}(y_j)\in\CO_{i,j}$. Similarly, we define
$\la_1\in\tCO_{i,j}$ and $\la_2\in\tCO_{i,j}$ as $r^1_{i,j}(\la)$ and
$r^2_{i,j}(\la)$.

Let $c_1\cl\CO_{i,j}\to\CO_{c(i),j}$ and $c_2\cl\CO_{i,j}\to\CO_{i,c(j)}$
be the isomorphisms induced by $(x_1,x_2)\mapsto (-x_1,x_2)$ and $(x_1,x_2)\mapsto (x_1,-x_2)$,
respectively.
Let $s_{12}\cl \CO_{i,j}\to\CO_{j,i}$ be the homomorphism induced by $\AAA^2\ni(x_1,x_2)\mapsto(x_2,x_1)\in\AAA^2$.

Let $\psi_{a,\nu}\cl\CO_{\nu_a,\nu_{a+1}}\To\HCO_{\AAA^n,x(\nu)}e(\nu)$
be the algebra homomorphism induced by
the projection $(x_1,\ldots,x_n)\longmapsto(x_a,x_{a+1})$.
We define $R_{i,j}\in\CO_{i,j}$ by the same formula \eqref{eq:R}.
Let us recall (see Remark~\ref{intertwin})
\begin{align*}
\vphi_a^2
&=F(x_a,x_{a+1}).
\end{align*}

\Lemma\label{lem:fr}
For $i,j\in J$ such that $\eps(i)=\eps(j)=0$, we have
\bnum
\item $F(x_1,x_2)^{-1}R_{i,j}$ belongs to $\CO_{i,j}{}^\times$,
\item if $i=j$, then
$F(x_1,x_2)^{-1}R_{i,j}-\left(\dfrac{x_1-x_2}{y_1-y_2}\right)^2$
belongs to $\CO_{i,j}(x_1-x_2)$.
\ee
\enlemma

Admitting this lemma for a while, we will construct $\ts_ae(\nu)$
and prove the theorem. By this lemma, we  can choose
$G_{i,j}\in\CO_{i,j}$ for $i,j\in J$, satisfying the following
conditions: \eq&&\left\{
\parbox{70ex}{
\bna
\item
$G_{i,j}\cdot s_{12}G_{j,i} = F(x_1,x_2)^{-1}R_{i,j}$,\label{eq:gsg}
\vs{1ex}
\item if $i=j$, then $G_{i,j}-\dfrac{x_1-x_2}{y_1-y_2}\in(x_1-x_2)\CO_{i,j}$,
\label{cond:gtau}
\vs{.5ex}
\item $c_1\bl(G_{i,j}\br)=G_{c(i),j}$ and
$c_2\bl(G_{i,j}\br)=G_{i,c(j)}$.\label{cond:cg}
\ee
}\right.\label{cond:Gd}
\eneq
Now we define
\eq&&\ts_ae(\nu)=\vphi _ae(\nu)\psi_{a,\nu}(G_{\nu_{a},\nu_{a+1}}).\eneq
By the construction, we have
\eqn
&&\text{$\ts_a\cdot e(\nu)f=s_a(f)e(s_a\nu)\cdot\ts_a$ for any $f\in \FRAC({\HAT{\CO}}_{{\AAA}^n,x(\nu)})$,}\\
&&\ts_aC_b=C_{s_a(b)}\ts_a,\\
&&\text{$\ts_a^2e(\nu)=R_{a,a+1}e(\nu)$.}
\eneqn

\noi Let us verify
 the condition \eqref{eq:hcrcd}\;\eqref{eq:taureg}:
$$\sigma_ae(\nu)\seteq\ts_ae(\nu)+f_{a,a+1}e(\nu)\in
t_ae(\nu)(\HCO_{\AAA^n,x(\nu)})^{\times}+\CO_n\langle
C_1,\cdots,C_n\rangle.$$
Setting $i=\nu_a$ and $j=\nu_{a+1}$ and
 writing $G_{i,j}$ for $\psi_{a,\nu}(G_{\nu_{a},\nu_{a+1}})$, we have
 \eqn
\tau_ae(\nu)&=&\Bigl(\ts_a-(y_a-y_{a+1})\delta_{i,j}+(y_a+y_{a+1})^{-1}C_aC_{a+1}\delta_{c(i),j}\Bigr)e(\nu)\\
&=&\Bigl(t_a+(x_a-x_{a+1})^{-1}-(x_a+x_{a+1})^{-1}C_aC_{a+1}\Bigr)e(\nu)G_{i,j}\\
&&\hs{10ex}-(y_a-y_{a+1})\delta_{i,j}e(\nu)+(y_a+y_{a+1})^{-1}C_aC_{a+1}\delta_{c(i),j}e(\nu)\\
&=&t_ae(\nu)G_{i,j}+\Bigl(G_{i,j}(x_a-x_{a+1})^{-1}-(y_a-y_{a+1})\delta_{i,j}\Bigr)e(\nu)\\
&&\hs{10ex}+C_a\Bigl(-(c_{a+1}G_{i,j})(-x_a+x_{a+1})^{-1}+(-y_a+y_{a+1})^{-1}\delta_{c(i),j}\Bigr)C_{a+1}e(\nu)\\
&=&t_ae(\nu)G_{i,j}+\Bigl(G_{i,j}(x_a-x_{a+1})^{-1}-(y_a-y_{a+1})\delta_{i,j}\Bigr)e(\nu)\\
&&\hs{10ex}+C_{a}\Bigl(G_{i,c(j)}(x_a-x_{a+1})^{-1}-(y_a-y_{a+1})^{-1}\delta_{i,c(j)}\Bigr)e(c_{a+1}\nu)C_{a+1}
\eneqn
belongs to $t_ae(\nu)(\HCO_{\AAA^n,x(\nu)})^{\times}+\CO_n\langle C_1,\cdots,C_n\rangle$
by \eqref{cond:Gd}\;\eqref{cond:gtau}.
Thus we have verified all the conditions in \eqref{eq:hcrcd},
and hence $\widehat{\RC}_n $ and $\Ohc$ are isomorphic.

\medskip
 \subsubsection{Proof of Lemma~\ref{lem:fr}}
 Now, we will prove Lemma~\ref{lem:fr} for $i,j\in J$ such that $\eps(i)=\eps(j)=0$.

 \smallskip
 Let us first derive (ii) admitting (i).
 Let $\xi$ be the generic point of  $\{x_1=x_2\}\subset\on{Spec}(\CO_{i,j})$, and let $(\CO_{i,j})_\xi$ be
 the localization of $\CO_{i,j}$ at $\xi$.
 By (i),
$F(x_1,x_{2})^{-1}R_{i,j}-\Bigl((x_1-x_2)/(y_1-y_2)\Bigr)^{2}$
belongs to $\CO_{i,j}$, and hence  it is enough to show that
$F(x_1,x_{2})^{-1}R_{i,j}-\Bigl((x_1-x_2)/(y_1-y_2)\Bigr)^{2}$
belongs to $(\CO_{i,j})_\xi(x_1-x_2)$.

We have
$$F(x_1,x_{2})^{-1}(x_1-x_2)^{-2}=\dfrac{(x_1+x_{2})^2}%
{(x_1^2-x_{2}^2)^2-2(x_1^2+x_{2}^2)^2}\equiv-1 \mod (\CO_{i,j})_\xi(x_1-x_2).$$
Since $R_{i,j}\equiv -(y_1-y_2)^{-2} \mod (\CO_{i,j})_\xi$,
we obtain
$$F(x_1,x_{2})^{-1}R_{i,j}\equiv \dfrac{(x_1-x_2)^2}{(y_1-y_2)^2} \mod (\CO_{i,j})_\xi(x_1-x_2).$$

It only remains to prove (i). We will use case-by-case check-up.
Recall that \eq&&\ba{rl} F(x_1,x_{2})^{-1}
&=\dfrac{(\la_1-\la_{2})^2(\la_1+\la_{2})^2}%
{(\lambda_1-\lambda_{2}-1)(\lambda_1-\lambda_{2}+1)(\lambda_1+\lambda_{2}-1)(\lambda_1+\lambda_{2}+1)}
.
\ea\label{eq:Rf}
\eneq

Note that (see \eqref{cond:hd}) for $i,j\in J$ such that
$\eps(i)=\eps(j)=0$, we have
$$
\parbox{60ex}
{
\bna
\item $\la(j)=\pm\la(i)$ implies $i=j$,

\vs{0.5ex}
\item $(\alpha_{\pr(i)},\alpha_{\pr(j)})<0$
implies $\la(j)=\la(i)\pm1$,

\vs{0.5ex}
\item
$\la(g(1/2))=1/2$.
\ee}
$$

\bigskip
\noindent$\bullet$ Case $\lambda(j)\ne \lambda(i),
\lambda(i)\pm 1$.

\smallskip
In this case, $F(x_1,x_{2})\in\CO^{\times}_{i,j}$ and $R_{i,j}=1$.

\bigskip
\noindent
$\bullet$ Case $\lambda(j)=\lambda(i)\pm1$ and
$\lambda(i),\lambda(j)\ne 0,1/2$.

\smallskip
Set $\la(j)=\la(i)+c$ with $c=\pm1$. In this case,
$y_1=\la_1-\la(i)$, $y_2=\la_2-\la(j)$, and
$R_{i,j}=\tQ_{i,j}(y_1,y_2)=\pm(y_1-y_2)=\pm(\la_1-\la_2+c)$.
Therefore,
\begin{align*}
F(x_1,x_2)^{-1}R_{i,j}
&=\pm\dfrac{(\la_1-\la_{2})^2(\la_1+\la_{2})^2}%
{(\lambda_1-\lambda_{2}-c)(\lambda_1+\lambda_{2}-1)
(\lambda_1+\lambda_{2}+1)}
\end{align*}
belongs to $\CO_{i,j}^{\times}$.

\bigskip
\noindent$\bullet$ Case $(\lambda(i),\lambda(j))=(1/2,3/2)$ or $(3/2,1/2)$.

\smallskip
Assume that $(\lambda(i),\lambda(j))=(1/2,3/2)$.
 First note that
$\lambda(j)=0$ if and only if $\CHAR \cor=3$.

\medskip
Let us first assume $\CHAR \cor=3$. In this case,
$\lambda_1-1/2=y_1^2$, $y_{2}=(\lambda_{2}-\la(j))^2$, and
$$R_{i,j}=\tQ_{i,j}(y_1,y_2)=\pm(y_1^4-y_2)
=\pm\bl((\la_1-1/2)^2-(\la_2-3/2)^2\br)
=\pm(\la_1-\la_2+1)(\la_1+\la_2+1).$$ Hence we have \eqn
F(x_1,x_{2})^{-1}R_{i,j}&= &
\pm\frac{(\la_1-\la_2)^2(\la_1+\la_2)^2}%
{(\lambda_1-\lambda_{2}-1)(\lambda_1+\lambda_{2}-1)}
\eneqn
belongs to $\CO_{i,j}^\times$.

\medskip
Assume that $\CHAR \cor\ne 3$. In this case,
$\lambda_1=(y_1-1/2)^2$, $y_{2}=\lambda_{2}-3/2$, and $R_{i,j}=
\tQ_{i,j}(y_1,y_2)=\pm(y_1^2-y_2)=\pm(\la_1-\la_2+1)$.
Therefore
\eqn
F(x_1,x_{2})^{-1}R_{i,j}&= &
\pm\frac{(\la_1-\la_2)^2(\la_1+\la_2)^2}%
{(\lambda_1-\lambda_{2}-1)(\lambda_1+\lambda_{2}-1)(\lambda_1+\lambda_{2}-1)}
\eneqn
belongs to $\CO_{i,j}^\times$.

\bigskip
\noindent$\bullet$ Case $(\lambda(i),\,\lambda(j))=(0,\pm1)$ or $(\pm1,0)$.

\smallskip

We may assume that $(\lambda(i),\,\lambda(j))=(0,c)$ with $c=\pm1$.
 In this case,
$y_1=\lambda_1^2$, $y_{2}=\lambda_{2}-c$,
$$R_{i,j}=\tQ_{i,j}(y_1,y_2)=\pm(y_2^2-y_1)=\pm\bl((\la_2-c)^2-\la_1^2\br)
=\pm(\la_1-\la_2+c)(\la_1+\la_2-c).$$ Hence \eqn
F(x_1,x_{2})^{-1}R_{i,j}&= &
\pm\frac{(\la_1-\la_2)^2(\la_1+\la_2)^2}%
{(\lambda_1-\lambda_{2}-c)(\lambda_1+\lambda_{2}+c)}
\eneqn
belongs to $\CO_{i,j}^\times$.

\bigskip
\noindent
$\bullet$ Case $i=j$ and
$\lambda(i)\ne 0,1/2$.

In this case, $y_k=\la_k-\la(i)$ ($k=1,2$)
and $R_{i,j}=-(y_1-y_2)^{-2}=-(\la_1-\la_2)^{-2}$.
Then by \eqref{eq:Rf}
\eqn
F(x_1,x_2)^{-1}R_{i,i}&=&\dfrac{-(\la_1+\la_2)^2}%
{(\lambda_1-\lambda_{2}-1)(\lambda_1-\lambda_{2}+1)(\lambda_1+\lambda_{2}-1)(\lambda_1+\lambda_{2}+1)}
\eneqn
belongs to $\CO_{i,j}^\times$.
$\Bigl(\dfrac{x_1-x_2}{y_1-y_2}\Bigr)^2$.

\bigskip
\noindent$\bullet$ Case $i=j$ and $\lambda(i)=1/2$.

In this case, $y_k^2=\la_k-1/2$ ($k=1,2$), and
 $$R_{i,j}=-2(y_i^2+y_2^2)(y_1^2-y_2^2)^{-2}=-2(\la_1+\la_2-1)(\la_1-\la_2)^{-2}.$$
Hence by \eqref{eq:Rf}, we have \eqn
F(x_1,x_{2})^{-1}R_{i,j}&= &\frac{(\la_1+\la_2)^2}%
{(\lambda_1-\lambda_{2}-1)(\lambda_1-\lambda_{2}+1)(\lambda_1+\lambda_{2}+1)}
\eneqn
and it belongs to $\CO_{i,j}^\times$.

\bigskip
\noindent$\bullet$ Case $i=j$ and $\lambda(i)=0$.

\smallskip
In this case,  $y_k=\la_k^2$ ($k=1,2$),
and $R_{i,j}=-(y_1-y_2)^{-2}=-(\la_1^2-\la_2^2)^{-2}$.
 Hence
\eqn
F(x_1,x_{2})^{-1}R_{i,j}&= &
\dfrac{-1}{(\lambda_1-\lambda_{2}-1)(\lambda_1-\lambda_{2}+1)
(\lambda_1+\lambda_{2}-1)(\lambda_1+\lambda_{2}+1)}
\eneqn
belongs to $\CO_{i,j}^\times$.

\bigskip
This completes the proof of Theorem~\ref{th:isom_aDhc}. \qed

\bigskip
Similarly to the case of cyclotomic Hecke-Clifford superalgebras,
we can define the notion of cyclotomic Sergeev superalgebras.

For $\Lambda=\sum_{i\in I_1}m_i\Lambda_i$ with $m_i\in\Z_{\ge0}$
(see \S\;\ref{sec:cyclq}), we set
$$f'_\Lambda(X_1)=\prod_{i\in J}(x_1-x(i))^{m_{\pr(i)}},$$
and define
$$\dAHCC_n^\Lambda=\dAHCC_n/\bl(\dAHCC_n f'_\Lambda(X_1)\dAHCC_n\br).$$
We call $\dAHCC_n^\Lambda$
the {\em cyclotomic Sergeev superalgebra}.

Then, taking $I$ large enough comparing to
$\supp(\Lambda)\seteq\set{i\in I}{m_i\not=0}$,
the cyclotomic
Sergeev superalgebra $\dAHCC_n^\Lambda$ is isomorphic to the
cyclotomic quiver Hecke-Clifford superalgebra $\RC_n^\Lambda$.

Note that the Sergeev algebra $\hc$ is isomorphic to
$\dAHCC_n^{\Lambda_{i_0}}$, where $i_0\in  I$ is a unique odd vertex $0\in\AAA^1_w$.
\bigskip

For a generalized Cartan matrix $A$, 
we denote by $\KLRR{}{n}(A)$ the Khovanov-Lauda-Rouquier algebra associated with $A$.
It is nothing but $\R_n(A)$ with $\IODD=\emptyset$.
The algebra $\KLRR{}{n}(A)$ depends on the parameter $Q=(Q_{i,j})_{i,j\in I}$
as in \S\ref{specRC}.
By the rescaling in Remark~\ref{rem:scale},
we can easily see that the $\cor$-algebra $\KLRR{}{n}(A)$ is unique up to isomorphism
when $A$ is of finite type or affine type
except the case when the Dynkin diagram has a cycle, i.e., when $A$ is of
$A^{(1)}_{n-1}$ type
($n\ge2$).

\Rem\label{comp_calc1} Consider the case $\CHAR \cor=3$ and
$A=A^{(2)}_{2}$ (see the picture (ii) (b) in \S\ref{sec:Dyndahc}).
Take a block subsuperalgebra $B$ of the affine Sergeev superalgebra $\dAHC{11}$
which categorifies $U^-(\mathfrak{g}(A))_{-\beta}$ with
$\beta=8\alpha_0+3\alpha_1$, where $\alpha_0$ is the short root.
Although $\KLR{\beta}(A)$ also categorifies
$U^-(\mathfrak{g}(A))_{-\beta}$, the set of irreducible objects
$\Irr(\MOD{\KLR{\beta}(A)})$ in the category $\MOD{\KLR{\beta}(A)}$
of finite-dimensional $\KLR{\beta}(A)$-modules
and
the set of irreducible objects $\Irr(\SMOD{B})$ in the category $\SMOD{B}$
of finite-dimensional $B$-supermodules correspond to different perfect bases.

Let us explain in more detail. By ~\cite{BK3} (see also ~\cite[part
II]{Kle}), the cyclotomic Sergeev superalgebras
$\CMH{n}\simeq \RC_n^{\Lambda_0}(A)$ categorify
the irreducible highest module $V(\Lambda_0)$. Namely,
denoting by $\SMOD{\RC_n^{\Lambda_0}(A)}$ the supercategory
of finite-dimensional $\RC_n^{\Lambda_0}(A)$-supermodules, we have
\eq&&
\label{spinparam}
\ba{l}
\soplus_{n\geq
0}K^\super(\SMOD{\RC_n^{\Lambda_0}(A)})_{\mathbb{C}} \simeq V(\Lambda_0),\quad\\[1ex]
\bigsqcup_{n\geq 0}\Irr^\super(\SMOD{\RC_n^{\Lambda_0}(A)})\simeq B(\Lambda_0)
\simeq\RP{}
\ea\eneq
where $\Irr^\super(\SMOD{\RC_n^{\Lambda_0}(A)})$
is the set of equivalence classes of irreducible $\RC_n^{\Lambda_0}(A)$-
supermodules by the equivalence relation $\sim$
given by $S\sim S'$$\Leftrightarrow $ $S\simeq S'$ or $S\simeq\Pi (S')$.
The first isomorphism in \eqref{spinparam} is a $U(\mathfrak{g}(A))$-module
isomorphism
and the second isomorphism is a
$U_v(\mathfrak{g}(A))$-crystal isomorphism.
Recall that $\RP{}$ is
the set of all {\it 3-restricted 3-strict partitions}, which is
equivalent to the set of {\it reduced Young walls of type
$A_{2}^{(2)}$} in \cite{Kang03}. A partition
$\lambda=(\lambda_1,\ldots,\lambda_r)$ is 3-restricted 3-strict if
the following conditions are satisfied.
\begin{itemize}
\item $\lambda_k=\lambda_{k+1}$ implies $\lambda_k\in 3\Z$,
\item $\lambda_k-\lambda_{k+1}<3$ if $\lambda_k\in 3\Z$,
\item $\lambda_k-\lambda_{k+1}\leq 3$ if $\lambda_k\not\in 3\Z$.
\end{itemize}

Let $\MSP{n}$ denote
the spin symmetric group superalgebra of order $n$. Then
the cyclotomic Sergeev superalgebra $\dAHC{n}^{\Lambda_0}$
is isomorphic to the tensor product $\MSP{n}\otimes\clif_n$
of the spin symmetric group superalgebra $\MSP{n}$ and the Clifford superalgebra $\clif_n$
by
Sergeev and Yamaguchi~\cite{Ser2,Yam}.
Hence \eqref{spinparam} holds also for
$\MSP{n}$.

For each $\lambda\in \RP{}\cong B(\Lambda_0)$, we denote by
$\SV{\lambda}$ the  corresponding isomorphism class of irreducibles
of $\MSP{|\lambda|}$.

On the other hand, by ~\cite{KK,LV} we have
\begin{align}
\label{KLRparam} \bigoplus_{n\geq
0}\KKK(\MOD{\KLRR{\Lambda_0}{n}(A)})_{\mathbb{C}} \cong
V(\Lambda_0),\quad \bigsqcup_{n\geq
0}\Irr(\MOD{\KLRR{\Lambda_0}{n}(A)}) \cong B(\Lambda_0),
\end{align}
where the left isomorphism is as $U(\mathfrak{g}(A))$-modules and
the right isomorphism is as $U_v(\mathfrak{g}(A))$-crystals. For
each $\lambda\in \RP{}\cong B(\Lambda_0)$, we denote by
$\HV{\lambda}$ the  corresponding isomorphism class of irreducibles
of $\KLRR{\Lambda_0}{n}(A)$.

If both $\Irr(\MOD{\KLR{\beta}(A)})$ and $\Irr(\SMOD{B})$ correspond to
the same perfect basis on $U(\mathfrak{g}(A))$-module
$L(\Lambda_0)$, then we must have
\begin{align*}
{\dim \SV{\lambda}}/{\dim \HV{\lambda}} =
2^{[(1+\gamma_1(\lambda))/2]}
\end{align*}
for any $\lambda\in \RP{}$ (see ~\cite[Lemma 22.3.8]{Kle}).
Here, $\gamma_1(\lambda)=\sum_{i\in I_\even}m_i=n-\ell$ if
the weight of $\lambda$ is $\Lambda_0-\beta$ and $\beta=\sum_{i\in I}m_i\alpha_i$
($n=\hgt(\beta)=\sum_{i\in I}m_i$ and $\ell=\sum_{i\in I_\odd}m_i$).

A computer calculation shows that for $\lambda=(6,4,1)$ ($\gamma_1(\lambda)=3$),
we have $\dim \HV{\lambda}=648$ while it is known%
\footnote{Historically, it was first
miscalculated as $\dim\SV{\lambda}=2592=4\times 648$ in ~\cite{MY}. If it were
correct, observing such a direct discrepancy between
the Khovanov-Lauda-Rouquier
algebras and the spin symmetric groups must become more difficult.} that
$\dim\SV{\lambda}=2880=4\times 720$.
\enrem

\Rem\label{comp_calc2}
Consider the case (ii) (d) in  \S\ref{sec:Dynahc}, $A=D^{(2)}_2=A^{(1)}_1$
indexed by $I=I_\odd=\{0,1\}$ and assume $q=\exp(2\pi\sqrt{-1}/8)\in \cor$ with $\CHAR \cor=0$.
By the isomorphism described in ~\cite[\S3.2.1]{Rou1} or Remark~\ref{rem:scale}, it is enough to
consider $\KLRR{\Lambda_0}{4}(A)$ for $Q_{0,1}(u,v)=u^2-a\mspace{1.5mu}uv+v^2$ for some $a\in\cor$.

By ~\cite{BK1,BK2,KK,Rou1,Tsu},
the family of superalgebras
$\{\RC_n^{\Lambda_0}(A)\}_{n\ge0}\simeq\{\AHCC^{\Lambda_0}_{n}(q)\}_{n\geq 0}$
(resp.\ $\{\KLRR{\Lambda_0}{n}(A)\}_{n\geq 0}$)
categorifies $U(\mathfrak{g}(A))$-module
(resp.\ $U_v(\mathfrak{g}(A))$-module) $V(\Lambda_0)$.
However, there is no Morita equivalence between
 $\RC^{\Lambda_0}_{4}(A)$ and $\KLRR{\Lambda_0}{4}(A)$ nor
weak Morita superequivalence between $\RC^{\Lambda_0}_{4}(A)$ and $\KLRR{\Lambda_0}{4}(A)$
whatever superalgebra structure we give to $\KLRR{\Lambda_0}{4}(A)$ and
for any choice of parameters $a\in\cor$.

The algebras $\RC^{\Lambda_0}_{4}(A)$ and $\KLRR{\Lambda_0}{4}(A)$
are not Morita equivalent
because we have (for any $a\in\cor$)
\begin{align*}
\dim Z(\RC^{\Lambda_0}_{4}(A))=4
\ne 5
=\dim Z(\KLRR{\Lambda_0}{4}(A))
\end{align*}
by a computer calculation. (See Remark~\ref{rem:center} for $Z$.)
Since it can be easily seen that $\Irr(\MOD{\KLRR{\Lambda_0}{4}(A)})$ consists of 2 irreducible modules of
dimension $1,4$ for any $a\in\cor$, 
these two irreducibles
are self-associate (see \S\;\ref{sec:selfa})
for any superalgebra structure on $\KLRR{\Lambda_0}{4}(A)$.
Moreover, $\Irr(\MOD{\RC^{\Lambda_0}_{4}(A)})$ also consists of two self-associate
 irreducible
modules. Hence
two supercategories $\MOD{\RC^{\Lambda_0}_{4}(A)}$ and
$\bl(\MOD{\KLRR{\Lambda_0}{4}(A)}\br)^\ct$
cannot be superequivalent for any superalgebra structure on
$\KLRR{\Lambda_0}{4}(A)$ by Lemma~\ref{lem:IMQ}.
Since $\KLRR{\Lambda_0}{4}(A)$ has no non-trivial block decomposition,
$\RC^{\Lambda_0}_{4}(A)$ and $\KLRR{\Lambda_0}{4}(A)$ cannot be
weakly Morita superequivalent.

\enrem

\end{document}